\documentclass[a4paper,reqno]{amsart}
\usepackage[foot]{amsaddr}

\usepackage[dvipsnames]{xcolor}
\usepackage{latexsym}

\usepackage{stmaryrd}

\usepackage{soul}

\usepackage{amsfonts,amssymb} 
\usepackage{stmaryrd}
\usepackage{amsmath}
\usepackage{mathtools}

\usepackage{cancel} 

\newtheorem{theorem}{Theorem}[section]
\newtheorem*{theorem*}{Theorem}
\newtheorem{lemma}[theorem]{Lemma}
\newtheorem{corollary}[theorem]{Corollary}
\newtheorem{proposition}[theorem]{Proposition}
\newtheorem{remark}[theorem]{Remark}
\theoremstyle{definition}
\newtheorem{definition}[theorem]{Definition}
\newtheorem{example}[theorem]{Example}

\theoremstyle{remark}

\usepackage{xcolor}
\usepackage{xparse}
\usepackage{tikz}
\usetikzlibrary{matrix}
\usetikzlibrary{cd}

\definecolor{amber}{rgb}{1.0, 0.49, 0.0}
\definecolor{amethyst}{rgb}{0.45, 0.31, 0.59}
\definecolor{applegreen}{rgb}{0.55, 0.71, 0.0}
\definecolor{asparagus}{rgb}{0.0, 0.42, 0.24}
\definecolor{darkbyzantium}{rgb}{0.36, 0.22, 0.33}
\definecolor{darklavender}{rgb}{0.45, 0.31, 0.59}

\begin{document}

\def\mathunderline#1#2{\color{#1}\underline{{\color{black}#2}}\color{black}}

\newcommand{\abs}[1]{\lvert#1\rvert}
\newcommand{\ID}{\mathbb{I}}
\newcommand{\complessi}{\mathbb C}
\newcommand {\interi}{\mathbb Z}
\newcommand{\reali}{I\!\!R}
\newcommand{\razionali}{\mathbb Q}
\newcommand{\pd}{\partial_x}
\newcommand\Pa{Painlev\'e}
\newcommand\Pt{{\rm P}_{\rm{\scriptstyle I}}}
\newcommand\Lr{{\mathcal L}}
\newcommand{\nn}{\nonumber}
\newcommand\Pol{{\mathcal P}}
\newcommand\Se{{\mathcal S}}
\newcommand\A{{\mathcal A}}
\newcommand\M{{\mathsf M}}
\newcommand {\KZ}{Knizhnik--Zamolodchikov }
\newcommand {\KZt}{Knizhnik--Zamolodchikov--type }
\newcommand {\KKSS}{Kirillov-Kostant-Souriau }
\newcommand {\KKS}{standard Lie--Poisson }

\newcommand {\KKZ}{_{\scriptscriptstyle{\operatorname{KZ}}}}
\newcommand {\nablak}{\nabla\KKZ}
\def\cwedge{\,\overset{\wedge}{,}\,}
\def\d{{\rm d}}
\def\ddt{\frac{\rm d}{{\rm d} t}}
\def\tr{{\rm Tr\,}}
\def\res{{\rm res}}
\def\ddx{\frac{\rm d}{{\rm d} x}}
\def\ddtt{\frac{\rm d}{{\rm d}\tau}}
\def\ddz{\frac{\rm d}{{\rm d} z}}
\def\RE{{\rm Re}}
\def\IM{{\rm Im}}
 \def\arg{{\rm arg}}
\def\iff{\Longleftrightarrow}
\def\U{{\mathcal U}}
\def\pderiv{\partial_x}
\def\orexp#1{#1^{\raise3pt\hbox{$\scriptstyle\!\!\!\!\!\!\!\to$}}\,{}}
\def\one#1{#1^{\raise5pt\hbox{$\scriptstyle\!\!\!\!1$}}\,{}}
\def\two#1{#1^{\raise5pt\hbox{$\scriptstyle\!\!\!\!2$}}\,{}}
\def\onetwo#1{#1^{\raise4pt\hbox{$\scriptstyle\!\!\!\!\!{12}$}}\,{}}
\def\twoone#1{#1^{\raise4pt\hbox{$\scriptstyle\!\!\!\!\!{21}$}}\,{}}
\def\ocomma{{\phantom{\Bigm|}^{\phantom {X}}_{\raise-1.5pt\hbox{,}}\!\!\!\!\!\!\otimes}\,}
\def\resA#1#2#3#4{\prescript{[#4]}{}{#1}^{(#2)}_#3}
\def\resApt#1#2#3{\prescript{[#3]}{}{#1_#2}}

\title[Isomonodromic deformations: Confluence, Reduction and Quantisation]
{Isomonodromic deformations: Confluence, Reduction and Quantisation}

\author{Ilia Gaiur}   
\email[Ilia Gaiur]{iygayur@gmail.com}
\author{Marta Mazzocco}
\email{M.Mazzocco@bham.ac.uk}
\address[Ilia Gaiur,Marta Mazzocco]{School of Mathematics, University of Birmingham, United Kingdom}

\author{Vladimir Rubtsov}
\email{volodya@univ-angers.fr}
\address[Vladimir Rubtsov]{LAREMA, UMR 6093 du CNRS, Mathematics Department, University of Angers, France, Institute of Geometry and Physics (IGAP), Trieste, Italy.}

\begin{abstract}
In this paper we study the isomonodromic deformations of systems of differential equations with poles of any order on the Riemann sphere as Hamiltonian flows on the product of co-adjoint orbits of the truncated current algebra, also called generalised Takiff algebra. Our motivation is to produce confluent versions of the celebrated Knizhnik--Zamolodchikov  equations and explain how their  quasiclassical solution can be expressed via the isomonodromic $\tau$-function. 
In order to achieve this, we study the confluence cascade of $r+ 1$ simple poles to give rise to  a singularity  of arbitrary Poincar\'e rank $r$ as a Poisson morphism and explicitly compute the isomonodromic Hamiltonians.
\end{abstract}
\maketitle

\noindent {\bf Statements and declarations.} This research was supported by the Engineering and Physical Sciences Research Council Grant $EP/P021913$. V.R. was partly supported by the project IPaDEGAN (H2020-MSCA-RISE-2017), Grant Number 778010 and by the Project ``GraNum"- 80 Prime du CNRS. The authors did not receive support from any other organization for the submitted work. The authors have no relevant financial or non-financial interests to disclose.

\section*{Introduction}
In this paper we study the theory of isomonodromic deformations for systems of differential equations with poles of any order on the Riemann sphere. Our initial motivation was to generalise an observation by Reshetikhin that the quasi--classical solution of the standard  Knizhnik--Zamolodchikov
equations (i.e. with simple poles) is expressed via the isomonodromic $\tau$-function arising in the case of Fuchsian systems \cite{Res}. Along the way of pursuing the project of extending this to poles of any order, we have found a number of interesting results, some of which were already known as folklore (i.e. either done as very specific examples or not really proved in detail), others completely original. 

The Knizhnik--Zamolodchikov (KZ) equations emerged in theoretical physics as the system of linear differential equations satisfied by the correlation functions in the two--dimensional Wess--Zumino--Witten model of conformal field theory associated to a genus $0$ curve \cite{KZ, BPZ}. In the case of $\mathfrak g=\mathfrak{gl}_m$, the KZ equations can be represented as 
a system of linear differential equations for a local section $\psi$ of the trivial bundle $B\times U(\mathfrak{gl_m}(\mathbb C))^{\otimes n}\to B$ over the base $B$ given by the configuration space of ordered $n$-uples of points in $\mathbb C$, namely $B:=\{(u_1,\dots,u_n)\in\mathbb C^n| u_i \neq u_j  \hbox{ for } i\neq j\}$:
\begin{equation}\label{KZ} 
d\psi = \sum_{i\neq j} \Pi^{ij}\frac{du_i -du_j}{u_i - u_j}\psi,
\end{equation}
where $ \Pi^{ij}\in{\rm End} (U(\mathfrak{gl}_m(\mathbb C)^{\otimes n})$ is the extension of the  non-degenerated symmetric tensor
$$\Pi \in \mathfrak{gl}_m(\mathbb C)\times \mathfrak{gl}_m(\mathbb C) = {\rm End (\mathfrak{gl}_m(\mathbb C))}$$
acting by left multiplication on the $i-$th and $j-$th components of the tensor product $U(\mathfrak{gl}_m(\mathbb C)^{\otimes n})$ and trivially on the others.
Geometrically one can think about \eqref{KZ} as a flat Hitchin connection in geometric quantisation \cite{Hit}.

As proved by Reshetikhin in \cite{Res} (see also \cite{Har} where this result was explained in terms of passing from Shr\"odinger to Heisenberg representation), the KZ equations can be also viewed as deformation quantisation of the Schlesinger system \cite{Sch} of non-linear differential equations 
\begin{equation}\label{eq:schl}
d A^{(i)} = \sum_{i\neq j} [A^{(i)}, A^{(j)}] \frac{du_i -du_j}{u_i - u_j},
\end{equation}
controlling the isomonodromic deformation of a Fuchsian system on $\mathbb P^ 1$,
\begin{equation}\label{SchlesConn}
 \frac{dY}{d \lambda} = \sum_{i=1}^n\frac{A^{(i)}}{\lambda-u_i}Y,
 \end{equation}
 with $n+1$ simple poles $u_1,\dots, u_n,\infty$. These equations are multi--time non-autonomous Hamiltonian systems with
 Hamiltonians
 \begin{equation}\label{SchlesH}
H_i : B\times \mathfrak{gl}_m(\mathbb C)^{n}\to \mathbb C
\end{equation}
given by 
\[
H_i := \sum_{i\neq j} \frac{{\rm Tr}(A^{(i)}A^{(j)})}{u_i - u_j}.
\]
Interestingly, if we treat the quantities $u_1\dots,u_n$ in the Hamiltonian as parameters rather than times, these Hamiltonians
form a family of autonomous Poisson commuting Hamiltonians  called  {\it Gaudin Hamiltonians}. This simple observation has been key to several efforts to introduce specific examples of confluent analogues of KZ: by first introducing confluent analogues of Gaudin, then quantising them and finally generating the non-autonomous versions. Let us give a summary of our understanding of these results here below.

The main idea for the quantisation of the Gaudin Hamiltonians was based on the standard point of view that for any finite dimensional Lie algebra $\mathfrak g$, the universal enveloping algebra $U(\mathfrak g)$ can be considered as a deformation of the symmetric algebra $S(\mathfrak g)$ via the  Poincar\'e-Birkhoff-Witt map. One then defines the quantum enveloping algebra as
\[
U_{\hbar}(\mathfrak g) = T(\mathfrak g)/(X\otimes Y - Y\otimes X - \hbar [X,Y]), \quad X,Y\in \mathfrak g,
\]
by naturally extending the symmetrisation map  to the map $ S(\mathfrak g)^{\otimes n} \to U_{\hbar}(\mathfrak g)^{\otimes n}$,  and then the functions $ {\rm Tr}(A^{(i)} A^{(j)})$ on $\mathfrak g^{\otimes n}$ are transformed to $\Pi^{ij}.$ 

To define a quantisation of the Gaudin Hamiltonians  it is necessary to describe the Hilbert space of the quantum model as
tensor product of some representations of ${\mathfrak g}^{\oplus n}$.  The quantised Hamiltonians $\widehat H_i$ act on this Hilbert space and the quantum problem consists in finding their spectrum, matrix elements and so on. Formulated rigorously, the quantum Gaudin Hamiltonians generate a large commutative subalgebra in $U(\mathfrak g)^{\otimes n}$ which can be easily completed to a maximal commutative subalgebra. This subalgebra is usually called {\it Gaudin} or {\it Bethe subalgebra}. The explicit formulae for the generators (namely the quantum Hamiltonians) were obtained in \cite{MTV,Tal}.

In the case of $\mathfrak{g}=\mathfrak{gl_m}$, one can fix a co-vector $\mu \in \mathfrak{g}^{*}$ and using  the standard basis of $\mathfrak{gl_m}$ one can re-write the quantised Gaudin Hamiltonians as
\begin{equation}\label{Gaudimham_gl}
\widehat H_i = \sum_{j\neq i} \sum_{r,s =1}^m  \frac{E^{(i)}_{rs} E^{(j)}_{sr}} {u_i - u_j} + \sum_{r,s =1}^m  \mu(E_{rs})E^{(i)}_{sr},
\end{equation}
where $E^{(i)}_{rs}$ means $E_{rs}$ (as the element of standard basis in $\mathfrak{gl_m}$) considering in the $i-$th tensor factor. We observe that even the case of {\it regular} 
$\mu \in \mathfrak{g}^{*}$ (i.e. semi-simple, when $\mu (E_{rs}) = \mu_r \delta_{rs}$ with distinct $\mu_r \in \mathbb C$), the point  $\infty$ is an order two pole. The case of semi--simple but not regular $\mu$ was treated in \cite{FFR}.  

The autonomous Gaudin model  \eqref{Gaudimham_gl} can be generalised in two directions: by allowing higher order singularities at the marked points $u_i \in \mathbb C$ thus giving rise to Gaudin models with irregular singularities in \cite{FFTL} or by taking an element  $\mu \in \mathfrak{g}^{*}$ that  is not semi-simple (i.e. has non-trivial Jordan blocks).  These two approaches were unified in the classical and in the quantum cases  in \cite{VY} where an analogue of the {\it bispectral dynamical duality} of \cite{FMTV} between the models was proved.

 The next important step consisted in deforming the quantum Gaudin Hamiltonian to obtain KZ. This was done in the case of the $A_n$ root system by de Concini and Procesi \cite{dCP} and generalised to any Lie algebra in \cite{MTL,FMTV}. More precisely,
for any complex simple Lie algebra $\mathfrak{g}$ with a Cartan subalgebra $\mathfrak{h}\subset \mathfrak {g}$ and a corresponding root system $\Delta\subset \mathfrak{h}^{*}$, Millson and Toledano-Laredo \cite{MTL} introduced 
 the following {\it Casimir connection:}
 \begin{equation}\label{casconn}
 \nabla_C \psi:= d\psi - \frac{\hbar}{2\pi i}\sum_{\alpha\in \Delta} C_{\alpha} \frac{d\alpha}{\alpha}\psi, 
 \end{equation}
 where for every $\alpha$ one takes the principal embedding of $\mathfrak{sl}_2$ so that $C_{\alpha} =  \frac{\langle\alpha,\alpha\rangle}{2}(e_{\alpha}f_{\alpha} + f_{\alpha}e_{\alpha} + \frac{1}{2}h_{\alpha}^2)$ is the Casimir in 3-dimensional subalgebra $\mathfrak{sl}_{2,\alpha}$
 with respect to the restriction of the fixed non-degenerated $ad-$invariant bilinear form $\langle -,-\rangle$ on $\mathfrak{sl}_{2,\alpha}$ and $\hbar\in \mathbb C.$ {A special class of quantum connections with one irregular singularity of Poincar\'e rank $2$ and several other simple poles appeared in \cite{FMTV} as dual to the standard KZ connection, and in \cite{Boalch0} was re--obtained as quantisation of Dubrovin's system (without the skew-symmetry condition). Dubrovin system was then generalised to simply laced Dynkin diagrams in \cite{Boalch} and quantised in \cite{Rem}. }
 
Confluent versions of the KZ equation, or in other words, KZ equations with irregular singular points of  arbitrary Poincar\'e rank were obtained for $\mathfrak{sl}_2$ by Jimbo, Nagoya and Sun
 \cite{JNS}.  In \cite{FFTL} a class of quantum integrable systems generalising the Gaudin model was introduced by considering non-highest weight representations of any simple Lie algebra. These {\it Gaudin models with irregular singularities}\/ are expected to give rise to {\it confluent KZ equations}\/ as the corresponding differential equations on conformal blocks. Such KZ equations have not been explicitly written and one of the purposes of this paper is to do this.

In order to achieve our aim, we first needed to find explicit formulae for the isomonodromic Hamiltonians and to introduce a good set of Darboux coordinates. We have succeeded in doing this for a class of isomonodromic connections which behave well under confluence. Let us describe this class in some detail here.
It is well known that
the isomonodromic deformation equations in the case of higher order poles have a co-adjoint orbit interpretation on a  current Lie algebra. In the case of the Painlev\'e equations, Harnad and Routhier \cite{HR} produced finite--dimensional parameterisations by introducing suitable truncations of the current Lie algebra 
- in this paper we call such truncated current Lie algebras {\it Takiff algebras}\footnote{In fact, in this paper, we deal with the so-called {\it generalised} Takiff algebras, while the original definition was only for degree $1$ truncation.} for brevity (see Section \ref{se:takiff} for the definition). Korotkin and Samtleben \cite{KS} then conjectured the \KKS bracket on the Takiff algebras {and later Boalch proved that indeed these brackets are preserved by the Jimbo-Miwa isomonodromic deformations \cite{Boalch2}}. In this paper, we unify these two approaches to study connections as elements of the product of co-adjoint orbits in the Takiff algebra. More precisely, we consider linear systems of ODEs with poles at $u_1,u_2,\dots,u_n,\infty$ of Poincar\'e rank $r_1,r_2,\dots,r_n,r_{\infty}$ respectively, in the form
\begin{equation}\label{SchlesConnPR}
 \frac{dY}{d \lambda} =A(\lambda) Y, \quad A(\lambda) = \sum_{i=1}^n \sum_{k=0}^{r_i} \frac{A^{(i)}_k}{(\lambda-u_i)^{k+1}}+\sum_{k=1}^{r_\infty} A^{(\infty)}_k z^{k-1},
 \end{equation}
where  $A(\lambda)$ is an element of the phase space 
\begin{equation}\label{eq:ph-sp}
M  :=\hat{\mathcal{O}}_{r_1}^{\star}\times\hat{\mathcal{O}}_{r_2}^{\star}\times\dots \hat{\mathcal{O}}_{r_n}^{\star}\times \hat{\mathcal{O}}_{r_\infty}^{\star}, \end{equation}
where $\hat{\mathcal{O}}_{r_i}^{\star}$ stands for the co-adjoint orbit of the complex Lie group $\widehat G_{r_i}$ corresponding to the Takiff algebra of degree $r_i$, for $r_i>0$, and for the standard Lie algebra $\mathfrak g$ co-adjoint orbit  for $r_i=0$.

Following the ideology of \cite{AdHarPr}, in Theorem \ref{th:QPKKS}, we show how to obtain the \KKS bracket  
\begin{equation}\label{KKSgen}
\{A^{(i)}_k\ocomma A^{(j)}_l\} =  \left\{\begin{array}{ll}
    -\delta_{ij}[{\Pi}, A^{(i)}_{k+l}\otimes\ID] & k+l \leq  r_i \\
    0 & k+l > r_i.
    \end{array}\right.
\end{equation}
on our phase space \eqref{eq:ph-sp} 
as Marsden--Weinstein reduction of the Poisson structure on 
$$
 \oplus_{i=1}^{n+1} (T^{\star}\mathfrak{gl}_m)^{r_i+1}= \oplus_{k=1}^{d} T^{\star}\mathfrak{gl}_m,
 $$
 obtained by endowing each copy of $T^{\star}\mathfrak{gl}_m$ with the canonical symplectic structure $\d P\wedge \d Q$. Here $d=\sum_{i=1}^{n+1} r_i+ n+1$ denotes the degree of the divisor $D$ of the connection \eqref{SchlesConnPR}.
 The Marsden--Weinstein reduction is obtained by the additional first integrals given by the moment maps of the inner group action by $\widehat G_{r_i}$ as in formulae \eqref{eq:innerTakiff}.

These coordinates $(Q_1,P_1,\dots,Q_d,P_d)$, that we call {\it lifted Darboux coordinates,} were first introduced by Jimbo, Miwa, Mori and Sato  in the case of linear systems of ODEs with $n$  simple poles and possibly a Poincar\'e rank one pole at $\infty$ \cite{JMMS}. Harnad generalised these coordinates to allow rectangular $m_1\times m_2$ matrices and used them to generalise Dubrovin duality \cite{Dub1} between two systems of linear ODEs: one of dimension $m_1$ and the other of dimension $m_2$ \cite{Har1} and \cite{W2}. 
Similar coordinates were also introduced and partly used in the context of non-autonomous Hamiltonian description of Garnier-Painlev\'e differential systems  by M. Babich and S. Derkachov \cite{Bab,BabDer}. However in these latter works, the authors restricted to the case of rational parametrisation of co--adjoint orbits of $Gl_n(\mathbb C)$ and other semi-simple Lie groups and  did not consider current Lie algebras. 

Interestingly, using the lifted Darboux coordinates, we can describe all possible isomonodromic systems with a fixed degree $d$ of the divisor of the poles of the connection \eqref{SchlesConnPR} 
as  Marsden--Weinstein
reductions of different inner group actions on 
the universal phase space $ \oplus_{k=1}^{d} T^{\star}\mathfrak{gl}_m$. These reductions give rise to symplectic leaves of dimension $(r_1+\dots+r_n+r_\infty+n)(m^2-m)$. We explain how to produce the Darboux coordinates, which we call {\it intermediate Darboux coordinates,} on such symplectic leaves. In the case of the Jimbo-Miwa isomonodromic problems associated to the fifth, fourth, third and second Painlev\'e equations the degree is always $d=4$,  the intermediate symplectic leaves have always dimension $6$ and are determined by the choice of $3$ spectral invariants giving a total dimension $9$ for the Poisson manifold. This  is the dimension of the moduli space of $SL_2(\mathbb C)$ connections with a given divisor $D$ of degree $4$ \cite{Kri}.

\begin{remark}
The problem of extending the Riemann-Hilbert symplectomorphism between the de Rahm moduli space of meromorphic connections on a Riemann surface $\Sigma$ with non-simple divisor (a divisor of points that can have multiplicity $>1$) and the analogous of the Betti moduli space of representations of the fundamental group of $\Sigma$, namely with the cusped character variety introduced in \cite{CMR,CMR1} is still open and is beyond the scope of the current paper. However, the Darboux coordinate description of the de Rahm moduli space achieved in this paper constitues an important first step towards that goal.
\end{remark}

\begin{remark}
It is worth mentioning here that the phase space \eqref{eq:ph-sp} is not a moduli space per se, however K. Hiroe and D. Yamakawa \cite{HY} showed that the sub--space  of stable connections admits a nice quotient with respect to the diagonal action of $GL_m(\mathbb C)$ on $M$: 
\[
M' =\{A(\lambda)\in M\vert \sum_{i=1}^{n+1}\pi(A^{(i)}_0) = 0, \text{``stable"
}\}/GL_m(\mathbb C),
\]
where 
\[
\pi :\widehat{\mathfrak g}_{r_i}^{*} \to \mathfrak{gl}_m^{*}
\]
is the moment map under the diagonal action of $GL_m(\mathbb C)$ on $M$, thus assuring that $M'$ is a smooth complex symplectic variety. The space $M'$ can be regarded as a certain moduli space for meromorphic connections on $\mathcal{O}_{\mathbb P^1}^{\oplus m}.$
Fix $n$ distinct points $u_1,\ldots,u_n \in \mathbb P^1$, and endow $\mathbb P^1$ with a coordinate $z$ for which $z(u_i)\neq \infty.$ The variable $z_i$ can be identified with $\lambda-u_i$ and $\widehat{\mathfrak g}_{r_i}^{*} $ can be embedded in $\mathfrak{gl}_m(\mathbb C[z^{-1}_i])\frac{dz_i}{z_i}$ via trace-residue pairing. Then each $A(\lambda)\in M$ determines a meromorphic connection $d-A(\lambda)$ on on $\mathcal{O}_{\mathbb P^1}^{\oplus m}$, having poles at $u_1,\dots,u_n,\infty.$ The condition $ \sum_{i=1}^{n+1}\pi(A^{(i)}_0) = 0$ singles out the connections which have no residue at infinity. 
\end{remark}

Our next result is the characterisation of the outer linear automorphisms of the Takiff algebra that preserve the \KKS structure \eqref{KKSgen} on the phase space \eqref{eq:ph-sp} (see Theorem \ref{th:aut0Tak} for a more articulated statement). 

\begin{theorem}\label{th:aut0Takint} Consider two elements $A(\lambda)$ and $B(\lambda)$ of the phase space \eqref{eq:ph-sp}, so that they both have the form \eqref{SchlesConnPR}:
$$
A(\lambda) = \sum_{i=1}^n \sum_{k=0}^{r_i} \frac{A^{(i)}_k}{(\lambda-u_i)^{k+1}}+\sum_{k=1}^{r_\infty} A^{(\infty)}_k \lambda^{k-1},\quad
B(\lambda) = \sum_{i=1}^n \sum_{k=0}^{r_i} \frac{B^{(i)}_k}{(\lambda-u_i)^{k+1}}+\sum_{k=1}^{r_\infty} B^{(\infty)}_k \lambda^{k-1}.
$$
Assume that $A(\lambda)$ and $B(\lambda)$ are related by a linear automorphism of the Lie algebra
$$
B^{(i)}_k =\sum_{l=0}^{r_i} T^{(i)}_{kl} A^{(i)}_l, \quad \forall i=1,\dots,n, \infty, 
$$
for some scalar quantities $T^{(i)}_{kl}$.  Then the Poisson condition
$$
 \{B^{(i)}_k\ocomma B^{(j)}_l\} =  \left\{\begin{array}{ll}
    \delta_{ij}[B^{(i)}_{k+l}\otimes\ID,{\Pi}] & k+l \leq  r_i \\
    0 & k+l > r_i
    \end{array}\right.\,\iff\,   \{A^{(i)}_k\ocomma A^{(j)}_l\} =  \left\{\begin{array}{ll}
    \delta_{ij}[ A^{(i)}_{k+l}\otimes\ID,{\Pi}],  & k+l \leq  r_i \\
    0 & k+l > r_i
    \end{array}\right.
$$
is satisfied if and only if the the following formulae are satisfied:
\begin{equation}\label{sect3:monom0}
B_k^{(i)}= \sum\limits_{j=k}^{r_i} A_j^{(i)}\mathcal{M}^{(r_i)}_{{k},j}(t_{1}^{(i)},t_2^{(i)},\dots t_{r_i}^{(i)}),
\end{equation}
where
\begin{equation}\label{sect3:monom}
\mathcal{M}^{(r_i)}_{k,j} = \sum\limits_{w(\alpha)=j}^{|\alpha|=k}\frac{k!}{\alpha_1!\alpha_2!\dots\alpha_{r_i}!}
\left(\prod\limits_{l=1}^{r_i}(t_l^{(i)})^{\alpha_l}\right),\quad |\alpha| = \sum\limits_{l=1}^{r_i} \alpha_l,\quad w(\alpha) = \sum\limits_{l=1}^{r_i} l \cdot \alpha_l.
\end{equation}
\end{theorem}

This  result allows us to introduce extra (i.e. in addition to the positions of poles) deformation parameters $t^{(i)}_1,\dots,t^{(i)}_{r_i}$, $i=1,\dots,n,\infty$ for any connection belonging to the phase space \eqref{eq:ph-sp}. In other words, we consider families of the form
$$
A(\lambda) = \sum_{i=1}^n \sum_{k=0}^{r_i} \frac{B^{(i)}_k}{(\lambda-u_i)^{k+1}}+\sum_{k=1}^{r_\infty} B^{(\infty)}_k \lambda^{k-1}
$$ 
where the elements $B^{(i)}_k$ contain explicitly
the deformation parameters $t^{(i)}_1,\dots,t^{(i)}_{r_i}$ as prescribed by formulae \eqref{sect3:monom0} and \eqref{sect3:monom}. The isomonodromic deformation equations will then impose a further implicit dependence of the matrices $A_k^{(i)}$ on the deformation parameters $t^{(i)}_1,\dots,t^{(i)}_{r_i}$ and on the position of the poles $u_1,\dots,u_n$.

\begin{remark}
    The set of parameters $t_1,\dots,t_r$ introduced in Theorem \ref{th:aut0Tak} may be replaced by the coefficients of the jet-expansion of local conformal changes of coordinates $z=z(\zeta)$ at the pole\footnote{We thank the referee for pointing out this nice geometric analogy.}. Such a description is natural in the framework of  irregular isomonodromic problems and provides a nice geometric intuition. However, the formulas for the coefficients of the jet-expansions are more complicated, especially for the study of the corresponding Hamiltonians.  In example \ref{example-auto-conf} we calculate the relation between the coefficients of such local conformal map and the parameters  parameters $t_1,\dots,t_r$ for $r=3$; this can be easily generalised to any $r$.
\end{remark}

\begin{remark} Let us stress that the class of connections we consider in this paper are
 elements of the space  (\ref{eq:ph-sp}). This class excludes some of the Jimbo-Miwa-Ueno connections. Indeed, our deformation parameters correspond to a subset of the Jimbo-Miwa-Ueno ones and this correspondence is $1:1$ only in the case of rank $m=2$. 
 For 
 example, the famous Dubrovin's system
$$
\frac{d Y}{d\lambda}= \left(U+\frac{V}{\lambda}\right) Y,
$$
where $U$ is a diagonal $n\times n$ matrix and $V\in \mathfrak{so}_n$,
is not an element of $\hat{\mathcal{O}}_{r_1}^{\star}\times\hat{\mathcal{O}}_{r_\infty}^{\star}$ for some $r_1,r_\infty$ because the diagonal elements of $U$ are independent deformation parameters.  Of course the isomonodromic deformation equations for $V$ as a function of $u_1,\dots,u_n$ can be written as a flow on a co--adjoint orbit 
$\mathcal{O}^{\star}$ of the Lie algebra $\mathfrak{so}_n$, but not as equations for the whole connection $U+\frac{V}{z}$ on the product of two co--adjoint orbits as our theory dictates. To include the Dubrovin system (and indeed all of Jimbo-Miwa-Ueno  deformation parameters) in our theory, one should either consider the extended coadjoint orbits introduced in 
\cite{Boalch1,Boalch2} or exploit the Laplace transform. In the latter setting, the confluence procedure distroys semi-simplicity, therefore it is a different process from the one considered by Cotti, Dubrovin and Guzzetti \cite{CGD, CGD1}.
\end{remark}

This is  the correct framework to study confluence of two or more poles. Indeed, we show that 
 the confluence cascade of $r+1$ simple poles at certain positions depending on $t^{(i)}_1,\dots, t^{(i)}_{r_i}$ gives rise to
 an element of the phase space \eqref{eq:ph-sp} which has a singularity of Poincar\'e rank $r$ and depends on $t^{(i)}_1,\dots,t^{(i)}_{r_i}$, $i=1,\dots,n,\infty$, as prescribed by formulae \eqref{sect3:monom0} and \eqref{sect3:monom}. The following theorem provides the inductive step to create the confluence cascade (we drop the index ${}^{(i)}$ for convenience).

\begin{theorem}\label{sect4:small_theorem}
Consider an $r$--parameter family of connections of the following form:
\begin{equation}\label{4.3_Conn}
    A = \sum\limits_{k=0}^{r}\frac{B_k(t_1,t_2\dots t_{r-1})}{(\lambda-u)^{k+1}} + \frac{C}{\lambda-v}+\hbox{holomorphic terms},
\end{equation}
where by holomorphic terms we mean terms holomorphic in $\lambda-u$ and $\lambda-v$, and each $B_k$ depends on the parameters $t_1,\dots,t_r$ as specified by \eqref{sect3:monom0}, \eqref{sect3:monom}. Assume
\begin{equation}\label{4.3_expan}
v = u + \sum\limits_{i=1}^r t_i \varepsilon^i = u + P_r(t,\varepsilon), 
\end{equation}
and that we have the following asymptotic expansions as $\varepsilon\to 0$
\begin{equation}\label{4.3_expanR}
     C \sim \sum\limits_{j=-r}^{\infty}W^{[j]}\varepsilon^j,\quad A_k\sim -\sum\limits_{l=1}^{r-k}\frac{W^{[-k-l]}}{\varepsilon^l}+A^{[k,0]}+\sum\limits_{l=1}^{\infty}A^{[k,l]}\varepsilon^l,
\end{equation}
for some matrices $W^{[-k-l]},A^{[k,l]}$.
Then the limit $\varepsilon\rightarrow 0$ the connection exists and is equal to
$$
 \tilde  A = \sum\limits_{i=0}^{r+1}\frac{\tilde {B}_i(t_1,t_2\dots t_r,t_{r+1})}{(\lambda-u)^{i+1}}+\hbox{holomorphic terms},
$$
where $\tilde{B}_i$'s are given by
\begin{equation}\label{sect3:IPA}
\tilde{B}_i(t_1\dots ,t_{r+1}) = \sum\limits_{k=i}^{r}\tilde A_k \mathcal{M}_{i,k}^{(r+1)}(t_1\dots t_{r+1}),\quad 
\tilde A_k=\left\{\begin{array}{ll}
    W^{[-k]}+A^{[k,0]}, & k<r+1. \\
    W^{[-r-1]}, & k=r+1.
\end{array}\right.
\end{equation}
\end{theorem}

We prove that the confluence procedure gives a Poisson morphism on the product of co-adjoint orbits and we calculate explicitly the confluent Hamiltonians, which define the correct isomonodromic deformations.

\begin{theorem}\label{th:mainHam}
Let $u$ be a pole of a connection $A$ with Poincar\'e rank $r$, which is the  result of confluence of $r$ simple poles with the simple pole $u$. Then the confluent Hamiltonians $H_1,\dots,H_{r}$ which correspond to the times $t_1,\dots t_r$ are defined as follows:
\begin{equation}\label{eq:ham-spec}
    \left(\begin{array}{c}
    H_1 \\
    H_2 \\
    \dots \\
    H_r
\end{array}\right) = \left(\mathcal{M}^{(r)}\right)^{-1}\left(\begin{array}{c}
    S^{(u)}_1 \\
    S^{(u)}_2 \\
    \dots \\
    S^{(u)}_r
\end{array}\right),
\end{equation}
where 
\begin{equation}
S_k^{(u)}= \frac{1}{2}\oint\limits_{\Gamma_u}(\lambda-u)^k\tr A^2 \d \lambda
\end{equation}
are spectral invariants of order $i$ in $u$ 
and the matrix $\mathcal{M}^{(r)}$ has entries $\mathcal{M}^{(r)}_{k,j}$ given by  (\ref{sect3:monom}). The Hamiltonian $H_u$ corresponding to the time $u$ is instead given by the standard formula
$$
{H}_{u_i} = \frac{1}{2}\underset{\lambda=u_i}{\res}\tr{A}(\lambda)^2.
$$
\end{theorem}

\begin{remark}
It is well known that the isomonodromic deformation equations are Hamiltonian,  namely that the flow is
Hamiltonian 
with respect to the Jimbo-Miwa-Ueno deformation parameters, see for example  \cite{Fed, Hurt, W1}. In \cite{Fed}, the isomonodromy equations have been described as integrable non-autonomous Hamiltonian systems. A symplectic fibre bundle whose base is the Jimbo-Miwa-Ueno deformation parameters space and the fibers are certain moduli spaces of unramified meromorphic connections was introduced in {\cite{Boalch2}}. This approach was extended by D. Yamakawa for any reductive Lie algebra $\mathfrak g$ \cite{Yam1} who removed some multiplicity restrictions and introduced a symplectic two-form on the fibration. Following the same geometric approach and Jimbo-Miwa-Ueno isomonodromic tau-function D. Yamakawa \cite{Yam2} has proven that the isomonodromy equations of Jimbo?Miwa?Ueno is a completely integrable non-autonomous Hamiltonian systems. He was also motivated by the quantisation theorem of Reshetikhin but he did not try to consider the quantisation of general isomonodromy equations\footnote{We are grateful to Prof. M.Jimbo who has drawn our attention and has sent a file of the paper \cite{Yam2}}. Recently, Bertola and Korotkin have derived a new Hamiltonian formulation of the Schlesinger equations (i.e. for the Fuchsian case) in terms of the dynamical $r$--matrix structure.
\end{remark}

\begin{remark}
The results of the theorem \ref{sect4:small_theorem} still hold true for the autonomous systems which are obtained by the confluence procedure from the Gaudin system. It was shown by Yu. Chernyakov in \cite{YuChern1} that the Poisson algebra which arises in the confluent elliptic and rational Gaudin systems coincides with the dual Takiff algebra equipped with the \KKS bracket (in \cite{YuChern1} the author use the word``fusion" instead of ``confluence").
\end{remark}

One of the main theorems of our paper gives a general formula for the confluent KZ Hamiltonians with singularities of arbitrary Poincar\'e rank in any dimension. 

\begin{theorem}\label{th:mainKZ}
Consider the differential operators: 
\begin{equation}\label{eq:KZopu}
\nabla_{u_j}:=
\frac{\partial}{\partial u_j}-\widehat H_{u_j},\quad j=1,\dots,n
\end{equation}
and
\begin{equation}\label{eq:KZopt}
\nabla_{k}^{(i)}:=\frac{\partial}{\partial t^{(i)}_k}-\widehat H_{k}^{(i)},\quad i=1,\dots,n,\infty, \quad k=1,\dots,r_i
\end{equation}
where the Hamiltonians $\widehat H_{u_j}$ which correspond to the positions of the poles $u_j$, $j=1\dots,n$, and $\widehat H_1^{(i)},\dots,\widehat H_{r}^{(i)}$ which correspond to the times $t_1^{(i)},\dots t_{r_i}^{(i)}$, for $i=1,\dots,n,\infty$, 
are given by the following elements of the universal enveloping algebra  $U\left(\hat{\mathfrak{g}}_{r_1}\oplus\dots  \oplus \hat{\mathfrak{g}}_{r_\infty}\right)$:
$$
\hat{H}_{u_j} = \frac{1}{2}\underset{\lambda=u_j}{\res}\tr_0 \hat{A}(\lambda)^2,
$$
and
$$
\mathcal{M}^{(r_i)}\left(\begin{array}{c}
    \hat{H}^{(i)}_1 \\
    \hat{H}^{(i)}_2 \\
    \dots \\
    \hat{H}^{(i)}_{r_i}
\end{array}\right) = \left(\begin{array}{c}
    \hat{S}^{(u_i)}_1 \\
    \hat{S}^{(u_i)}_2 \\
    \dots \\
    \hat{S}^{(u_i)}_{r_i}
\end{array}\right),\quad \hat{S}_k^{(u_i)} = \frac{1}{2}\oint\limits_{\Gamma_{u_i}}(\lambda-u_i)^k\tr_0 \hat{A}(\lambda)^2 \d \lambda,
$$
where
$$
\hat{A}(\lambda) = \sum\limits_{i}^n\left(\sum\limits_{j=0}^{r_i}\frac{\hat{B}^{(i)}_j\left(t^{(i)}_1,t^{(i)}_2\dots t^{(i)}_{r_i}\right)}{(\lambda-u_i)^{j+1}}\right),\
$$
with $\hat{B}^{(i)}$'s given by
$$
\hat{B}^{(i)}_j(t^{(i)}_1,\dots t^{(i)}_{r_i}) = \sum\limits_{k=j}^r \hat{A}^{(i)}_k\mathcal{M}^{(r_i)}_{j,k}(t^{(i)}_1,t^{(i)}_2\dots t^{(i)}_{r_i}),\quad \hat{A}_k = \sum\limits_{\alpha}e^{(0)}_{\alpha}\otimes e_{\alpha}^{(i)}\otimes z_i^k,
$$
and $e_\alpha^{(0)}$ corresponds to the quantisation of $\mathfrak g^*$ to $\mathfrak g$ while
$$
e_{\alpha}^{(i)} = 1  \otimes\dots\otimes \underset{i}{e_{\alpha}}\otimes\dots\otimes 1.
$$
Then the differential operators commute
$$
[\nabla_{u_j},\nabla_{u_s}]=[\nabla_{k}^{(i)},\nabla_{u_s}]=[\nabla^{(i)}_{k},\nabla_l^{(a)}]=0,
$$
$\forall\, j,s=1,\dots,n, \, i,a=1,\dots,n,\infty, \, k=1,\dots,r_i,\, l=1,\dots, r_a$.
We call the system of differential equations 
$$
\nabla_{u_j}\Psi=0,\quad \nabla_{k}^{(i)}\Psi=0,\qquad j=1,\dots,n, \, i=1,\dots,n,\infty, \, k=0,\dots,r_i,
$$
confluent KZ equations.
\end{theorem}

Moreover, we express the isomonodromic Hamiltonians in terms of the lifted Darboux coordinates and  show that the quasiclassical solutions of the confluent  KZ equations is expressed via the isomonodromic $\tau$-function.

\begin{theorem}\label{th-semicl}
Given a solution $(P_1,\dots,P_d,Q_1,\dots,Q_d)$ of the classical isomonodromic deformation equations, the corresponding semi-classical solution $\Psi_{sc}$ (see subsection \ref{suse:semi-cl}) of the confluent KZ equations
\begin{equation*}
\hbar \frac{\partial\Psi}{\partial u_j}=\widehat H_{u_j}\Psi,\quad j=1,\dots,n
\end{equation*}
and
\begin{equation*}
\hbar  \frac{\partial\Psi}{\partial t^{(i)}_k}=\widehat H_{k}^{(i)}\Psi,\quad i=1,\dots,n,\infty, \quad k=1,\dots,r_i
\end{equation*}
evaluated along the solution $(P_1,\dots,P_d,Q_1,\dots,Q_d)$, admits the following WKB expansion
\begin{equation}\label{eq:asyWBK}
    \Psi_{sc}(Q(t), t) \sim \tau^{\frac{i}{\hbar}}\left(1+O(\hbar)\right), \quad \hbar\to 0.
\end{equation}
in terms of the classical isomonodromic $\tau$-function
$$
\d \ln (\tau) := \sum\limits_{i}\left(H^{(i)}_{u_i}\d u_i+\sum_{k=1}^{r_i}H^{(i)}_{k}\d t_k^{(i)} \right).
$$
The asymptotic expansion \eqref{eq:asyWBK} is valid for $u_1,\dots,u_n$, $t_k^{(i)}$, $i=1,\dots,n,\infty$, $k=1,\dots,r_i$ in a poly-disk that does not contain the zeroes of the action functional evaluated along the given solution $(P_1,\dots,P_d,Q_1,\dots,Q_d)$.
\end{theorem}

This statement was mentioned in \cite{Res} for the case of the standard KZ, namely with simple poles. We also discuss the quantisation of the reduced Darboux coordinates and provide the quantised reduced systems in some examples.

\begin{remark}
Most of our results extend to the case of isomonodromic deformations for meromorphic connections on principal $G-$bundles over the Riemann sphere for any arbitrary complex reductive group $G$ - this is for example the situation of the famous Fuji-Suzuki $D_{2n+2}^{(1)}-$higher Painlev\'e hierarchies and matrix Painlev\'e equations \cite{BCR}. Only the results about the so-called ``lifted Darboux coordinates coordinates" can't immediately be generalised to any connected complex reductive group case. However, we have decided to restrict to the $GL_m(\mathbb C)$ case having in mind a wider audience. For the same reasons, we often do not use the language of sheaves and schemes. To extend our results to higher genus Riemann surfaces is instead a rather serious job. First, one can extend the ``rational" truncated polynomial currents to their trigonometric and elliptic analogues and to define a proper pairing and basis. For $g=1$ case such job can be done probably using the results of \cite{PRS1, PRS2, EPR} and we postpone to subsequent papers.
\end{remark}

This paper is organised as follows. In Section 1, we recall the case of Fuchisan connections, we discuss the lifted Darboux coordinates and the Marsden--Weinstein reduction to the phase space \eqref{eq:ph-sp} in the case of $r_1=\dots=r_n=r_\infty=0$ and remind the Hamiltonian formulation. In Section 2, we collect facts about the Takiff algebras;  we discuss the lifted Darboux coordinates at each separate pole for any choice of the Poincar\'e rank $r$ and show 
how to obtain the \KKS bracket  
\eqref{KKSgen}
as Marsden--Weinstein reduction of the Poisson structure on $\sum_{k=0}^{r}T^{\star}\mathfrak{gl}_m$. We also 
prove Theorem \ref{th:aut0Takint} and discuss the inner group action on the universal phase space.
 Finally, we show how to obtain the intermediate Darboux coordinates and discuss some examples.
In Section 3, we discuss the isomonodromic deformations. In Section 4, we discuss the confluence procedure. We first carry out the confluence of two simple poles,  explain how to obtain confluence cascades, prove Theorems \ref{sect4:small_theorem} and \ref{th:mainHam}. Is Section 5, we apply the theory to the case of the Painlev\'e equations. In Section 6, we deal with quantisation. We give a general formula for the confluent KZ equations with singularities of arbitrary Poincar\'e rank and prove Theorem \ref{th-semicl}.\\

We conclude this introduction with a discussion about further research directions emerging from our work. In the case of Fuchsian systems, the analytic continuation of the solutions of the isomonodromic equation is described by the Artin braid group $B_n$ realised as the fundamental group of the configuration space $B$ of $n$ points quotiented by the natural action of the symmetric group $S_n$ \cite{DM}. At quantum level, due to T. Kohno and V. Drinfeld, the universal $R$--matrix of $U_{\hbar}\mathfrak g$ gives a representation of $B_n$ as the monodromy of the KZ equation with values in $V^{\otimes n}$, where $V$ is a finite $U_{\hbar}\mathfrak g$--module. This is based on the fact that KZ is realised as a $S_n$-equivariant flat connection on the topologically trivial vector bundle over $X_n$ with fibre $V^{\otimes n}$. It would be interesting to understand how to modify this picture under confluence. In particular, it is not yet clear what happens if one braids two punctures in the Fuchsian system, say $u_1$ and $u_2$, and then confluence $u_2$ with a third puncture $u_3$.

Even more interesting is the problem to confluence the duality between KZ and the Casimir connection by Millson and Toledano-Laredo and to study its effect on the monodromy. Indeed,
given a Lie algebra $\mathfrak g$  with Cartan sub-algebra $\mathfrak h$ and Weyl group $W$, this Casimir connection is a $W$-equivariant flat  connection on $\mathfrak h$ with simple poles along the root hyperplanes and values in any finite-dimensional $\mathfrak g$-module $V$. 
In the case of $\mathfrak g=\mathfrak{sl}_n$, Toledano-Laredo \cite{TL} proved that the monodromy of this family of connections is equivalent to the quantum Weyl group action of the generalised braid group $B_{\mathfrak g}$ of type $\mathfrak g$ on $V$ obtained by regarding the latter as a module over the quantum group $U_{\hbar}\mathfrak g$. While the $R$--matrix representation is a deformation of the natural action of the symmetric group $S_n$ on $V^{\otimes n}$, the representation of $B_{\mathfrak g}$ deforms the action of a finite extension of $W$  on any finite-dimensional ${\mathfrak g}$--module. In the same paper \cite{TL} Toledano-Laredo showed that the duality between $\mathfrak{sl}_{k}$ and $\mathfrak{sl}_l$ derived from their joint action on the space of $k\times l$ matrices exchanges KZ for $\mathfrak{sl}_{k}$ with the Casimir connection for $\mathfrak{sl}_l$. The current paper opens the problem  of confluencing this duality.\\

\noindent{\bf Acknowledgements}
The authors are grateful to M. Babich, M. Bershtein, M. Bertola, Y. Chernyakov, H. Desiraju, B. Dubrovin, P. Gavrilenko, J. Harnad, J. Hurtubise, M. Jimbo, D. Korotkin, I. Sechin, B. Malgrange, N. Nikolaev, V. Poberezhny, G. Rembado, D. Talalaev, V. Tarasov, V. Toledano--Laredo, A. Veselov and A. Zotov for interesting discussions  {and to P. Boalch for useful bibliographic remarks}. We also thank both referees for their careful reading and their useful critical remarks which definitely contributed to the improvement of our paper.

\section{Fuchsian systems}
The aim of this section is to review the Poisson and symplectic aspects of the deformation equations for connections over the $n+1$-holed sphere with simple poles at the punctures. Starting from the linear system with simple poles at $\lambda=u_1,\dots,u_n,\infty$, 
\begin{equation}\label{intro:fuchs}
    \frac{d}{{d\lambda}} \Psi =\sum\limits_{i=1}^{n}\frac{A^{(i)}}{\lambda - u_i}\Psi,\quad \lambda\in \Sigma,\quad u_i\neq u_j,
\end{equation}
where $A^{(1)},\dots, A^{(n)}$ are non-resonant elements in $\mathfrak{sl}_m(\mathbb C)$ such that $A^{(\infty)}:= -\sum A^{(i)}\neq 0$,
we consider the following one-form
\begin{equation}\label{sect1:def}
    \Omega = (\d_{u} \Psi)\Psi^{-1},\quad \d_{u} \Psi := \sum\limits_{i}\partial_{u_i}\Psi \d u_i.
\end{equation}
Since we consider only isomonodromic deformations, $\Omega$ is a one-form valued meromorphic function in the variable $\lambda$ with simple poles at $u_1\dots,u_n,\infty$.

Using the local solutions of (\ref{intro:fuchs}) in the neighbourhood of the poles $u_i$'s and applying Liouville theorem, this form may be written as
\begin{equation}\label{sect1:defForm}
    \Omega = - \sum\limits_{i} \frac{A^{(i)}}{\lambda - u_i} \d u_i.
\end{equation}
The compatibility condition for (\ref{intro:fuchs}) and (\ref{sect1:def}), also called zero-curvature condition,
\begin{equation}
\d_{u} A - \frac{d}{{d\lambda}}\Omega + [A, \Omega] = 0,    
\end{equation}
gives the Schlesinger equations \eqref{eq:schl}. 

\subsection{Phase space}

The Schlesinger equations are Hamiltonian, with natural phase space  given by the direct product of co-adjoint orbits which are symplectic leaves of the \KKS bracket:
$$
\left(A^{(1)},A^{(2)},\dots A^{(n)}\right)\in \mathcal{O}^{\star}_1\times\mathcal{O}^{\star}_2\times\dots \times \mathcal{O}^{\star}_n.
$$
In the case when $\mathfrak{g}$ is a Lie algebra  with a non-degenerate bi-linear form (i.e. Killing form), we may identify the co-adjoint orbits with the adjoint orbits. The Poisson brackets  may then  be written as
\begin{equation}\label{sect1:FuchsPB}
\left\{A^{(i)}\ocomma A^{(j)}\right\} = \delta_{ij}[{\Pi},1\otimes A^{(i)}] \quad \Longleftrightarrow \quad \left\{A^{(i)}_{\alpha}, A^{(j)}_{\beta}\right\}= -\delta_{ij}\sum\limits_{\gamma} \chi^{\gamma}_{\alpha\beta} A^{(i)}_{\gamma}
\end{equation}
where the lower indices $\alpha, \beta$ and $\gamma$ correspond to the Lie co-algebra basis, $\chi_{\alpha\beta}^\gamma$ are the structure constants of the Lie algebra and ${\Pi}$ is a quadratic Casimir element. In the case of $\mathfrak{gl}_m$ it acts as a permutation operator, i.e. 
$$
{\Pi} (A\otimes B) {\Pi}=B\otimes A.
$$
In the case of a Lie algebra with an orthogonal basis $e_{\alpha}$ with respect to the Killing form, the quadratic Casimir $\Pi$ writes as
$$
\Pi = \sum_{\alpha}e_{\alpha}\otimes e_{\alpha}.
$$
Such bracket may be rewritten as an $r$-matrix bracket for the connection, i.e.
\begin{equation}\label{sect1:r-br}
    \{A(\lambda)\ocomma A(\mu)\} = \left[\frac{{\Pi}}{\lambda-\mu}, A(z)\otimes \ID + \ID \otimes A(\mu)\right].
\end{equation}
The isomonodromic Hamiltonians for the Schlesinger equations are
\begin{equation}\label{sect1:SchlesingerHam}
    H_i = \underset{\lambda=u_i}{\operatorname{Res}}\operatorname{Tr} \frac{A(\lambda)^2}{2} = \sum\limits_{j\neq i}\frac{\operatorname{Tr}(A^{(i)}A^{(j)})}{u_i-u_j}.
\end{equation}
In the case of $3$ co-adjoint orbits in $\mathfrak{sl}_2$, the Schlesinger equations can be reduced to the
 Painlev\'e VI equation which is a non-autonomous Hamiltonian system with 1 degree of freedom.
 
For $n$ co-adjoint orbits, the fully reduced dimension can be computed using the spectral type technique introduced by Katz \cite{Katz}. When all matrices $A^{(i)}$, $i=1,\dots,n$ are semi-simple, 
the spectral type approach gives the dimension of the fully reduced phase space as a function of the eigenvalues multiplicities of the residues, including the residue at infinity given by the Fuchs condition:
\begin{equation}\label{sect1:FuchRel}
\sum\limits_{i=1}^n A^{(i)} = - A^{(\infty)}.
\end{equation} Katz' formula is:
\begin{equation}\label{sect1:KatzForm}
N = 2 - (1- {n}) {m^2} -\sum\limits_{i=1}^{ {n}}\sum\limits_{j=1}^{ {l_i}} {(m^i_j)^2}-\sum\limits_{j=1}^{ {l_\infty}} {(m^\infty_j)^2},
\end{equation}
where $l_i$ is the cardinality of the set of eigenvalues for the residue $A^{(i)}$ and $m_j^i$ is the multiplicity of the $j$-th eigenvalue of the residue $A^{(i)}$ for $i=1,\dots,n,\infty$.

The Fuchs condition \eqref{sect1:FuchRel} may be viewed as the moment map of the Hamiltonian group action of conjugation by $z$-independent invertible matrices and $A_\infty$ is a constant of motion for the Schlesinger equations.

Formula \eqref{sect1:KatzForm} coincides with the dimension of the reduced space under symplectic reduction as follows:
\begin{equation}\label{sect1:KatzSymp}
N = \sum\limits_{i=1}^n \operatorname{dim}\mathcal{O}_i^{\star} - \operatorname{dim}G - \operatorname{stab}\mathcal{O}_{\infty}^{\star},
\end{equation}
 where $\operatorname{stab}\mathcal{O}_{\infty}^{\star}$ is the dimension of the stabilizer for the Jordan form of the residue at $\infty$. When $A^{(\infty)}$ is an element of the co-adjoint orbit of generic form ($ad$-regular), we have that $\operatorname{stab}\mathcal{O}_{\infty}^{\star} = \operatorname{dim}\mathfrak{h}$, so  formula \eqref{sect1:KatzSymp} simplifies to
$$
N = \sum\limits_{i=1}^n \operatorname{dim}\mathcal{O}_i^{\star} - \operatorname{dim}G - \operatorname{dim}\mathfrak{h}.
$$
For example, in the case of the Painlev\'e VI equation, we deal with the coadjoint orbits of the $\mathfrak{sl}_2(\mathbb{C})$ and formula (\ref{sect1:KatzSymp}) gives
$$
N = 3\cdot \operatorname{dim}\mathcal{O}_{\mathfrak{sl}_2} - \operatorname{dim}(SL_2) -\operatorname{dim}\mathfrak{h}_{\mathfrak{sl}_2} = 3\cdot 2 - 3 - 1 = 2,
$$
which is exactly the dimension of the Painlev\'e VI equation phase space. In some sense the multiplicity of the eigenvalues tells us that the Jordan form may be written as the tensor product of identity matrices of sizes corresponding to the 
the multiplicities. The stabilizer of such matrix is the set of  block diagonal matrices, so the dimension is greater then the dimension of the Cartan torus and finally we obtain the smaller phase space.\\

Our first goal is to describe this full reduction as a Hamiltonian reduction and a Marsden-Weinstein quotient. To this aim, we will need first to extend the phase space to $T^{\star}\mathfrak{gl}_m$ and show that the Darboux coordinates on this cotangent bundle reduce to the \KKSS form on the co-adjoint orbits. We will then discuss how the invariants of the co-adjoint orbits correspond to moment maps with respect to different Hamiltonian group actions on the extended phase space.

\subsection{Extended phase space and its Darboux coordinates}
In this subsection, we start by working locally, namely we restrict to the case of a single co-adjoint orbit $\mathcal{O}^{\star}$ of $\mathfrak{gl}_m$ and identify $\mathfrak{gl}_m^{\star}$ with $\mathfrak{gl}_m$ via the Killing form. In the last part of this subsection we extend to the product of $n$ co-adjoint orbits.

We consider $T^{\star}\mathfrak{gl}_m$ with the standard Darboux coordinates $(Q,P)$ and the canonical symplectic structure:
\begin{equation}\label{eq1:Darboux}
\omega = \tr\left(\d P \wedge \d Q\right) = \sum\limits_{i,j} \d P_{ij}\wedge \d Q_{ji}.
\end{equation}

Following \cite{AdHarPr,AHR,AHR1}, we explain how to obtain the \KKS bracket \eqref{sect1:FuchsPB} on $\mathfrak{g}^{\star}$ as Marsden--Weinstein reduction of the Poisson structure on $T^{\star}\mathfrak{gl}_m$. There is a direct way to see this reduction by a straightforward computation (see  \cite{JMMS}) that we resume in the next proposition:

\begin{proposition}\label{sect1:prop1} Consider the canonical symplectic structure on $T^{\star}\mathfrak{gl}_m$:
\begin{equation}\label{sect1:Darboux}
\omega = \tr\left(\d P \wedge \d Q\right) = \sum\limits_{i,j} \d P_{ij}\wedge \d Q_{ji}.
\end{equation}
Let 
\begin{equation}\label{sect1:DarbouxA}
A=QP,
\end{equation}
where we use the ring structure of $\mathfrak{gl}_m$ to justify the multiplication of  $Q$ and $P$. Then $A$ satisfies the
\KKS bracket \eqref{sect1:FuchsPB} for $\mathfrak{gl}_m$. 
\end{proposition}

\begin{proof}
The Poisson bracket which corresponds to the symplectic form in (\ref{sect1:Darboux}) may be written in the following way
$$
\{P\ocomma Q\} = {\Pi},\quad \{P \ocomma P\}=\{Q \ocomma Q\} = 0.
$$
Inserting this relation to the bracket between $A$'s we obtain
\begin{multline*}
\{A\ocomma A\} = \{QP\ocomma QP\} = (Q\otimes\ID) \{P\ocomma Q\} (\ID \otimes P) + 
(\ID \otimes Q) \{Q\ocomma P\} (P\otimes \ID) = \\ = (Q\otimes\ID) {\Pi} (\ID \otimes P) -
(\ID \otimes Q) {\Pi} (P\otimes \ID) = [{\Pi}, \ID \otimes QP] =[{\Pi}, \ID \otimes A]
\end{multline*}
As we wanted to prove.
\end{proof}

\begin{definition}
We call $T^{\star}\mathfrak{gl}_m$ {\em extended phase space} and the canonical coordinates $P,Q$  {\em lifted Darboux coordinates}.\end{definition}

To restrict  to the co-adjoint orbit, we have to fix the invariants of the co--adjoint action, i.e. the Jordan form of matrix $QP=A$. Such a procedure leads to some additional non-linear equations for the entries of $Q$ and $P$, and there is no hope to derive the explicit symplectic structure on the co-adjoint orbit from such a perspective. Therefore, we follow the construction of \cite{AdHarPr}  to obtain the co-adjoint orbits via Hamiltonian reduction.

The space $T^{\star}\mathfrak{gl}_m \simeq \mathfrak{gl}_m\times \mathfrak{gl}_m$ carries two natural commuting symplectic actions of $GL_m$ which we call {\em inner and outer}:
\begin{equation}
g\underset{\operatorname{inner}}{\times}(P,Q) = (gP,Qg^{-1}),\qquad h\underset{\operatorname{outer}}{\times}(P,Q) = (Ph,h^{-1}Q),\qquad h,g\in GL_m.
\end{equation}

\begin{lemma}\label{lm:inoutga}
the inner and outer actions are Hamiltonian with equivariant moment maps given by
\begin{equation}
    \begin{array}{ll}\mu_{\operatorname{inner}}:&T^{\star}\mathfrak{gl}_m\to \mathfrak{gl}_m^\star \\
& (P,Q)\mapsto \Lambda = PQ\end{array}\quad \begin{array}{ll}\mu_{\operatorname{outer}}:&T^{\star}\mathfrak{gl}_m\to \mathfrak{gl}_m^\star \\
& (P,Q)\mapsto A = QP\end{array}.
\end{equation}
\end{lemma}

Let us restrict to the open affine subset of $T^{\star}\mathfrak{gl}_m$ where at least one of the two matrices $Q$ and $P$ is invertible. For example $Q$. Then, resolving the moment map for $\Lambda$ we obtain
$$
P = \Lambda Q^{-1},\quad A = QP=Q\Lambda Q^{-1}.
$$
As a consequence, $A$ and $\Lambda$ belong to the same co-adjoint orbit.

Since the inner and outer actions commute, $A$ is invariant under the inner action, while $\Lambda$ is invariant under the outer action. Therefore we use the inner group action to fix $\Lambda$  in  Jordan normal form  without changing $A$. In other words, we take  the Jordan normal form $\Lambda_0$ of $A$ and select $\Lambda=\Lambda_0$. This gives
$$
T^{\star}\mathfrak{gl}_m \underset{\Lambda_0}{\sslash} G  = \mu_{\operatorname{inner}}^{-1}(\Lambda_0) \slash G ,
$$
here we denote by $ \underset{\Lambda_0}{\sslash}$ the quotient with respect to the inner action of  $GL_m$  on $T^{\star}\mathfrak{gl}_m$. 
We may resume these results in the following:

\begin{lemma}\label{sect1:lemma1}
The map
$$
\begin{array}{ccc}
T^{\star}\mathfrak{gl}_m \underset{\Lambda_0}{\sslash} G_{\operatorname{inner}}& \to &\mathcal{O}^{\star}\\
(Q,P)&\mapsto&  A:=QP \\
\end{array}
$$
is a rational symplectomorphism and the  Jordan normal form $\Lambda_0$ of $A$ is given by
$$
\Lambda_0=PQ.
$$
\end{lemma}

\begin{remark}
When $A$ is a full-rank matrix, both $P$ and $Q$ must be invertible. So we may embed $(P,Q)$ into the group $GL_m$ and $P$ and $Q$ can be seen as
left and right eigenvector matrices for the matrix $A$. In the case when $A$ may be diagonilized, the action of the Cartan torus (i.e. the stabilizer of $\Lambda$) leads to a well known fact from linear algebra - the eigenvectors are defined up to multiplication by non-zero constant.
When $A$ is not a full-rank matrix, we may choose $Q$ to be an invertible matrix (so it may be viewed as an element of $GL_n$). Then the rank of $P$ must equal to the rank of $A$. The the moment map $\Lambda$ will inherit the rank of $A$ automatically. Since $P$ in this case not invertible, the reduced coordinates take the form
$$
P = \Lambda Q^{-1},\quad A=Q\Lambda Q^{-1},\quad \operatorname{det}\Lambda = \operatorname{det} A = \operatorname{det} P =0.
$$
This means that instead of considering $T^{\star}\mathfrak{gl}_m$ as lifted space, we could take $T^{\star}GL_m \ni(Q,\Lambda Q^{-1}) $. Such consideration is closely related to the approach introduced in \cite{BerKor}. However, this approach is not very useful for our purposes, since we wish to work with polynomial unreduced parametrisation, rather then rational.
\end{remark}

\begin{remark}
In the case when we consider $\mathfrak g$ to be any reductive Lie algebra and $A\in \mathfrak g^\ast$, then we expect that Lemma \ref{sect1:lemma1} is still valid if we fix the value $\Lambda$ of the moment map in $ \mathfrak g^\ast$ and $Q$ and $P$ (or just $Q$ in the case of degenerate orbit) as the elements from the corresponding Lie group $G$.
\end{remark}

Let us now consider the case of the product of many co-adjoint orbits. Since the Poisson brackets  \eqref{sect1:FuchsPB} are {\it local,}\/ namely the residues at different  points commute, the facts we summarised so far easily extend to this case. Indeed, we can apply the above construction to the co-adjoint orbit  at each pole of the Fuchsian system (except $\infty$) and define:
$$
A^{(i)} = Q_iP_i.
$$
In this case we have that inner and outer actions can be lifted to the direct sum of n copies $T^{\star}\mathfrak{gl}_m$ in a natural way 
$$
g\underset{\operatorname{inner}}{\times}(P_1,P_2,\dots P_n,Q_1,Q_2,\dots Q_n) = (g_1P_1,\dots g_n P_n,Q_1g_1^{-1},\dots Q_ig_i^{-1},\dots Q_n g_n^{-1}),\quad g\in \underset{n}{\times} GL_m
$$
$$
h\underset{\operatorname{outer}}{\times}(P_1,P_2,\dots P_n,Q_1,Q_2,\dots Q_n) = (P_1h_1,\dots  P_nh_n,h_1^{-1}Q_1,\dots h_i^{-1}Q_i,\dots h_n^{-1}Q_n ),\quad h\in \underset{n}{\times} GL_m
$$
and the lemma \ref{lm:inoutga} is easily generalised as follows:
\begin{lemma}\label{lm:inoutgaGlob}
These inner and outer actions are Hamiltonian with equivariant moment maps given by
$$
    \mu_{\operatorname{inner}}:\begin{array}{rcl}\underset{n}{\oplus}T^{\star}\mathfrak{gl}_m & \to & \underset{n}{\oplus}\mathfrak{gl}_m^\star \\
(P_1,\dots P_n;Q_1,\dots Q_n) & \to & {(P_1Q_1,P_2Q_2,\dots P_nQ_n)}\end{array}
$$
$$
\mu_{\operatorname{outer}}:\begin{array}{rcl}\underset{n}{\oplus}T^{\star}\mathfrak{gl}_m & \to &  \underset{n}{\oplus}\mathfrak{gl}_m^\star \\
(P_1,\dots P_n;Q_1,\dots Q_n) & \to & {(Q_1P_1,Q_2P_2,\dots Q_nP_n)}\end{array}
$$

\end{lemma}
\proof
Let us prove it for the inner action only. The vector field generated by the group action (via element $\xi = (\xi_1, \xi_2,...\xi_n) \in\oplus_n \mathfrak{gl}_m$
is given by
$$
X_\xi (P_i, Q_i) = \frac{d}{dt}(e^{-t\xi_i}P_i, Q_ie^{t\xi_i})\Big\vert_{t=0} = (-\xi_iP_i, Q_i\xi_i) = 
\sum\limits_{i=1}^n\Big(\sum\limits_{k,j} -(\xi_iP_i)_{kj}\frac{\partial\,}{\partial P_{i_{kj}}}+(Q_i\xi_i)_{kj}\frac{\partial\,}{\partial Q_{i_{kj}}}\Big).
$$
Inserting $X_\xi$ into the symplectic form we obtain
$$
\omega(X_\xi, \circ) = \sum\limits_i^n\sum\limits_{k,j}\Big[ -(\xi_iP_i)_{kj}dQ_{i_{jk}} - (Q_i\xi_i)_{kj}dP_{i_{jk}} \Big] = - \sum\limits_i^n\operatorname{Tr}\left(\xi_iP_idQ_i+Q_i\xi_idP_i\right) =- \sum\limits_i^n d\operatorname{Tr}\left(\xi_iP_iQ_i\right),
$$
so the corresponding Hamiltonian is
$$
h_\xi\left(m \right) = \langle \mu\left( m\right), \xi \rangle = \sum\limits_i^n\operatorname{Tr}\left(\xi_i\mu\left(m\right)_i\right) = \operatorname{Tr}\left(\xi_iP_iQ_i\right)
$$
where $m=(P_1,P_2,\dots P_n,Q_1,\dots Q_n)$. So the moment map is given by
$$
\mu\left(m\right) = \left(P_1Q_1,P_2Q_2,\dots P_iQ_i,\dots P_nQ_n\right),
$$
which is equivariant
$$
\mu\left(g\circ m\right) =  \left(g_1^{-1}P_1Q_1g_1,g_2^{-1}P_2Q_2g_2,\dots g_i^{-1}P_iQ_ig_i,\dots g_n^{-1}P_nQ_ng_n\right) = g^{-1}\mu(m)g = \operatorname{Ad}^{\star}_{g^{-1}}(\mu(m))
$$
\endproof

Then the following result is a straightforward computation

\begin{lemma}
A Hamiltonian system on the phase space 
$$
\mathcal{O}^{\star}_1\times\mathcal{O}^{\star}_2\times\dots \times \mathcal{O}^{\star}_n\ni \left(A^{(1)},A^{(2)},\dots A^{(n)}\right)
$$
can be lifted to the extended phase space
$$
 T^{\star}\mathfrak{gl}_m\times T^{\star}\mathfrak{gl}_m \times\dots \times T^{\star}\mathfrak{gl}_m\ni \left(Q_1,P_1,Q_2,P_2\dots Q_n,P_n\right)
$$
with additional first integrals given by the moment maps of the inner group action
$$
\mu_{\operatorname{inner}} := P_iQ_i = \Lambda^{(i)},
$$
where the inner group action is given by
$$
(g_1,g_2,\dots g_n) \underset{\operatorname{inner}}{\times}(P_1,Q_1,P_2,Q_2,\dots P_n,Q_n) = (g_1P_1,Q_1g_1^{-1},\dots g_iP_i,Q_ig_i^{-1},\dots g_n P_n,Q_n g_n^{-1}).
$$
Moreover, if $\Lambda^{(i)}_0$ is the Jordan normal form of $A^{(i)}$, we can fix   $\Lambda^{(i)}=\Lambda^{(i)}_0$.
\end{lemma}

In particular, the Schlesinger Hamiltonians (\ref{sect1:SchlesingerHam}) can be lifted to the extended phase space $T^{\star}\mathfrak{gl}_m$ as follows
\begin{equation}\label{eq:LiftedSchles}
    H_i=\sum\limits_{j\neq i}\frac{\operatorname{Tr}(Q_iP_iQ_jP_j)}{u_i-u_j},
\end{equation}
and it can be checked directly that they Poisson commute with the moment maps of the inner group action.

\subsection{Outer group action and the gauge group}
We have seen that the inner group action allows us to restrict from $T^{\star}\mathfrak{gl}_m$ to $\mathcal{O}^{\star}_1\times\mathcal{O}^{\star}_2\times\dots \times \mathcal{O}^{\star}_n$. Now we consider the {\it outer group action } that will allow us to reduce further. This is given by 
$$
(g_1,g_2,\dots g_n) \underset{\operatorname{outer}}{\times}(P_1,Q_1,P_2,Q_2,\dots P_n,Q_n) = (P_1g_1,g_1^{-1} Q_1,\dots P_ig_i,g_i^{-1}Q_i,\dots  P_n g_n, g_n^{-1}Q_n )
$$
and is also Hamiltonian (see Lemma \ref{lm:inoutga}).

Because inner and outer group actions commute,  their moment maps Poisson commute too. However, the Schlesinger Hamiltonians are generally not invariant under outer action, unless the outer action is restricted to be a diagonal action, i.e.
$$
g_1=g_2=\dots=g_n=g.
$$
In this case, the outer action reduces to the standard $GL_m$-action on $\mathcal{O}^{\star}_1\times\mathcal{O}^{\star}_2\times\dots \times \mathcal{O}^{\star}_n$, or equivalently to the constant gauge group action:
$$
g  \underset{\operatorname{outer}}{\times} A = \sum\limits \frac{g^{-1}A^{(i)}g}{z-u_i}.
$$
The moment map of such diagonal action is
\begin{equation}\label{sect1:GaugeMoment}
\sum\limits_{i=1}^nQ_iP_i=\sum\limits_{i=1}^n A^{(i)} = -A^{(\infty)},
\end{equation}
which is the Fuchs relation.

In order to describe the reduction procedure induced by the outer diagonal action in terms of the Marsden-Weinstein reduction, following  Proposition 2.2.7 of \cite{Audin} (see also  \cite{Hitchin}) we further extend the phase space by adding another copy of $T^{\star}\mathfrak{gl}_m$:
\begin{equation}
    (P_1,Q_1 \dots P_n,Q_n;P_{\infty},Q_{\infty})\in \bigoplus\limits_{i=1}^{n+1}T^{\star}\mathfrak{gl}_m,\quad \omega = \sum\limits_{i=1}^{n}\tr \d P_i\wedge \d Q_i + \tr \d P_{\infty}\wedge \d Q_{\infty},
\end{equation}
with the outer group action of the form
$$
g \underset{\operatorname{outer}}{\times} (P_1,Q_1 \dots P_n,Q_n;P_{\infty},Q_{\infty}) = (P_1g,g^{-1} Q_1,\dots P_ig,g^{-1}Q_i,\dots  P_n g, g^{-1}Q_n; P_{\infty} g, g^{-1}Q_{\infty} ).
$$
The corresponding extended space given by the reduction with respect to the inner group action takes form 
$$
\left(A^{(1)},A^{(2)},\dots A^{(n)}; A^{(\infty)}\right)\in \mathcal{O}^{\star}_1\times\mathcal{O}^{\star}_2\times\dots \mathcal{O}^{\star}_n\times\mathcal{O}_{\infty}^{\star}.
$$
The reduction with respect to the relation (\ref{sect1:GaugeMoment}) on the extended phase space may be viewed as the Marsden-Weinstein quotient
$$
\bigoplus\limits_{i=1}^{n+1}T^{\star}\mathfrak{gl}_m \sslash G = \mu^{-1}(0)/G,\quad \mu = \sum\limits_{i=1}^n Q_iP_i + Q_{\infty}P_{\infty},
$$
that corresponds to the Fuchsian relation on the phase space reduced with respect to the inner group action.

Finally, the fully reduced phase space then has form 
$$
M:=\mathcal{O}^{\star}_1\times\mathcal{O}^{\star}_2\times\dots \mathcal{O}^{\star}_n\times\mathcal{O}_{\infty}^{\star} \sslash G\simeq \bigoplus\limits_{i=1}^{n+1}\left(T^{\star}\mathfrak{gl}_m \underset{\Lambda^{(i)}}{\sslash} G  \right)\sslash G,
$$
Where $\simeq$  denotes the symplectomorphism between symplectic manifolds.\\

Moreover, the Hamiltonians are homogeneous polynomials in the lifted Darboux coordinates. Such dependence plays crucial role in the quantisation of the isomonodromic systems as we will discuss in  section \ref{se:quant}.\\

In this paper we extend this scheme for the isomonodromic problems with irregular singularities and will introduce a well defined confluence procedure that creates an irregular singularity of Poincar\'e rank $r$ as a result of collision of $r+1$ simple poles. In the next section, we study the case of the irregular singularities  along the same lines of the regular one.

\section{Takiff algebras and associated symplectic manifolds.}\label{se:takiff}

It is well known that
the isomonodromic deformation equations in the case of higher order poles also have a co-adjoint orbit interpretation on a current algebra. In the case of the Painlev\'e equations, Harnad and Routhier \cite{HR} produced finite dimensional parameterisations that can be interpreted as  introducing suitable truncations of the current algebra. Korotkin and Samtleben \cite{KS} then conjectured the \KKS structure on truncated current algebras - also called Takiff algebras. In this section, we unify these two approaches and classify the linear Takiff algebra automorphisms that preserve the \KKS structure. As a consequence, we obtain a general formula that prescribes the way to introduce independent deformation parameters in generic connections with poles of any Poincar\'e rank.

Loosely speaking, the Takiff algebra of degree $r_i$, is the Taylor part of a current algebra quotiented by the ideal generated by $z^{r_i}$ where $r_i$ is the order  of the pole at $u_i$ and $z$ is the local coordinate at $u_i$.   
For a general system with poles at $u_1,u_2,\dots,u_n,\infty$ of Poincar\'e rank $r_1,r_2,\dots,r_n,r_{\infty}$ respectively, the phase space is
$$
M :=\hat{\mathcal{O}}_{r_1}^{\star}\times\hat{\mathcal{O}}_{r_2}^{\star}\times\dots \hat{\mathcal{O}}_{r_n}^{\star}\times \hat{\mathcal{O}}_{r_\infty}^{\star} \sslash G_{\operatorname{gauge}},
$$
where $\hat{\mathcal{O}}_{r_i}^{\star}$ stands for the co-adjoint orbit of the Takiff algebra of degree $r_i$.

In this section we remind generalities about Takiff algebras and describe the Poisson structure on their co-adjoint orbits. Moreover, we  explain the lifted Darboux parametrisation for the co-adjoint orbits of Takiff algebras.  We show that the lifted space is always the same and the way to distinguish between different isomonodromic systems is the Hamiltonian group action we choose to obtain the reduced phase space. In Section \ref{se:CONF} we will show that that the  Takiff algebras algebras naturally arise during the confluence procedure.

The Takiff algebra $\hat{\mathfrak{g}}_{{r}}$ of the Lie algebra $\mathfrak{g}$ is the Lie algebra of polynomials of given degree $r$ in an indeterminate variable  $z$ with the following Lie bracket
\begin{equation}
    \left[\sum_{i=0}^{{r}} A_i z^i,\sum_{j=0}^{{r}} B_j z^j\right] = \sum_{i=0}^{{r}}
\left(\sum_{j=0}^i[A_i,B_{i-j}]\right) z^i.
\end{equation}
This algebra may be viewed as a double quotient of the loop algebra $\mathfrak{g}[[z]]$ as follows. Denote by $\mathfrak{g}[z]^{+}$ the subalgebra of the elements which has a finite limit when $z$ goes to the origin. Then $\hat{\mathfrak{g}}_{{r}}$ is defined as
$$
\hat{\mathfrak{g}}_{{r}} =\mathfrak{g}[z]^{+}\slash z^{{r+1}}\mathfrak{g}[z]^{+},\quad \mathfrak{g}[z]^{+}= \mathfrak{g}[[z]] \Big\slash \mathfrak{g}[z]^{-},\quad  \mathfrak{g}[z]^{-} = \left\{f\in \mathfrak{g}[[z]]:\,\lim\limits_{z\rightarrow \infty}f(z) = 0\right\},
$$
These algebras are known in the Integrable Systems community as   truncated loop algebras or truncated current
algebras. The variable $z$ is usually called spectral parameter and, as we will illustrate here below, it induces a grading on the Takiff algebra.

In the case when $\mathfrak{g}$ admits an invariant non-degenerate bi-linear form (Killing form), we may define the co-algebra $\hat{\mathfrak{g}}_r^{\star}$ in the following way
$$
\hat{\mathfrak{g}}_r^{\star} = \mathfrak{g}[z]^{-}/ z^{-(r+1)-1} \mathfrak{g}[z]^{-} = \left\{ A=\frac{A_{r}}{z^{r+1}}\dots+\frac{A_0}{z} \,\Big\vert \, A_i\in \mathfrak{g}\right\}.
$$
The pairing between $\hat{\mathfrak{g}}_r$ and $\hat{\mathfrak{g}}_r^{\star}$ is given by
\begin{equation}
\label{eq:pairing}
\langle A, B\rangle = \oint\limits_{S^1}\tr{(AB)}dz = \sum\limits_{i=0}^{{r}}\tr{A_iB_i}.
\end{equation}
Let us assume that the Lie algebra $\mathfrak{g}$ is given by
$$
\mathfrak{g}=\operatorname{Span}\left\{X_1,\dots X_{m}\right\},\quad [X_i,X_j] = C_{ij}^{k}X_k, \quad \langle X_i, X_j \rangle = \delta_{ij},
$$
then for the Takiff algebra $\hat{\mathfrak{g}}_{{r}}$ we have the following basis and structure equations
$$
X_{\alpha,i}=X_iz^{\alpha},\quad [X_{i,\alpha},X_{j,\beta}] = \left\{\begin{array}{ll}
    C_{ij}^k X_{k,\alpha+\beta}, & \alpha+\beta \leq r  \\
    0 & \alpha+\beta > r .
\end{array}\right.
$$
For the dual algebra $\mathfrak{g}_{{r}}^{\star}$, we use the following basis
$$
X^{\alpha,i} = X^{i}z^{-\alpha{-1}}, \quad \langle X^{i},X_j\rangle = \delta_{ij},
$$
so that the pairing is given by
$$
\langle X^{i,\alpha}, X_{j,\beta}\rangle = \delta_{\alpha \beta} \langle X^{i},X_j\rangle = \delta_{\alpha \beta}\delta_{ij}.
$$

The details about Takiff algebras or truncated current algebras and their \KKS bracket may be found in \cite{FadTak} (see part 2, chap. 4 \S1). In the following sub-section we recall the essentials of this construction.

\subsection{Standard Lie--Poisson bracket for the Takiff algebras} 
Let us remind the reader that the \KKS bracket on the dual Lie algebra $\mathfrak{g}^{\star}$ is given by
$$
\{f,g\}(L) = -\langle L, [\d f(L),\d g(L)]\rangle,\quad f,g \in C^{\infty}(\mathfrak{g}^{\star}),\quad \d f(L),\d g(L) \in \mathfrak{g}.
$$
The coadjoint orbits $\mathcal{O}^{\star}$ are symplectic leaves of the \KKS structure on $\mathfrak{g}^{\star}$. The vector fields on $\mathcal{O}^{\star}$ may be identified with the elements of Lie algebra $\mathfrak{g}$ and the symplectic form takes the form
$$
\omega_{\operatorname{KKS}}(X,Y)(L) =  -\langle L, [X,Y]\rangle.
$$
Following \cite{FadTak}, we now describe the \KKS structure on the dual $\hat{\mathfrak{g}}_n^{\star}$ of the Takiff algebra. Let's consider the following element of the dual $\hat{\mathfrak{g}}_n^{\star}$
$$
A = \sum\limits_{\alpha=1}^{{r}}\sum_{i} A_{\alpha, i}X^{\alpha,i}\in \hat{\mathfrak{g}}_n^{\star},
$$

The coefficients $A_{\alpha,i}$ are functions on the coadjoint orbit, with $\d A_{\alpha,i}=X_{i,\alpha}$ so that the \KKS bracket is given by
\begin{equation}\label{eq:TakKKS}
\{A_{\alpha, i},A_{\beta, j}\} = - \langle A , [X_{i,\alpha},X_{j,\beta}] \rangle = -\langle A , C_{ij}^{k}X_kz^{\alpha+\beta} \rangle = \left\{\begin{array}{ll}
    -C_{ij}^{k}A_{\alpha+\beta, k}, & \alpha+\beta \leq r  \\
    0 & \alpha+\beta < r. 
\end{array}\right.
\end{equation}
This is a graded Poisson structure of degree $1$, and the Takiff co-algebra inherits the grading:
$$
\hat{\mathfrak{g}}_{{r}}^{\star}: = \bigoplus\limits_{i=0}^{r}\hat{\mathfrak{g}}^{\star, i}_{{r}},\quad \{\hat{\mathfrak{g}}^{\star,i}_{{r}},\hat{\mathfrak{g}}^{\star,j}_{{r}}\}\subseteq \hat{\mathfrak{g}}^{\star,{i+j}}_{{r}},
$$
where $\hat{\mathfrak{g}}^{\star,i}_{{r}}= \left\{ A=\frac{A_{{i}}}{z^{{i+1}}} \,\Big\vert \, A_i\in \mathfrak{g}^{ \star}\right\}$. The same grading is induced on the co-adjoint orbit $\hat{\mathcal{O}}_{{r}}^{\star}$.

\begin{remark}
Note that the degree of the grading is due to the choice of the pairing \eqref{eq:pairing} in the Takiff algebra. If we had chosen a different measure, say $\frac{\d z}{z^k}$, then the degree would have been $k$. 
\end{remark}

In the case when $\mathfrak{g}$ is $\mathfrak{gl_{m}}$ we have the following Poisson structure
\begin{equation}
\{(A_{\alpha})_{ij}, (A_{\beta})_{kl}\} = \left\{\begin{array}{ll}
    (A_{\alpha+\beta})_{il}\delta_{jk}-(A_{\alpha+\beta})_{kj}\delta_{il} & \alpha+\beta \leq r \\
    0 & \alpha+\beta > r,
\end{array}\right.
\end{equation}
which may be written in the $r$-matrix form
\begin{equation}\label{1.1}
\{A_{\alpha}\ocomma A_{\beta}\} =  \left\{\begin{array}{ll}
    -[{\Pi}, A_{\alpha+\beta}\otimes\ID] & \alpha+\beta \leq  r \\
    0 & \alpha+\beta > r.
    \end{array}\right.
\end{equation}

As mentioned before, the co-adjoint orbits of the Takiff algebra form the phase space of the isomonodromic deformation equations in the case of irregular singularities while in the Fuchsian $\hat{\mathcal{O}}^{\star}_{{0}}= \mathcal{O}^{\star}$.

\subsection{Lifted Darboux coordinates} 
As shown in the previous section, the lifted Darboux coordinates for the co-adjoint orbits of an {ordinary} Lie algebra are given by a symplectic reduction from $T^{\star} \mathfrak{gl}_m$. We prove the same result for the Takiff algebras. Our construction  follows ideas introduced by Chervov and Talalaev in \cite{TalCh} to parametrize the space of the irregular Gaudin systems.

We start from the following space
$$
\mathfrak{g}=\mathfrak{gl_m},\quad T^{\star}\hat{\mathfrak{g}}_{r} = \left\{ \left(P,Q\right)\Big\vert\,P = \sum\limits_{i=0}^{r} P_iz^{i},\,Q=\sum\limits_{i=0}^{r} Q_iz^{-i-1},\quad P_i,Q_i\in \mathfrak{gl}_m \right\}.
$$
The symplectic form on  $T^{\star}\hat{\mathfrak{g}}_{n}$ is given by the differential of the Liouville form:
\begin{equation}\label{sect3:symplForm}
    \omega = \d \langle P, \d Q \rangle = \oint\limits_{S^1}\tr\left( \d P\wedge \d  Q\right) d z = \d \sum\limits_{i=0}^{r} \tr\left(P_i \wedge \d Q_i\right),
\end{equation}
here $d$ is the differential on the space of the spectral parameter $z$, while $\d$ is the differential on the phase space.

\begin{lemma}\label{lm:uni-sp}
The map
$$
\begin{array}{ccc}
\bigoplus\limits_{i=0}^{{r}} T^{\star}\mathfrak{gl}_m& \to &T^{\star}\hat{\mathfrak{g}}_{{r}}\\
(P_0,\dots,P_r,Q_0,\dots,Q_r)&\mapsto &( P, Q)\\
\end{array}
$$
is a symplectomorphism.
\end{lemma}

The proof of this result is a straightforward consequence of the fact that $T^{\star}\hat{\mathfrak{g}}_{n}$ and  $\bigoplus\limits_{i=1}^n T^{\star}\mathfrak{gl}_m$
are isomorphic as vector spaces and formula (\ref{sect3:symplForm}) shows that they are  symplectomorphic to each other. However, we have enphasised this simple fact into a Lemma because $\bigoplus\limits_{i=1}^n T^{\star}\mathfrak{gl}_m$ provides the ambient space for the confluence procedure.\\

We now want to construct the Lie group $\hat G_{{r}}$ of the Takiff algebra. Its elements are given by:
$$
g(z) = g_0+\sum\limits_{i=1}^{r} g_{i}z^i,\quad g_0 \in GL_m,\quad g_i \in \mathfrak{gl}_m,
$$
where, in order to be able to multiply  both on the left and on the right, $\mathfrak{gl}_m$ is considered as a bi-module of $GL_m$.
The group structure {of $\hat G_n$} is given by   $GL_m$ multiplication mod $z^n$, i.e.
$$
g(z)\cdot h(z) = g(z)h(z) \operatorname{mod} z^{{r+1}} = g_{0}h_{0} + \sum\limits_{i=1}^{{r}}\left( \sum\limits_{j=0}^{i} g_{i-j}h_{j}\right)z^i.
$$
The inverse is given by
$$
g^{-1} = g_{0}^{-1}\left[1+\sum\limits_{i=1}^{r}g_{0}^{-1}g_{i+1}z^{i}\right]^{-1}= g_{0}^{-1}(1+\tilde{g}(z))^{-1} = g_{0}^{-1}\sum\limits_{i=0}^{\infty}(-1)^{i}\tilde{g}(z)^{i} \operatorname{mod} z^{r+1},
$$  
and the neutral element is given by the identity matrix. The induced inner and outer actions on $T^{\star}\hat{\mathfrak{g}}^{{r}}$ are given by
\begin{equation}\label{eq:outerTakiff}
    g \underset{\operatorname{outer}}{\times} \left( P , Q\right) = \left( \left[ P \circ g \right]\operatorname{mod} z^{{r+1}};\pi_{-} \left[g^{-1}\circ  Q \right] \right) 
\end{equation}
\begin{equation}\label{eq:innerTakiff}
    g \underset{\operatorname{inner}}{\times} \left( P , Q \right) = \left(\left[g\circ  P \right] \operatorname{mod} z^{{r+1}};\pi_{-}\left[ Q  \circ g^{-1}\right]\right) 
\end{equation}
where $\pi_{-}$ is a projection to the Laurent part with respect to spectral parameter $z$, i.e.
$$
\pi_{-}\left[\sum\limits_{i=-\infty}^{\infty}T_iz^i\right]=\sum\limits_{i=-\infty}^{-1}T_iz^i
$$

\begin{lemma}\label{lemma:liftedTakif}
Both inner and outer actions are Hamiltonian with the moment maps respectively
\begin{equation}
    \begin{array}{ll}\mu_{\operatorname{inner}}:& T^{\star}\hat{\mathfrak{g}}_{{r}} \to \hat{\mathfrak{g}}_{{r}}^\star \\
& ( P , Q )\mapsto \Lambda(z) = \pi_{-}\left[ P  Q \right]\end{array}\quad \begin{array}{ll}\mu_{\operatorname{outer}}:& T^{\star}\hat{\mathfrak{g}}_{r} \to \hat{\mathfrak{g}}_{r}^\star \\
& ( P , Q )\mapsto A(z) = \pi_{-}\left[ Q  P \right].\end{array}
\end{equation}
\end{lemma}
These two moment maps are dual in a sense of Adams--Harnad--Previato duality \cite{AdHarPr}. Since inner and outer group actions commute, $A(z)$ and $\Lambda(z)$ Poisson commute with respect to Poisson bracket induced by (\ref{sect3:symplForm}). As in the Fuchsian case, $A(z)$ is  an element of the co-adjoint orbit for the Takiff algebra. On the other hand, $\Lambda(z)$ becomes an invariant of the orbit after quotient via the inner group action. 

This fact gives us the opportunity to generalise the statement of Lemma \ref{sect1:lemma1} to the case of Takiff algebras:
\begin{lemma}\label{sect2:lemma1}
The map
$$
\begin{array}{ccc}
T^{\star}\hat{\mathfrak{g}}_{{r}} \underset{\Lambda_0}{\sslash} \hat{G}_{{r}}& \to &\hat{\mathcal{O}}_{{r}}^{\star}\\
( Q , P )&\mapsto&  A(z):=\pi_{-}\left[ Q  P \right] \\
\end{array}
$$
where $ \underset{\Lambda_0}{\sslash} $ denotes the Hamiltonian reduction w.r.t. the inner action in which the moment map has value  $\Lambda_0$,
is a rational symplectomorphism and the Jordan normal form $\Lambda_0$ of $A$ is given by
$$
\Lambda_0(z)=\pi_{-}\left[ P  Q\right].
$$
The explicit form of $A(z)$ is
{\begin{equation}\label{param}
A(z)=\frac{A_{{r}}}{z^{{r+1}}}\dots+\frac{A_{{0}}}{z},\quad A_{k} = \sum\limits_{i=0}^{r-k}\chi_{i,i+k},\quad \chi_{i,j} = Q_iP_j.
\end{equation}
while $\Lambda_0(z)$ takes form
\begin{equation}\label{eq:mom-prod}
\Lambda_0(z)=\frac{\Lambda_{{r}}}{z^{{r+1}}}\dots+\frac{\Lambda_{{0}}}{z},\quad \Lambda_{k} = \sum\limits_{i=0}^{{{r}}-k}P_{i+k}Q_{i},
\end{equation}}
\end{lemma}

\begin{remark} According to Lemma \ref{lm:uni-sp}, all co-adjoint orbits, i.e the ones for the ordinary Lie algebras and the one for the Takiff algebras, are reductions of the same phase space. Systems with different orders of poles are obtained by different choices of the group realising the reduction: 
 in the Fuchsian case we considered the action of the direct product of $GL_m$, while in the case of the Takiff algebra we use the inner action of $\hat G_m^{{r}}:=GL_m[z]\slash z^{{r+1}}GL_m[z]$.
\end{remark}

The parametrisation (\ref{param}) allows a nice combinatorial description which is presented on Fig. \ref{fig1:PQ}. 

\begin{figure}
    \centering
    \begin{tikzpicture}[baseline= (a).base]
    \node[scale=0.65] (a) at (0,0){
    \begin{tikzcd}
    \arrow[rd, color=blue]Q_{{0}}P_{{0}} & & Q_{{1}}P_{{1}} \arrow[ld, color=red] \arrow[rd, color=blue] & & \arrow[ld, color=red]Q_{{2}}P_{{2}} \arrow[rd, color=blue]& \dots & \dots & Q_{{r-1}}P_{{r-1}} \arrow[rd, color=blue]&  & \arrow[ld, color=red]Q_{{r}}P_{{r}} \\
    & Q_{{0}}P_{{1}} & & Q_{{1}}P_{{2}} &  & Q_{{2}}P_{{3}} & \dots & \dots & Q_{{r-1}}P_{{r}} & \\
    &  &\dots  & \dots & \dots & \dots & \dots & \dots & \\
    & &   & & Q_{i}P_{i+{k}-1}\arrow[rd, color=blue]& &\arrow[ld, color=red] Q_{i+1}P_{i+{k}} &  &  & \\
    & & &  & &  Q_{i}P_{i+{k}} & \, & & \arrow[ll, color=green] A_{{k}} \\
    & & & & & \dots &  & & \\
    & & & & & Q_{{0}}P_{{r}} & & &  
\end{tikzcd}};
\end{tikzpicture}
    \caption{Lifted Darboux coordinates for the Takiff algebra of degree $r$. In this diagram we have $r+1$ rows, and we number them starting at the top with row $0$, all the way down to row $r$. 
The sum of the elements in row $k$ gives the coefficient $A_{{k}}$ of the power of $z^{-{{k}}-1}$, the blue arrow  follows each $Q_i$ matrix from the formula above to the one below, while the red one follows $P_i$.}
    \label{fig1:PQ}
\end{figure}

\begin{theorem}\label{th:QPKKS}
The Poisson bracket induced by the Darboux coordinates $Q_{i},P_{i}$ on the space of matrices $A_{{k}}$, $k=0,\dots,r$  coincides with the graded Poisson structure (\ref{1.1}).
\end{theorem}
\begin{proof}
This statement is a straightforward corollary of the Lemma \ref{sect2:lemma1}. However, here we prove it directly for the sake of clarity.
The Poisson bracket on the elements $\chi_{ij}$ in \eqref{param} is given by
\begin{multline*}
\{\chi_{ij}\ocomma \chi_{kl}\} = \{Q_{i}P_{j}\ocomma Q_{k}P_{l}\} = \delta_{jk}(Q_{i}\otimes1){\Omega}( \ID\otimes P_{l}) - \delta_{il}  ( \ID\otimes Q_{k}){\Omega}(P_{j}\otimes  \ID) = \\ = \delta_{jk}(Q_{i}P_{l}\otimes \ID){\Omega}-\delta_{il}{\Omega}(Q_{k}P_{j}\otimes  \ID) = \delta_{jk}(\chi_{il}\otimes \ID){\Omega}-\delta_{il}{\Omega}(\chi_{kj}\otimes  \ID)
\end{multline*}
which is the same as
$$
\left\{\left(\chi_{ij}\right)_{\alpha\beta}, \left(\chi_{kl}\right)_{\gamma\delta}\right\} = \delta_{jk}\delta_{\gamma \beta}\left(\chi_{il}\right)_{\alpha\delta} - \delta_{il}\delta_{\alpha\delta}\left(\chi_{kj}\right)_{\gamma\beta}.
$$
By direct computation
\begin{multline*}
\{A_k\ocomma A_l\} = \sum\limits_{i,j}\{\chi_{i,i+k}\ocomma \chi_{j,j+l}\} = \sum\limits_{i,j}\delta_{j, i+k}(\chi_{i,j+l}\otimes\ID)\Pi - \delta_{i,j+l}\Pi(\chi_{j,i+k}\otimes \ID) = \\
= \sum\limits_{i}(\chi_{i,i+k+l}\otimes\ID)\Pi - \sum\limits_{j}\Pi(\chi_{j,j+k+l}\otimes \ID) = -[\Pi, A_{k+l}\otimes\ID]
\end{multline*}
we obtain the proof of the statement. When $k+l>{{r}}$ the Poisson bracket is automatically zero.
\end{proof}

In the next lemma, we show that the quadratic Casimir elements for the Takiff algebra are given by functions of the spectral invariants of the co-adjoint orbit:

\begin{lemma}
For the Takiff algebra of degree $r$, the following quantities are Casimirs
\begin{equation}\label{eq:quadCas}
    I_k = \underset{z=0}{\res} \left(z^{r+k}\tr A^2\right),\quad 0<k<r.
\end{equation}
\end{lemma}

\proof
The fact that $I_k$ are Casimir functions may be checked by  direct computation. Here we demonstrate it for $k=1$ since we use this fact later in the  text. Explicitly $I_1$ writes as follows
$$
I_1 = \sum\limits_{j=0}^r\tr A_j A_{r-j}.
$$
The Poisson bracket with an arbitrary generator of the Poisson algebra defined via Lie-Poisson bracket for the Takiff algebra gives
\begin{multline*}
\sum\limits_{j=0}^r\left\{(A_i)_{\alpha}, \tr A_j A_{r-j}\right\} = \sum_{j=0}^{r}[A_{i+j}, A_{r-j}]_{\alpha} + \sum_{j=i}^r[A_{r-j+i},A_j]_{\alpha} = \\ =\sum_{l=i}^{r}[A_l,A_{r-l+i}]_{\alpha} + \sum_{j=i}^r[A_{r-j+i},A_j]_{\alpha} = 0.
\end{multline*}
In the same way we may prove that $I_k$ are the Casimirs for $k>1$.
\endproof

\subsection{Poisson automorphisms of the Takiff algebra and independent deformation parameters.}
In this subsection, we describe the class of linear automorphisms of the Takiff algebra which preserve the Poisson bracket, namely linear maps
\begin{equation}\label{sect2:lin_comb}
B_i = \sum\limits_{j={{0}}}^r T_{ij}A_j,\qquad T_{ij}\in\mathbb C,\quad  i,j={{0}},\dots {{r}},
\end{equation}
such that 
\begin{equation}\label{sect2:Poisson_constr}
    {\{B_i\ocomma B_j\} = [\Pi, \ID\otimes B_{i+j}]\quad\iff\quad    \{A_i\ocomma A_j\} = [\Pi, \ID\otimes A_{i+j}]}
\end{equation}
In the next theorem we describe explicitly the constraints on the coefficients $T_{ij}$.

\begin{theorem}\label{th:aut0Tak}
The coefficients $T_{ij}$ characterising the class of linear automorphisms of the Takiff algebra which satify the Poisson condition (\ref{sect2:Poisson_constr}) define an ideal 
$\mathcal{P}$ in the ring $\mathbb{C}[T_{11}\dots T_{rr}]$ given by the equations
\begin{equation}\label{sect2:ideal}
   \mathcal{P} =\left\{
   \begin{array}{ll}
        T_{00} = 1 & \\
        T_{0k}=0, & k>{{0}} \\
        T_{k0}=0, & k>{{0}}\\
        T_{ik} = 0, &   k {{<}} i \\
        T_{sl} =\sum\limits_{i,j>{0}}^{i+j=l} T_{pi}T_{mj} & \forall p,m>{0}:\, {p+m=s}.
    \end{array}\right.
\end{equation}
Moreover we have the following ring isomorphism for the quotient
\begin{equation}\label{eq:birat}
    \mathcal{Q}:\quad \mathbb{C}[T_{00}\dots T_{rr}]\slash \mathcal{P} \rightarrow \mathbb{C}[t_1\dots t_{r}]
\end{equation}
such that
\begin{equation}\label{eq:tki}
     T_{1i} = t_{i},\quad T_{ki} =  \frac{1}{i!}\frac{d^i}{d\varepsilon^i} P_r(t,\varepsilon)^k\Big\vert_{\varepsilon=0},\quad P_r(t,\varepsilon) = \sum\limits_{i=1}^r\varepsilon^it_i,
\end{equation}
so that $T_{ki}$ is just the coefficient of the $\varepsilon^i$ term in the polynomial $P_r(t,\varepsilon)^k.$
\end{theorem}
\begin{remark}
The equations which define the ideal $\mathcal{P}$ do not depend on the  specific form of $\Pi$, i.e. on the structure constants of a Poisson bracket. Therefore,  the classification of the automorphisms is a consequence of the grading structure and not a property of the specific Lie co-algebra.
\end{remark}

\begin{proof}
Assume the matrices $A_i$ and $B_i$ satisfy the Poisson relations \eqref{sect2:Poisson_constr} and prove the relations for the coefficients $T_{ij}$. Let us start from the relation for $B_1$
\begin{equation}\label{eq:prB1PB}
    \{B_{{0}}\ocomma B_{{0}}\} = [\Pi, \ID\otimes B_{{0}}].
\end{equation}
Substituting \eqref{sect2:lin_comb} in \eqref{eq:prB1PB} and expanding, we obtain
\begin{multline}
     \{B_{{0}}\ocomma B_{{0}}\} =\sum\limits_{i,j=0}^rT_{0i}T_{0j}\{A_i\ocomma A_j\} = \sum\limits_{k={{0}}}^r \left(\sum\limits_{i=0}^{k}T_{0i}T_{0,k-i}\right) \left[\Pi,\ID\otimes A_k\right] = \\
    = [\Pi, \ID\otimes B_{{0}}] = \sum\limits_{k=1}^rT_{0k}[\Pi, \ID\otimes A_k].
\end{multline}
This relation defines a system  of equations for the coefficients $T_{0j}$, which takes the form
$$
    T_{00}T_{00}= T_{00},\quad 2T_{00}T_{0k} + \sum\limits_{i=1}^{k-1}T_{0i}T_{0,k-i} = T_{0k},
$$
that, by recursion, leads to the first set of equations which generate the ideal $\mathcal{P}$:
$$
T_{00}=1,\quad T_{0k}=0,\quad k>{{0}}.
$$
The next statement we want to prove is that $T_{{k0}}=0$ for $k>1$. We use 
\begin{equation}
    \{B_{{1}}\ocomma B_k\} = [\Pi, \ID\otimes B_{k+1}]\qquad k=1,\dots,{{r}}.
\end{equation}
Again, substituting \eqref{sect2:lin_comb} and expanding, we obtain
$$
\sum\limits_{i,j={{0}}}^{{r}}T_{1i}T_{kj} [\Pi, \ID\otimes A_{i+j}]=\sum_{j={{0}}}^{r} T_{k+1,j} [\Pi, \ID\otimes A_j],
$$
and collecting all coefficients of $[\Pi, \ID\otimes A_1]$,
we have that 
$$
T_{k+1,0} = T_{10}T_{k0},
$$
that is solved by
$$
T_{k+1,0} = \left(T_{10}\right)^{k}.
$$
On the other hand substituting \eqref{sect2:lin_comb} in
\begin{equation}
    \{B_{{1}}\ocomma B_r\}  = 0,
\end{equation}
we obtain
$$
T_{10}T_{r0} = 0 = \left(T_{10}\right)^r \quad\Rightarrow\quad T_{10}=0,
$$
as we wanted. Now to demonstrate the statement that $T_{ik}=0$ for $k{{<}} i$ we use the relation
\begin{equation}\label{sect2:th1}
    \{B_{{1}}\ocomma B_{{1}}\} = [\Pi, \ID\otimes B_{{2}}].
\end{equation}
By substituting \eqref{sect2:lin_comb} we see that the left hand side of  (\ref{sect2:th1})
\begin{equation}\label{sect2:th2}
    \{B_{{1}}\ocomma B_{{1}}\} = \sum_{i,j>2}T_{1i}T_{1j} [\Pi, A_{i+j}]  = T_{11}T_{11} [\Pi, A_{{2}}] + \sum_{i=4}\kappa_i[\Pi, A_i] ,
\end{equation}
does not contain terms in $A_{{0}}$ or $A_{{1}}$, it contains only one term that depends on $A_{{2}}$, given by $ T_{11}T_{11} [\Pi, A_{{2}}]$ and all other terms depend on $A_{{3}},\dots,A_r$. Expanding the right hand side of  (\ref{sect2:th1}) we obtain
\begin{equation}\label{sect2:th3}
    [\Pi, \ID\otimes B_{{2}}] = T_{21}[\Pi, \ID\otimes A_{{1}}]+\sum\limits_{i={{2}}}T_{2i}[\Pi, \ID\otimes A_i].
\end{equation}
Therefore $ T_{21} = 0.$ Similarly, 
applying the $\{B_{{1}}\otimes \circ \}$  to $B_{{2}}\dots B_{r}$ and using the same approach we obtain that $T_{ik}=0$ for $k< i$.
The last relation in \eqref{sect2:ideal} is obtained by imposing \eqref{sect2:Poisson_constr}, substituting  \eqref{sect2:lin_comb}  and expanding as before, and then by imposing all other conditions we have obtained so far.

We now prove the second part of the Theorem. First of all, we observe that thanks to relations \eqref{sect2:ideal},
the coefficients $t_j:=T_{1j}$  for $j>0$ form a basis in the quotient ring $ \mathcal{Q}:\quad \mathbb{C}[T_{00}\dots T_{rr}]\slash \mathcal{P}$. Then, because each 
$T_{ik}$ must be given by a polynomial $P_k^{(i)}$of $t_1,\dots,t_{r}$, we just need to check the degree and the form of the coefficients.
To this aim we use the last relation of  \eqref{sect2:ideal} for $ T_{ij}$ by induction on $j$ from $i$ to $r$. We omit this computation as it is straightforward.
\end{proof}

In the next section we will see how such dependence on the parameters $t_i$'s arises during the confluence procedure. In some sense, the irregular deformation parameters are just the deformation of the representation for the Takiff algebra.

\begin{example}\label{ex.tak234} In order to give a taste of how the general elements of the Takiff co-algebra depend on the Poisson automorphism parameters $t_i$, we provide a few examples of low degree. We consider an element of the Takiff co-algebra as a polynomial in $\frac{1}{z}$. In the case of $\hat{\mathfrak{g}}_1^{\star}$ Theorem \ref{th:aut0Tak} gives
\begin{equation}
B(z) = \frac{t_1A_1}{z^2}+\frac{A_0}{z}.
\end{equation}
In this case, we see that the invariant space of the action of $A_0$ is defined up to multiplication by a constant, so this example is quite trivial. Let us look at   $\hat{\mathfrak{g}}_2^{\star}$. In this case, the general element writes as
\begin{equation}
    B(z) = \frac{t_1^2A_2}{z^3}+\frac{t_1A_1+t_2 A_2}{z^2}+\frac{A_0}{z}.
\end{equation}
\end{example}

\begin{example}\label{example-auto-conf} The next example is the case of $\hat{\mathfrak{g}}_3^{\star}$ where the element of the co-algebra writes as
\begin{equation}\label{eq-ex-conf}
        B(z) = \frac{t_1^3A_3}{z^4}+\frac{t_1^2A_2+2t_1t_2 A_3}{z^3}+\frac{t_1A_1+t_2A_2+t_3A_3}{z^2}+\frac{A_0}{z}.
\end{equation}
Let us now see how to obtain these formulae for $B(z)$ starting from the following connection with a pole of order $4$ at zero
    $$
    A(z)dz = \left(\frac{A_3}{z^4}+\frac{A_2}{z^3}+\frac{A_1}{z^2}+\frac{A_0}{z} + O(1)\right) dz
    $$
    and performing a local conformal change of coordinates
    $$
    z = f(\zeta)= \frac{\zeta }{\tau_1}+\zeta^2\tau_2+\zeta^3\tau_3+O(\zeta^4).
    $$
    Then 
    \begin{multline*}
    B(z) = A(f(\zeta))f'(\zeta)d\zeta = \\=\left(\frac{\tau_1^3A_3}{\zeta^4}+\frac{\tau_1^2A_2-2\tau_1^4\tau_2A_3}{\zeta^3}+\frac{\tau_1A_1-\tau_1^3\tau_2A_2+(2\tau_1^5\tau_2^2-\tau_1^4\tau_3)A_3}{\zeta^2}+\frac{A_0}{\zeta} + O(1)\right) d\zeta.
    \end{multline*}
    It may explicitly checked that such transformation preserves the Poisson bracket. The parameters $\tau_1, \tau_2$ and $\tau_3$ are related to $t_1,t_2,t_3$ in \eqref{eq-ex-conf} by the following bi-rational map
    $$
    \tau_1=t_1, \quad \tau_2 = -t_2/t_1^3, \quad \tau_3 = -t_3/t_1^4-2t_2^2/t_1^5.
    $$   
\end{example}

\subsection{Direct product of the co-adjoint orbits and outer group action}

To study systems with more than one pole, we will need to consider the symplectic space given by the direct product of different co-adjoint orbits of the Lie algebra for simple poles, or of the appropriate Takiff algebra for higher order poles. We use here a unified notation, in which we understand that for poles of order $1$, the Poincar\'e rank is $r=0$ and
$\hat{\mathfrak g}_{0}$ is $\mathfrak g$, $\hat{\mathcal{O}}_{0}^{\star}$ is $\mathcal{O}_i^\star$, $T^{\star}\hat{\mathfrak{g}}_{0}$ is $T^{\star}{\mathfrak{g}}$ and $\hat{G}_0$ is $G$. With this notation in mind, the symplectic space we consider is 
\begin{equation}\label{eq:prod-coadj-many}
   \hat{\mathcal{O}}_{r_1}^{\star}\times\hat{\mathcal{O}}_{r_2}^{\star}\times\dots \hat{\mathcal{O}}_{r_n}^{\star}\times  \hat{\mathcal{O}}_{r_\infty}^{\star},
\end{equation}
where we always assume to have a pole at infinity like in the Fuchsian case.
This product of co-adjoint orbits may be viewed as the reduction of the {\it universal symplectic space} $\bigoplus\limits_{i=1}^{n} T^{\star}\hat{\mathfrak{g}}_{r_i}$
with respect to the inner action of the group $\mathcal G^{(n)}:=\hat{G}_{r_1}\times \hat{G}_{r_2} \times \dots \times \hat{G}_{r_n}\times \hat{G}_{r_\infty}:$
$$
\hat{\mathcal{O}}_{r_1}^{\star}\times\hat{\mathcal{O}}_{r_2}^{\star}\times\dots \hat{\mathcal{O}}_{r_n}^{\star}\times  \hat{\mathcal{O}}_{r_\infty}^{\star} = \left(\bigoplus\limits_{i=1,\dots,n,\infty} T^{\star}\hat{\mathfrak{g}}_{r_i}\right) \underset{\otimes\Lambda_{r_i}}{\sslash} \mathcal G^{(n)}.
$$
Since we have the following symplectomorphism
$$
\bigoplus\limits_{i=1}^{n} T^{\star}\hat{\mathfrak{g}}_{r_i}\simeq \bigoplus\limits_{i=1}^{r_1+r_2+\dots+r_n+r_\infty} T^{\star}\mathfrak{gl}_m,
$$
we obtain that
$$
\hat{\mathcal{O}}_{r_1}^{\star}\times\hat{\mathcal{O}}_{r_2}^{\star}\times\dots \hat{\mathcal{O}}_{r_n}^{\star}\times  \hat{\mathcal{O}}_{r_\infty}^{\star} \simeq \left(\bigoplus\limits_{i=1}^{r_1+\dots +r_n+r_\infty} T^{\star}\mathfrak{gl}_m\right) \underset{\otimes\Lambda_{r_i}}{\sslash} \mathcal G^{(n)}
$$
where we denote by $\underset{\otimes\Lambda_{r_i}}{\sslash}$ the Hamiltonian reduction with respect to the inner action in which the value of the moment map is given by the product of values $\Lambda_{r_i}$ of  the inner moment map for each 
$ \hat{G}_{r_i}$.\\

We now take into account the outer action on each co-adjoint orbit; similarly to the Fuchsian case, in order to have a well defined action on the whole connection, we again restrict to the diagonal case 
$$
g_1=g_2=\dots = g_n=g_\infty = g,
$$
where $g_i$ doesn't depend on the spectral parameter $z_i$.  Therefore, this constant diagonal action is the constant gauge group $G$ action as in the Fuchsian case. 

The moment map of this constant diagonal outer action takes the form
$$
\mu = \sum_{j=1,\dots,n,\infty}\sum\limits_{i=0}^{r_j} Q_i^{(j)}P_i^{(j)},
$$
which may be again seen as the sum of residues at poles. Finally, the fully reduced space takes the form
\begin{equation}\label{sect3:irr_space}
    M = \hat{\mathcal{O}}_{r_1}^{\star}\times\hat{\mathcal{O}}_{r_2}^{\star}\times\dots \hat{\mathcal{O}}_{r_n}^{\star} \times\hat{\mathcal{O}}_{r_\infty}^{\star} \sslash \mathcal G^{(n)} \simeq \left[\left(\bigoplus\limits_{i=1}^{r_1+\dots +r_n+r_\infty} T^{\star}\mathfrak{gl}_m\right) \underset{\Lambda}{\sslash} \mathcal G^{(n)}\right] \sslash \mathcal G^{(n)}.
\end{equation}
The quotient with respect to the diagonal outer action has the same effect as in the Fuchsian case - it specifies the residue at the infinity. However, differently from the Fuchsian case, where this was enough to fully characterise the Fuchsian singularity at infinity, here we have a pole of arbitrary Poincare rank $r$ at infinity, where the connection takes the form
$$
A(\lambda) = \lambda^{r-1}A^{(\infty)}_{r+1}+\sum_{k=0}^{r-1}A^{(\infty)}_{k+2}\lambda^k + \text{regular terms at }\infty.
$$
We may view the moment map as fixing the term $A^{(\infty)}_1$. In the next section we will study the isomonodomic deformations of irregular connections that are elements of the space (\ref{sect3:irr_space}).

\subsection{Fixing the spectral invariants. Reduction with respect to the inner action.}

In this section we compute explicitly the reduced coordinates for the co-adjoint orbits of the quotient of Takiff algebras with respect to the inner group action on the lifted Darboux coordinates in the case of degrees 1, 2, 3 and 4 - this choice is motivated by the fact that in the Painlev\'e confluence scheme the maximal pole order we have is 4. However, the described procedure can be easily expanded for the Takiff algebra of any degree - we give a hint and some explanation in the discussion after  the examples. In each example we give explicit results in the case of $\mathfrak{sl}_2$, since this is the case of the isomonodromic problems for the Painlev\'e equations. We also provide the coordinates in the diagonal gauge - the case when the leading term is diagonal by using the additional outer action of the gauge group $G$.

\subsubsection{First order pole. Takiff algebra of degree 1} In this case Takiff algebra coincide with the ordinary Lie algebra. The parametrisation in such situation was obtained in works \cite{Bab,BabDer}.

\subsubsection{Second order pole. Takiff algebra of degree $1$}\label{suse:deg2}

The Darboux parametrisation is given by
$$
A(z) = \frac{Q_0P_1}{z^2}+\frac{Q_0P_0+Q_1P_1}{z},\quad \omega = \d\Theta,\quad  \Theta = \tr\left(P_1\d Q_1 + P_0\d Q_0 \right),
$$
so that the extended phase space is of dimension $4m^2$. We now want to reduce this dimension by solving the moment map conditions
$$
P_1Q_0=\Lambda_1,\quad P_0Q_0+P_1Q_1=\Lambda_0
$$
w.r.t. $P_0$ and $P_1$. To do this, we only need to assume that $Q_0$ is invertible, namely  $(Q_0,Q_1)\in \bigoplus{Gl}_m\times \mathfrak{gl}_m$.
This inversion sends the Liouville form to
$$
\theta = \tr\left(\Lambda_1Q_0^{-1}\d Q_1+\Lambda_1Q_0^{-1}\d Q_0 - \Lambda_1Q_0^{-1}Q_1Q_0^{-1}\d Q_0 \right),
$$
while 
$$
A(z) \mapsto \frac{Q^{\vphantom{-1}}_{0}\Lambda^{\vphantom{-1}}_{1}Q_{0}^{-1}}{z^2}+\frac{Q_0\Lambda_0Q_0^{-1}+[Q_{1}Q_0^{-1},Q_{0}\Lambda_{1}Q_{0}^{-1}]}{z}.
$$
We now want to reduce the dimension by $2m$ via the torus action $Q_i\to Q_i D_i$, where $D_i$ is a diagonal matrix, that fixes the invariants of the co-adjoint orbit $\Lambda_0,\Lambda_1$.
To this aim, we find the Darboux coordinates  $p_1,\dots p_{m(m-1)},q_1,\dots q_{m(m-1)}$ explicitly in such a way that
\begin{equation}
    \Theta = \tr\left(\Lambda_1Q_0^{-1}\d Q_1+\Lambda_0Q_0^{-1}\d Q_0 - \Lambda_1Q_0^{-1}Q_1Q_0^{-1}\d Q_0 \right) = \sum\limits_{i=1}^{m(m-1)}p_i\d q_i.
\end{equation}
The number of unknown functions also equals to $2m(m-1)$, due to the factorisation of the torus action. There are many possible choices 
for the Darboux coordinates $p_1,\dots p_{m(m-1)},q_1,\dots q_{m(m-1)}$ in this situation, our aim to find one good choice; it is convenient to use the following change
$$
L_1 = Q_0^{-1}Q_1,
$$
then Liouville form transforms to
$$
\Theta = \tr \left[\Lambda_1 \d L_1 +\left(\Lambda_0 + [L_1,\Lambda_1] \right)Q_0^{-1}\d Q_0 \right].
$$
The Liouville form is always defined up to a closed form. Since $\Lambda_1$ is an invariant of the co-adjoint orbit (i.e. is a constant) the term
$$
\Lambda_1 \d L_1 = \d \left(\Lambda_1L_1\right)
$$
is exact, so we may drop it. The equation for the differential form therefore simplifies to
$$
\tr \left[\left(\Lambda_0 + [L_1,\Lambda_1] \right)Q_0^{-1}\d Q_0\right] = \sum\limits_{i=1}^{m(m-1)}p_i\d q_i,
$$
which allows us to pick our Darboux coordinates $p_1,\dots p_{m(m-1)},q_1,\dots q_{m(m-1)}$ in such a way that $Q_0$ depends only on $q_1,\dots q_{m(m-1)}$ (i.e. $Q_0$  is a section of a principal bundle over the Lagrangian sub-manifold), while the entries of $L_1$ are given by the solutions of $m(m-1)$ linear equations. For example
we may take the off-diagonal entries of $Q_0$ as the coordinates on the Lagrangian sub-manifold. By using the torus action, we can make the following choice for $Q_0$:
$$
Q_0 = \left(\begin{array}{ccccc}
    1 & q_1& \dots &\dots&  q_{m-1} \\
    0& 1& q_{m}&\dots&q_{2m-3}\\
      \vdots&0 &\ddots &\ddots& \vdots \\
       0&\dots &0&1&q_{\frac{m(m-1)}{2}}\\
           0&\dots &\dots &0&1\\
\end{array}\right)\left(\begin{array}{ccccc}
    1 &  0 &\dots & \dots &0\\
    q_{\frac{m(m-1)}{2}+1} & 1 &  0 &\dots & 0\\
     q_{\frac{m(m-1)}{2}+2} & q_{\frac{m(m-1)}{2}+3}  &  1&\ddots & 0\\
       \vdots&\vdots &\ddots &\ddots& \vdots \\
           q_{(m-1)^2}&\dots &\dots &q_{m(m-1)}& 1\\
\end{array}\right).
$$
For $\mathfrak{sl}_2$ we have
$$
\Lambda_i = \left(\begin{array}{cc}
    \theta_i & 0 \\
    0& -\theta_i
\end{array}\right),\quad
Q_0 = \left(\begin{array}{cc}
    1 & q_1 \\
    0& 1
\end{array}\right)\left(\begin{array}{cc}
    1 &  0 \\
    q_2 & 1 
\end{array}\right),\quad L_1 =  \frac{1}{2\theta_1}\left(\begin{array}{cc} 
    0 & -p_{{2}}\\ 
    p_{{2}}q_{2}^{2}-2\theta_{{0}}q_{{2}}+p_{{1}} & 0
    \end{array}\right)
$$
and the matrix $A(z)$ takes the following form
\begin{multline}
A(z) = \frac{2\theta_1}{z^2}\left(\begin {array}{cc} q_{{1}}q_{{2}}+1/2&- \left(q_{{1}}q_{{2}}+
1 \right) q_{{1}}\\ \noalign{\medskip}q_{{2}}&-q_{{1}}q_{{2}}-1/2
\end {array}
\right) + \\ +\frac{1}{z}\left(\begin {array}{cc} p_{{1}}q_{{1}}-q_{{2}}p_{{2}}+\theta_{{0}}
&-p_{{1}}q_{1}^{2}+ \left( 2q_{{1}}q_{{2}}+1 \right) p_{{2}}-2
\theta_{{0}}q_{{1}}\\ \noalign{\medskip}p_{{1}}&-p_{{1}}q_{{1}}+q_{{2
}}p_{{2}}-\theta_{{0}}\end {array}
\right).
\end{multline}
If we take into account the outer action of $SL_2$, the leading term can be chosen in diagonal form and we have
$$
Q_1^{-1}A(z)Q_1 = \frac{\theta_1}{z^2}\left(\begin{array}{cc}
    1 & 0 \\
    0 & -1
\end{array}\right)+\frac{1}{z}\left(\begin {array}{cc} \theta_{{0}}&p_{{2}}\\ \noalign{\medskip}p
_{{2}}q_{2}^{2}-2\,\theta_{{0}}q_{{2}}+p_{{1}}&-\theta_{{0}}
\end {array} 
\right).
$$

\subsubsection{Third order pole. Takiff algebra of degree  $2$} In this case, the parametrisation in terms of lifted Darboux coordinates is given by
$$
A(z) = \frac{Q_0P_2}{z^3}+\frac{Q_0P_1+Q_1P_2}{z^2}+\frac{Q_0P_0+Q_1P_1+Q_2P_2}{z},
$$
so that the extended phase space is of dimension $6m^2$. The moment map is given by the equations
$$
P_2Q_0=\Lambda_2,\quad P_1Q_0+P_2Q_1=\Lambda_1,\quad 
P_0Q_0+P_1Q_1+P_2Q_2=\Lambda_0.
$$
Here we again use the following change of variables
$$
L_1=Q_0^{-1}Q_1,\quad L_2=Q_0^{-1}Q_2
$$
that maps the Liouville form to
$$
\Theta =\tr\left(\Lambda_2 \d L_2 + \Lambda_1 \d L_1 - \Lambda_2L_1\d L_1 + \left(\Lambda_0+[L_2,\Lambda_2]+[L_1,\Lambda_1-\Lambda_2 L_1]\right)Q_0^{-1}\d Q_0\right).
$$
As in the previous case, the first $2$ terms are closed differential forms, so we can drop them. The dimension of the reduced phase space equals to $3m(m-1)=3N$ and we consider the following parametrisation
$$
\tr(- \Lambda_2L_1\d L_1 ) = \sum\limits_{i=1}^{N/2} p_i\d q_i,\quad \tr\left[\left(\Lambda_0+[L_2,\Lambda_2]+[L_1,\Lambda_1-\Lambda_2 L_1]\right)Q_0^{-1}\d Q_0\right] = \sum\limits_{i=N/2+1}^{3N/2} p_i\d q_i.
$$
For simplicity, let us denote 
$$
\Theta_1=\tr(- \Lambda_2L_1\d L_1 ),\quad \Theta_2= \tr\left[\left(\Lambda_0+[L_2,\Lambda_2]+[L_1,\Lambda_1-\Lambda_2 L_1]\right)Q_0^{-1}\d Q_0\right], 
$$ 
so that  $\Theta = \Theta_1+\Theta_2$. Now if we will find the right parametrisation of $L_1$, we may choose $Q_0$ to be a matrix which depends only on $q_{N/2+1},\dots q_{3N/2}$ (i.e. again $Q_0$ depends only on the coordinates of the Lagrangian sub-manifold) and then obtain $L_2$ by solving a system of linear equations. 
In the non-degenerate case, when $\Lambda_2$ is a semi-simple matrix with distinct eigenvalues $\zeta_i$, we have
\begin{multline*}
\Theta_1= \sum\limits_{i<j} -\zeta_i (L_1)_{ij}\d (L_1)_{ji}-\zeta_j (L_1)_{ji}\d (L_1)_{ij} =\\ \sum\limits_{i<j}(\zeta_i-\zeta_j)(L_1)_{ji}\d (L_1)_{ij} - \d (\zeta_i(L_1)_{ij}(L_1)_{ji}) \sim \sum\limits_{i<j}(\zeta_i-\zeta_j)(L_1)_{ji}\d (L_2)_{ij}
\end{multline*}
and we see that a natural choice of the Darboux coordinates are the off-diagonal entries of $L_1$, such that
$$
\{(L_1)_{ij} , (L_1)_{kl}\} =\operatorname{sgn}(j-i)\delta_{kj}\delta_{li} (\zeta_i-\zeta_j).
$$
In the case of $\mathfrak{sl_2}$ we have
$$
\Lambda_2 = \left(\begin{array}{cc}
    \theta_2 & 0 \\
    0 & -\theta_2
\end{array}\right),\quad L_1 = \left(\begin{array}{cc}
    \dots & q_1 \\
    \frac{p_1}{2\theta_2} & \dots
\end{array}\right).
$$
Here the diagonal part of $L_1$ is irrelevant, since it does not contribute to $\Theta_1,\Theta_2$  and it may be chosen to be zero by the torus action. Solving the linear equations for the Cartan form $\Theta_2$ we obtain
$$
\Lambda_i = \left(\begin{array}{cc}
    \theta_i & 0 \\
    0& -\theta_i
\end{array}\right),\quad
Q_0=\left(\begin{array}{cc}
    1 & q_1 \\
    q_2 & 1
\end{array}\right)
$$
$$
L_1=\frac{1}{2\theta_2}\left(\begin{array}{cc}
    0 & (p_{2} q_{2}+p_3 q_3-\theta_{0}) q_{1}-p_{2}+\frac{\theta_1}{\theta_2}p_3 \\
    p_{1}-p_{1} q_{1} q_{2}+(p_{3} q_3-\theta_{0}) q_{2}-2 \theta_{1} q_3
 & 0
\end{array}\right)
$$
Here we take in a slightly different form  of $Q_0$ respect to in the previous example for the sake of obtaining a neater final formula. 
The matrix $A(z)$ takes form
\begin{multline}
  A(z)= \frac{1}{z^3} \frac{1}{1-q_1q_2}\left( \begin {array}{cc}  \theta_{{2}} \left( q_{{1}}q_{{2}}+1
 \right) &-2\, \theta_{{2}}q_{{1}}\\ 2\,q_{{2}}
 \theta_{{2}}&- \theta_{{2}} \left( q_{{1}}q_{{2}}+1 \right) 
\end {array} \right) + \\ +
 \frac{1}{z^2} \frac{1}{1-q_1q_2}
 \left( \begin {array}{cc}  \theta_{{1}}q_{{1}}q_{{2}}+2\, \theta_{{2}
}q_{{1}}q_{{3}}-q_{{2}}p_{{3}}+ \theta_{{1}}&-2\,q_{1}^{2}q_{{3}}
 \theta_{{2}}-2\, \theta_{{1}}q_{{1}}+p_{{3}}\\ -q_
{2}^{2}p_{{3}}+2\, \theta_{{1}}q_{{2}}+2\, \theta_{{2}}q_{{3}}&-
 \theta_{{1}}q_{{1}}q_{{2}}-2\, \theta_{{2}}q_{{1}}q_{{3}}+q_{{2}}p_{{
3}}- \theta_{{1}}\end {array}
 \right)  +\\+\frac{1}{z}
  \left( \begin {array}{cc} p_{{1}}q_{{1}}-q_{{2}}p_{{2}}-p_{{3}}q_{{3}
}+ \theta_{{0}}&-p_{{1}}{q_{{1}}}^{2}+p_{{3}}q_{{1}}q_{{3}}- \theta_{{0
}}q_{{1}}+p_{{2}}\\  -p_{{2}}{q_{{2}}}^{2}-p_{{3}}q_{
{2}}q_{{3}}+ \theta_{{0}}q_{{2}}+p_{{1}}&-p_{{1}}q_{{1}}+q_{{2}}p_{{2}
}+p_{{3}}q_{{3}}- \theta_{{0}}\end {array} \right) 
\end{multline}
The diagonal gauge gives
\begin{multline}
Q_0^{-1}A(z)Q_0 =\frac{1}{z^3}\left(\begin{array}{cc}
    \theta_3 &  0\\
    0 & -\theta_3
\end{array}\right) + \frac{1}{z^2}\left( \begin {array}{cc}  \theta_{{1}}&p_{{3}}\\ 2
\, \theta_{{2}}q_{{3}}&- \theta_{{1}}\end {array} \right) \\+\\\frac{1}{z}\left( \begin {array}{cc} -p_{{3}}q_{{3}}+ \theta_{{0}}&-p_{{2}}q_{{1
}}q_{{2}}-p_{{3}}q_{{1}}q_{{3}}+ \theta_{{0}}q_{{1}}+p_{{2}}
\\ \noalign{\medskip}-p_{{1}}q_{{1}}q_{{2}}+p_{{3}}q_{{2}}q_{{3}}-
 \theta_{{0}}q_{{2}}+p_{{1}}&p_{{3}}q_{{3}}- \theta_{{0}}\end {array}
 \right) 
\end{multline}
Choosing a different parameterisation for $Q_0$, i.e.
$$
Q_0= \left(\begin{array}{cc}
    1 & q_1 \\
    0 & 1
\end{array}\right)\left(\begin{array}{cc}
    1 & 0 \\
    q_2 & 1
\end{array}\right)
$$
the system takes the form
$$
Q_0^{-1}A(z)Q_0 = \frac{\theta_3}{z^3}\left(\begin{array}{cc}
    1 & 0 \\
    0 & -1
\end{array}\right)+\frac{1}{z^2}\left(\begin {array}{cc} \theta_{{2}}&-2\theta_{{3}}q_{{1}}
\\ \noalign{\medskip}p_{{1}}&-\theta_{{2}}\end {array} 
\right)+\frac{1}{z}\left(\begin {array}{cc} q_{{1}}p_{{1}}+\theta_{{1}}&p_{{3}}
\\ \noalign{\medskip}p_{{3}}q_{3}^{2}+ \left( -2q_{{1}}p_{{1}}-2
\,\theta_{{1}} \right) q_{{3}}+p_{{2}}&-q_{{1}}p_{{1}}-\theta_{{1}}
\end {array}
\right).
$$

\subsubsection{Fourth order pole. Takiff algebra of degree $3$} Here we provide only the result 
\begin{multline}
Q_1^{-1}A(z)Q_1 = \frac{\theta_4}{z^4}\left(\begin{array}{cc}
    1 & 0 \\
    0 & -1
\end{array}\right)+\frac{1}{z^3}\left(\begin {array}{cc} \theta_{{3}}&-2\theta_{{4}}q_{{3}}
\\ \noalign{\medskip}2\theta_{{4}}q_{{4}}&-\theta_{{3}}
\end {array} 
\right)+\\ + \frac{1}{z^2}\left( \begin{array}{cc}
  2\theta_{{4}}q_{{3}}q_{{4}}+\theta_{{2}
}&-\theta_{{4}}q_{3}^{3}q_{4}^{2}+ \left( \theta_{{3}}-4
\theta_{{4}} \right) q_{{4}}q_{3}^{2}-\theta_{{4}}q_{{3}}+p_{{4}
}\\ \noalign{\medskip}-\theta_{{4}}q_{3}^{2}{q_{{4}}}^{3}+
 \left( \theta_{{3}}-4\theta_{{4}} \right) q_{4}^{2}q_{{3}}+
 \left( 2\theta_{{3}}-\theta_{{4}} \right) q_{{4}}+p_{{3}}&-2
\theta_{{4}}q_{{3}}q_{{4}}-\theta_{{2}}
\end{array}\right)+\\+\frac{1}{z}\left(\begin {array}{cc} q_{{3}}p_{{3}}-q_{{4}}p_{{4}}+\theta_{{1}}
&p_{{2}}\\ \noalign{\medskip}p_{{2}}q_{2}^{2}-2\,p_{{3}}q_{{2}}q_{
{3}}+2\,p_{{4}}q_{{2}}q_{{4}}-2\,\theta_{{1}}q_{{2}}+p_{{1}}&-q_{{3}}
p_{{3}}+q_{{4}}p_{{4}}-\theta_{{1}}\end {array} 
\right)
\end{multline}

\begin{remark}
There is an interesting difference between poles of odd or even order. Indeed, when the order of pole is even $r+1 = 2k$, then the reduced phase space dimension is divisible by $4$, and we have a kind of {\it polarisation.} Indeed, for poles of order $2k$, the connection can be locally written as
$$
\frac{A_0}{z}+\dots\frac{A_{2k-1}}{z^{2k}},
$$
and the matrices $A_{k},\dots,A_{2k-1} $ form a Poisson commuting family whose dimension is half of the total dimension. Therefore they define a Lagrangian sub-manifold in the phase space. We can then assume that these matrices are parameterized by 
 $Q_0,\dots,Q_{k-1},P_k,\dots P_{2k-1}$ only. This hints at a hidden quaternionic (hyper-K\"ahler) structure. In the case of poles of odd order, we will still have that $A_{k+1},\dots,A_{2k-1} $ form a Poisson commuting family, but now the dimension of the subspace they define is not of half the dimension of the total space. In this case, we may expect an analog of Sasakian structure.
\end{remark}

\section{Isomonodromic deformations}\label{se:IDM}

Let us discuss an important consequence of Theorem \ref{th:aut0Tak}. Suppose we consider a connection on the Riemann sphere with $n+1$ poles of Poincar\'e ranks $r_1,\dots,r_n,r_\infty$ and ask about how to deform it by keeping the monodromy data constant. To answer, we have to choose some independent deformation variables and then impose that all other quantities depend on those according to the  {isomonodromicity condition}. When all poles are simple, their positions give us enough independent variables for generic isomonodromic deformations, because the number of the isomonodromic Hamiltonians equals half of the dimension of the space of accessory parameters. When higher order poles are present, their positions don't give enough independent variables. Theorem \ref{th:aut0Tak} allows us to introduce further $r-1$ independent variables for every singularity of Poincar\'e rank $r$, or in other words we have the following:

\begin{corollary}\label{cor:gen-form}
The general element in the Takiff algebra co-adjoint orbit $\widehat{\mathcal O}^\star_r$ has the  form 
\begin{equation}\label{4.3_Conn_0}
     \sum\limits_{i=0}^{r}\frac{B_i(t_1,t_2\dots t_{r})}{(\lambda-u)^{i+1}},
\end{equation}
with
$$
B_i(t_1,t_2,\dots t_{r}) = \sum\limits_{j=i}^r A_j\mathcal{M}^{(r)}_{i,j}(t_{1},t_2,\dots t_{r}),\quad 
\mathcal{M}^{(r)}_{i,j}  = \frac{1}{j!}\frac{d^j}{d\varepsilon^j} P_r(t,\varepsilon)^i\Big\vert_{\varepsilon=0},\quad P_r(t,\varepsilon) = \sum\limits_{i=1}^{r}\varepsilon^it_i,
$$
and the coefficients $A_j$ satisfy the Takiff algebra Poisson bracket (\ref{1.1}).
\end{corollary}

In this paper,  we therefore consider the isomonodoromic deformations of connections of the form
\begin{equation}\label{eq:genIso}
    \frac{d}{d\lambda}\Psi = \sum\limits_{i=0}^n\left(\sum\limits_{j=0}^{r_i}\frac{B^{(i)}_j\left(t^{(i)}_1,t^{(i)}_2\dots t^{(i)}_{r_i-1}\right)}{(\lambda-u_i)^{j+1}}-\sum_{i=1}^{r_{\infty}}\lambda^{i-1}B_{i}^{(\infty)}\left(t^{(\infty)}_1,t^{(\infty)}_2\dots t^{(\infty)}_{r_\infty-1}\right)\right)\Psi,
\end{equation}
where the deformation parameters are the locations of the poles $u_1\dots u_n$ and the coefficients of the Poisson Takiff algebra automorphisms $t_j^{(i)}$, for $i=1,\dots,n,\infty$ and $j=1,\dots,r_i-1$. The isomonodromic deformation condition means that the matrix differential one from
\begin{equation}
\Omega = \d_{u,t} \Psi \Psi^{-1} = \sum\limits_{i=1}^{n}\left[\Omega^{(0)}_i\d u_i + \sum\limits_{j=1}^{r_i-1} \Omega^{(j)}_i\d t^{(i)}_j\right],
\end{equation}
is a single valued holomorphic one form on $\mathbb{CP}^1\setminus\{u_1\dots u_n\}$.
In general, the explicit form of $\Omega$ may be obtained by studying the local solutions of the equation (\ref{eq:genIso}) as in the celebrated paper by Jimbo, Miwa and Ueno \cite{JMUI}.

In this paper we consider the isomonodromic deformations of the connections \eqref{eq:genIso} as  non-autonomous Hamiltonian systems written on a suitable set of  co-adjoint orbits. The zero curvature condition splits the isomonodromic equation into two parts: a Lax equation that defines the dynamics on the co-adjoint orbits, and  an additional relation 
between the partial derivative of $\Omega$ w.r.t. $\lambda$ and the partial derivative of the connection with respect to deformation parameters
$$
\frac{d}{d t^{(i)}_j} A - \frac{\partial}{\partial \lambda} \Omega^{(i)}_j + \left[A,\Omega^{(i)}_j \right] = \underbrace{\left(\frac{\partial}{\partial t^{(i)}_j} A -\frac{\partial}{\partial \lambda} \Omega^{(i)}_j \right)}_{0} + \underbrace{\left(\left(\frac{d}{d t^{(i)}_j}-\frac{\partial}{\partial t^{(i)}_j}\right) A + \left[A,\Omega^{(i)}_j \right]\right)}_{0} = 0.
$$
Thanks to this, we may define the coefficients of the one form $\Omega$ through the following formula:
\begin{equation}\label{eq:Omega-def}
    \Omega^{(i)}_j= \int \frac{\partial A}{\partial t^{(i)}_j}  d\lambda.
\end{equation}
The matrix $\Omega^{(i)}_j$ is defined up to the addition of a matrix which does not depend on $\lambda$. Different choices of the gauge result in different constant terms - we will see how to fix this constant term in the examples section \ref{suse:PIIJM}.

As mentioned before, the deformation parameters  $t_1^{(i)},\dots t_{r_i}^{(i)}$, $i=1,\dots,n,\infty$ appear as the result of confluence and may be seen as avatars of the Schlesinger system deformation parameters we start with. If we consider the divisor of singularities (where we denote $\infty$ by $u_{n+1}$)
$$
D:= \sum_{1}^{n+1} (r_i+1)u_i,
$$
we see that the total number of deformation parameters we introduce is given via the degree of such divisor, i.e.
$$
d = \underset{\text{\# of singularities}}{n+1} + \underset{\text{\# irregular times}}{\sum r_i}.
$$
In this paper, the idea is that the number of deformation parameters doesn't change during the confluence procedure, or, in other words, the degree $d$ is fixed. 

Here, we want to answer an important question raised by Bertola and Harnad: what is the relation between our deformation parameters and the Jimbo-Miwa-Ueno ones?
In \cite{JMUI}, the number of deformation parameters  depends on the degree of the singularity divisor {\it as well as on the rank of the connection.} The number of Jimbo-Miwa deformation parameters is not preserved during the confluence cascade. Each coalescence leads to the appearance of additional $m-2$ parameters, where $m$ is the rank of isomonodromic problem. Here we refer to the rank of a Lie algebra as the dimension of any of its Cartan subalgebras $\mathfrak{h}$. Obviously in the case of $\mathfrak{sl}_2$ connection, this number equals to zero and the number of Jimbo-Miwa-Ueno deformation parameters coincides with ours.

Let's dwell on the $\mathfrak{sl}_2$ case in more detail to explain the relation between our parameters and the ones by Jimbo-Miwa-Ueno. Consider a connection with a pole at $u$ of Poincar\'e rank $r$, i.e.
\begin{equation*}
A \underset{\lambda\to u}{\sim} \frac{B_r}{z^{r+1}}+\frac{B_{r-1}}{z^r}+\dots \frac{B_0}{z}+O(1)\quad \in \mathfrak{sl}_2,
\end{equation*}
where $z=\lambda - u$ is the local coordinate and the matrices $B_k$ are linear combinations of the bare co-adjoint orbit coordinates $A_j$ and  contain our deformation parameters as specified in formula \eqref{sect3:monom0}.

The Jimbo-Miwa-Ueno deformation parameters $w_j$ are the exponents of asymptotic behaviour of the formal solution at the irregular pole:
$$
\Psi \underset{\lambda\to u}{\sim} P(z)\left(\ID + o(z)\right)z^{\Lambda_0}\exp\left[-\sum_{j=1}^r\frac{w_j}{jz^j}\sigma_3\right],\quad \sigma_3=\left(\begin{array}{cc}
    1 & 0 \\
    0 & -1
\end{array}\right).
$$
These $w_j$ can in fact be seen as the  spectral invariants associated to the matrices $B_k$. 
Thanks to this fact, in the case of $\mathfrak{sl}_2$ there is a rational map which sends the  Jimbo-Miwa-Ueno deformation parameters to ours. To obtain this map explicitly, 
we diagonalise at the pole $\lambda= u$ and  obtain the following correspondence between Jimbo-Miwa-Ueno deformation parameters $w_i$ and our $t_j$ via
$$
\begin{array}{l}
    w_r = \theta_r t_1^r  \\
    w_{r-1} = \theta_{r-1} t_1^{r-1} +  (r-1)\theta_rt_1^{r-2}t_2 \\
    \dots \\
    w_{k} = \sum\limits_{j=k}^r \theta_j \mathcal{M}_{k,j}^{(r)}(t_1,t_2\dots,t_r) \\
    \dots \\
    w_1 = \sum\limits_{i=1}^{r}\theta_i t_i.
\end{array}
$$
Here the $\theta_i$'s are the spectral invariants of the matrices $A_j$, so we separate the non-autonomous part (dependence on deformation parameters) from the spectral invariants that determine the symplectic leaf in the phase space.
Roughly speaking, this map is a map between 2 phase spaces
$$
\hat{\mathfrak{g}}_r \,\rightarrow \, \hat{\mathcal{O}}_r\times \mathbb{C}^r,
$$
which is not bi-rational - starting from the irregular point of Poincar\'e rank 2 we have to deal with square roots if when we write $t_1\dots t_r$ via Jimbo-Miwa parameters $w_j$'s.

For  higher rank, we may think about our times as a special sub-family of the Jimbo-Miwa-Ueno isomonodromic deformations. The local solution writes as
$$
\Psi \underset{\lambda\to u}{\sim} P(z)\left(\ID + o(z)\right)z^{\Lambda_0}\exp\left[-\sum_{j=1}^r\frac{1}{jz^j}\left(\begin{array}{cccc}
    w^{(j)}_1 & 0 & \dots & 0    \\
    0 & w^{(j)}_2 & \dots & 0 \\
    &\dots  & \dots &  \\
    0 & \dots & 0 & w^{(j)}_m
\end{array}\right)\right]
$$
and $w^{(j)}_k$ are the Jimbo-Miwa-Ueno deformation parameters. Then our deformation parameters are given by the following special trajectory
$$
\frac{w^{(j)}_k}{w^{(j)}_l} = \operatorname{const},
$$
and may be considered as the deformation along a projective line in a space of Jimbo-Miwa-Ueno parameters.

In the next section we will see how the general form \eqref{eq:genIso} of the isomonodromic problem with irregular singularities naturally arises during the confluence procedure.

\section{Confluence procedure}\label{se:CONF}

\subsection{Coalescence of two simple poles}\label{suse:1p1}
Without loss of generality, we consider confluence of  $u_n:=v_1$ and $u_{n-1}:=w$, which is given by the following change of deformation parameters
\begin{equation}\label{sect2.1:1+1confl}
u_{i}=u_{i},\quad i=1\dots n-1,\quad v_1=w+\varepsilon t_1.
\end{equation}
Taking the limit $\varepsilon\rightarrow 0$ the deformation parameter  $v_1$ tends to $w$,  this is what is meant by {\it coalescence.}\/ We rewrite matrix $A(\lambda)$ as
$$
A(\lambda) = \sum\limits_{i=1}^{n-2}\frac{A^{(i)}}{\lambda-u_i} + \frac{B}{\lambda-w}+\frac{C}{\lambda - w-\varepsilon t_1},\quad B= A^{(n-1)}, \quad C=A^{(n)},
$$
where $B$ and $C$ are introduced as a convenient notation to avoid too many indices. We want to assume some $\varepsilon$ expansions for the matrices $B$ and $C$ in order that the limit of $A(\lambda) $ as $\varepsilon\mapsto 0$ is well defined and the resulting system has a double pole at $w$. To this end, observe that by rewriting the last two terms in $A(\lambda)$ as
$$
 \frac{B}{\lambda-w}+\frac{C}{\lambda - w-\varepsilon t_1}=\frac{B}{\lambda-w}+\frac{1}{\lambda - w}C\left(1-\frac{\varepsilon t_1}{\lambda - w}\right)^{-1} 
$$
and expanding $\left(1-\frac{\varepsilon t_1}{\lambda - w}\right)^{-1} $ in $\varepsilon$ we obtain
$$
 \frac{B}{\lambda-w}+\frac{C}{\lambda - w-\varepsilon t_1} =  \frac{C+B}{\lambda - w} + \frac{\varepsilon t_1}{(\lambda - w)^2}C+O(\varepsilon^2).
$$
In order to produce  a second order pole, we need the following two limits to be finite:
$$
\underset{\varepsilon\rightarrow 0}{\lim} (\varepsilon C) := A^{(n-1)}_1\neq 0,\quad \underset{\varepsilon\rightarrow 0}{\lim}( C+B) := A^{(n-1)}_0,
$$
Assuming that $A^{(i)}$'s, $B$ and $C$ may be expanded in the Laurent series in $\varepsilon$ we obtain expansions
\begin{equation}\label{sect2.1:expansion}
    A^{(i)} = \tilde{A}^{(i)} + O(\varepsilon),\quad C = \frac{1}{\varepsilon}A_1^{(n-1)}+C_0+O(\varepsilon),\quad B = -\frac{1}{\varepsilon}A_1^{(n-1)} + B_0 + O(\varepsilon),\quad C_0+B_0 = A_0^{(n-1)}.
\end{equation}
Note that we have called these limits $ A_0^{(n-1)}$ and $A_1^{(n-1)} $ respectively to adhere to the notation of section 3.

Under these hypotheses, we can take the limit as $\varepsilon\rightarrow 0$ and define
\begin{equation}\label{limA}
\tilde{A}(\lambda):=\underset{\varepsilon\rightarrow 0}{\lim} A(\lambda) = \sum\limits_{i=1}^{n-2}\frac{\tilde{A}^{(i)}}{\lambda - \tilde{u}_i} + t_1\frac{{A}^{(n-1)}_1}{(\lambda-w)^2}+\frac{{A}^{(n-1)}_0}{\lambda-w}.
\end{equation}

\begin{remark} Observe that the number of deformation parameters has not changed after the confluence, $n-1$ of them have remained as positions of poles, but one of them has become part of the leading term at the second order pole - this is compatible with Theorem \ref{th:aut0Tak}. Indeed, in the next Proposition \ref{prop:PB} we will prove that the matrices ${A}_1^{(n-1)}$ and ${{A}^{(n-1)}_0}$ satisfy the Takiff algebra Poisson brackets. We will see that as we increase the Poincar\'e rank of the poles in the confluence procedure, more and more deformation parameters will appear in the numerators of pole expansions exaclty in the way predicted by Theorem \ref{th:aut0Tak}.
\end{remark}

Now let us focus on the deformation equations. The change of variables (\ref{sect2.1:1+1confl}) transforms the deformation 1-form (\ref{sect1:defForm}) to
$$
\Omega = -\sum\limits_{i=1}^{n-2}\frac{{A}^{(i)}}{\lambda-u_i}\d u_i - \frac{A^{(n-1)}}{\lambda-w}\d w - \frac{A^{(n)}}{\lambda - w-\varepsilon t_1} (\d w+\varepsilon\d t_1). 
$$
Applying the expansion (\ref{sect2.1:expansion}), we obtain
\begin{equation}\label{limO}
\tilde{\Omega} = \underset{\varepsilon\rightarrow 0}{\lim} \Omega =  -\sum\limits_{i=1}^{n-2}\frac{\tilde{A}^{(i)}}{\lambda-u_i}\d u_i - \left(t_1\frac{{A}^{(n-1)}_1}{(\lambda-w)^2}+\frac{{A}^{(n-1)}_0}{\lambda-w}\right)\d w-\frac{{A}^{(n-1)}_1}{\lambda-w}\d t_1.
\end{equation}
The deformation 1-form $\tilde\Omega$ satisfies equation \eqref{eq:Omega-def} with $\tilde A$ in place of $A$.

\begin{definition}
We call the process of taking the expansions  \eqref{sect2.1:expansion} and the limits \eqref{limA}, \eqref{limO}, {\it 1$+$1 confluence procedure.}
\end{definition}

The the connection  $A$ and the deformation one form $\Omega$ are linear in $A^{(i)}$'s so the $O(\varepsilon)$ terms vanish during the limiting procedure. Since the Poisson structure and the Schlesinger Hamiltonians are  quadratic structures the limiting procedure becomes more complicated. Now we explain how to tackle their confluence.

\begin{proposition}\label{prop:PB}
The 1$+$1 confluence procedure gives a Poisson morphism between the direct product of the co-adjoint orbits to the Lie algebra and the co-adjoint orbit of the Takiff algebra:
$$
\mathcal{O}^{\star}_1\times\mathcal{O}^{\star}_2\times\dots \mathcal{O}^{\star}_n\times\mathcal{O}_{\infty}^{\star} \xrightarrow{\operatorname{confluence}} \mathcal{O}^{\star}_1\times\mathcal{O}^{\star}_2\times\dots \mathcal{O}^{\star}_{n-2}\times \hat{\mathcal{O}}^{\star}_{2,n-1}\times\mathcal{O}_{\infty}^{\star}.
$$
Namely, if the matrices $A^{(i)},B,C$ satisfy the \KKS brackets \eqref{sect1:FuchsPB}, then the matrices $\tilde A^{(i)},{A}^{(n-1)}_0,{A}^{(n-1)}_1$ satisfy  the Poisson algebra of the coefficients for the Takiff algebra \eqref{1.1}, i.e.
\begin{multline}\label{conf-PBT}
 \left\{\tilde{A}^{(i)}_{\alpha}, \tilde{A}^{(j)}_{\beta}\right\}= -\delta_{ij}\sum\limits_{\gamma} \chi^{\gamma}_{\alpha\beta} \tilde{A}^{(i)}_{\gamma},\quad
  \left\{\tilde{A}^{(i)}_{\alpha}, {A}^{(n-2)}_{0,\beta}\right\}= \left\{\tilde{A}^{(i)}_{\alpha}, {A}^{(n-2)}_{1,\beta}\right\}=0\quad  i,j=1,\dots n-2,\\
\left\{{A}^{(n-2)}_{1,\alpha},{A}^{(n-2)}_{1,\beta}\right\}=0,\quad \left\{{A}^{(n-2)}_{1,\alpha},{A}^{(n-2)}_{0,\beta}\right\}=-\chi_{\alpha\beta}^{\gamma}{A}^{(n-2)}_{1,\gamma},\quad \left\{{A}^{(n-2)}_{0,\alpha},{A}^{(n-2)}_{0,\beta}\right\}=-\chi_{\alpha\beta}^{\gamma}\left({A}^{(n-2)}_{0,\gamma}\right),
\end{multline}
\end{proposition}

\proof  The Poisson structure \eqref{conf-PBT} for the coefficients of the connection near the irregular singularity coincides with the standard Lie-Poisson bracket \eqref{1.1} for the co-adjoint orbit $\tilde{\mathcal{O}}_2^{\star}$ of the Takiff algebra $\mathfrak{g}_2 \backsimeq \mathfrak{g}[z]/(z^2\mathfrak{g}[z])$, where $ \mathfrak{g}[z]$ is a Lie algebra of the polynomials with coefficients in $\mathfrak{g}$. Therefore, we need to prove that \eqref{conf-PBT} arises as the $1+1$ confluence from the standard Lie-Poisson bracket on 
$\mathcal{O}^{\star}_{n-1}\times\mathcal{O}^{\star}_n$.
The first row relations are straightforward and we omit the proof. To prove the relations in the second row of \eqref{conf-PBT}, let us consider the Poisson relations  (\ref{sect1:FuchsPB})  for $B$ and $C$
$$
\{C_\alpha, C_\beta\} = -\sum\limits_{\gamma} \chi^{\gamma}_{\alpha\beta} C_{\gamma},\quad \{B_\alpha, B_\beta\} = -\sum\limits_{\gamma} \chi^{\gamma}_{\alpha\beta} B_{\gamma},\quad \{C_\alpha, B_\beta\} = 0.
$$
Inserting the expansion (\ref{sect2.1:expansion}) and
expanding the Poisson relations in $\varepsilon$, we obtain
\begin{multline*}
\frac{1}{\varepsilon^2}\left\{A^{(n-1)}_{1,\alpha},A^{(n-1)}_{1,\beta}\right\}+\frac{1}{\varepsilon}\left(\left\{A^{(n-1)}_{1,\alpha}, C_{0,\beta}\right\}+\left\{C_{0,\alpha}, A^{(n-1)}_{1,\beta}\right\}\right)+ \\ +\left\{C_{0,\alpha},C_{0,\beta}\right\} + \left\{A^{(n-1)}_{1,\alpha}, C_{1,\beta}\right\}+\left\{C_{1,\alpha}, A^{(n-1)}_{1,\beta}\right\}  = -\chi^{\gamma}_{\alpha \beta}\left(\frac{1}{\varepsilon}A^{(n-1)}_{1,\gamma} + C_{0, \gamma} \right) + o(\varepsilon)
\end{multline*}
\begin{multline*}
\frac{1}{\varepsilon^2}\left\{A^{(n-1)}_{1,\alpha},A^{(n-1)}_{1,\beta}\right\}-\frac{1}{\varepsilon}\left(\left\{A^{(n-1)}_{1,\alpha}, B_{0,\beta}\right\}+\left\{B_{0,\alpha}, A^{(n-1)}_{1,\beta}\right\}\right)+ \\ +\left\{B_{0,\alpha},B_{0,\beta}\right\} - \left\{A^{(n-1)}_{1,\alpha}, B_{1,\beta}\right\} - \left\{B_{1,\alpha}, A^{(n-1)}_{1,\beta}\right\}  = \chi^{\gamma}_{\alpha \beta}\left(\frac{1}{\varepsilon}A^{(n-1)}_{1,\gamma} - B_{0, \gamma} \right) +o(\varepsilon)
\end{multline*}
\begin{multline*}
-\frac{1}{\varepsilon^2}\left\{A^{(n-1)}_{1,\alpha},A^{(n-1)}_{1,\beta}\right\}+\frac{1}{\varepsilon}\left(\left\{A^{(n-1)}_{1,\alpha}, B_{0,\beta}\right\}-\left\{C_{0,\alpha}, A^{(n-1)}_{1,\beta}\right\}\right)  +  \\ +\left\{C_{0,\alpha},B_{0,\beta}\right\} + \left\{A^{(n-1)}_{1,\alpha}, B_{1,\beta}\right\} - \left\{C_{1,\alpha}, A^{(n-1)}_{1,\beta}\right\}=  o(\varepsilon).
\end{multline*}
Collecting different terms in $\varepsilon$, we obtain
\begin{eqnarray}\label{ConfPA}
\varepsilon^{-2}:&& \left\{A^{(n-1)}_{1,\alpha},A^{(n-1)}_{1,\beta}\right\}=0, \nn \\
\varepsilon^{-1}:&& \left\{A^{(n-1)}_{1,\alpha}, C_{0,\beta}\right\}+\left\{C_{0,\alpha}, A^{(n-1)}_{1,\beta}\right\} =  -\chi^{\gamma}_{\alpha \beta}A^{(n-1)}_{1,\gamma},\nn\\
\varepsilon^{-1}: && \left\{A^{(n-1)}_{1,\alpha}, B_{0,\beta}\right\}+\left\{B_{0,\alpha}, A^{(n-1)}_{1,\beta}\right\} = -\chi^{\gamma}_{\alpha \beta}A^{(n-1)}_{1,\gamma},\nn\\
\varepsilon^{-1}: && \left\{A^{(n-1)}_{1,\alpha}, B_{0,\beta}\right\}-\left\{C_{0,\alpha}, A^{(n-1)}_{1,\beta}\right\}=0, \\
\varepsilon^{0}:&&\left\{C_{0,\alpha},C_{0,\beta}\right\} + \left\{A^{(n-1)}_{1,\alpha}, C_{1,\beta}\right\}+\left\{C_{1,\alpha}, A^{(n-1)}_{1,\beta}\right\} = -\chi^{\gamma}_{\alpha \beta}C_{0, \gamma}, \nn \\
\varepsilon^{0}:&& \left\{B_{0,\alpha},B_{0,\beta}\right\} - \left\{A^{(n-1)}_{1,\alpha}, B_{1,\beta}\right\} - \left\{B_{1,\alpha}, A^{(n-1)}_{1,\beta}\right\} = -\chi^{\gamma}_{\alpha \beta} B_{0, \gamma}  \nn \\
\varepsilon^{0}:&& \left\{C_{0,\alpha},B_{0,\beta}\right\} + \left\{A^{(n-1)}_{1,\alpha}, B_{1,\beta}\right\} - \left\{C_{1,\alpha}, A^{(n-1)}_{1,\beta}\right\}=0. \nn
\end{eqnarray}
The term of order  $\varepsilon^{-2}$ in (\ref{ConfPA}) proves the first relation in the second row of \eqref{conf-PBT}. Let us prove the second relation. Take the $1/\varepsilon$ term
$$
\left\{{A}^{(n-2)}_{1,\alpha}, B_{0,\beta}\right\}-\left\{C_{0,\alpha}, {A}^{(n-2)}_{1,\beta}\right\}=0 \quad \Longleftrightarrow \quad  \left\{C_{0,\alpha},
{A}^{(n-2)}_{1,\beta}\right\}=\left\{ {A}^{(n-2)}_{1,\alpha}, B_{0,\beta}\right\}
$$
and put it in the Poisson relation between $ {A}^{(n-1)}_1$ and $C_0$. We get
\begin{multline*}
-\chi_{\alpha\beta}^{\gamma} {A}^{(n-2)}_{1,\gamma}= \left\{ {A}^{(n-2)}_{1,\alpha}, C_{0,\beta}\right\}+\left\{C_{0,\alpha}, {A}^{(n-2)}_{1,\beta}\right\} 
= \left\{ {A}^{(n-2)}_{1,\alpha}, C_{0,\beta}\right\}+\left\{{A}^{(n-2)}_{1,\alpha}, B_{0,\beta}\right\} = \\ = \boxed{\left\{{A}^{(n-2)}_{1,\alpha}, C_{0,\beta} + B_{0,\beta}\right\} = -\chi_{\alpha\beta}^{\gamma}{A}^{(n-2)}_{1,\gamma}}
\end{multline*}
which proves the second relation. Now let us compute the last Poisson bracket
$$
\left\{C_{0,\alpha}+B_{0,\alpha},C_{0,\beta}+B_{0,\beta}\right\} = \left\{C_{0,\alpha},C_{0,\beta}\right\}+\left\{C_{0,\alpha},B_{0,\beta}\right\}+\left\{B_{0,\alpha},C_{0,\beta}\right\}+\left\{B_{0,\alpha},B_{0,\beta}\right\}.
$$
Using the $\varepsilon^{0}$-terms from (\ref{ConfPA}) for $\left\{C_{0,\alpha},C_{0,\beta}\right\}$ and $\left\{B_{0,\alpha},B_{0,\beta}\right\}$, we obtain
\begin{multline}\label{cbPB}
\left\{C_{0,\alpha}+B_{0,\alpha},C_{0,\beta}+B_{0,\beta}\right\} = -\chi_{\alpha\beta}^{\gamma}(C_{0,\beta}+B_{0,\beta}) - \left\{{A}^{(n-2)}_{1,\alpha},C_{1,\beta}\right\}-\left\{C_{1,\alpha},{A}^{(n-2)}_{1,\beta}\right\}+\\+\left\{{A}^{(n-2)}_{1,\alpha},B_{1,\beta}\right\}+\left\{B_{1,\alpha},{A}^{(n-2)}_{1,\beta}\right\}
+\left\{C_{0,\alpha},B_{0,\beta}\right\}+\left\{B_{0,\alpha},C_{0,\beta}\right\}.
\end{multline}
The last $\varepsilon^{0}$-term in (\ref{ConfPA}) leads to the following relations
$$
\left\{C_{0,\alpha},B_{0,\beta}\right\} = \left\{C_{1,\alpha}, {A}^{(n-2)}_{1,\beta}\right\}-\left\{ {A}^{(n-2)}_{1,\alpha},B_{1,\beta}\right\}
$$
$$
\left\{B_{0,\alpha},C_{0,\beta}\right\} = \left\{ {A}^{(n-2)}_{1,\alpha},C_{1,\beta}\right\}-\left\{B_{1,\alpha}, {A}^{(n-2)}_{1,\beta}\right\}
$$
which cancel all terms in the right-hand side of (\ref{cbPB}) except the first term, so we obtain
$$
\left\{C_{0,\alpha}+B_{0,\alpha},C_{0,\beta}+B_{0,\beta}\right\} = -\chi_{\alpha\beta}^{\gamma}(C_{0,\gamma}+B_{0,\gamma}),
$$
which concludes the proof.
 \endproof

Observe that the relations \eqref{ConfPA} contain more information than we need, and that one could actually try to come up with a Poisson algebra involving all coefficients $B_k$, $C_k$ in the expansion \eqref{sect2.1:expansion}. However we are only interested in the Poisson subalgebra generated by  ${A}^{(n-1)}_1$, ${A}^{(n-1)}_0=C_0+B_0$ and $\tilde A^{(i)}$ for $i=1,\dots,n-2$. The main feature of this subalgebra is that it does not depend on a choice of a Poisson algebra for the coefficients $B_k$ and $C_k$ for $k>1$. We call this sub-algebra {\it Isomonodromic Poisson Algebra (IPA),} since these are the only elements which survive in the isomonodromic problem after the confluence procedure.

\begin{proposition}
The $1+1$ confluence procedure produces the isomonodromic Hamiltonians giving the zero curvature condition
$$
{\rm d}_u \tilde A -\frac{d}{d\lambda}\tilde\Omega+[\tilde A,\tilde\Omega]=0
$$
as equation of motion. 
\end{proposition}

\proof To prove this, we start from the extended symplectic form for the Schlesinger equations:
$$
\omega_{\hbox{\tiny{KKS}}}+ \sum\limits_{i=1}^n\d u_i \wedge \d H_i.
$$
Here $\omega_{\hbox{\tiny{KKS}}}$ is the symplectic form which corresponds to the \KKS structure on the direct product of the co-adjoint orbits. Thanks to Proposition \ref{prop:PB}, the \KKS bracket tends to the Takiff algebra Poisson bracket, therefore $\omega_{\hbox{\tiny{KKS}}}$ tends to the corresponding symplectic form. Let us concentrate on the $ \sum\limits_{i=1}^n\d u_i \wedge \d H_i$ part. This part transforms to
$$
\sum\limits_{i=1}^n\d u_i \wedge \d H_i\rightarrow \sum\limits_{i=1}^{n-2}\d u_i \wedge \d H_i + \d w\wedge \d \left(H_{n-1}+H_{n}\right) + \d t_1 \wedge \d\left( \varepsilon H_n\right).
$$
Since we are working on a symplectic leaf of the \KKS bracket, the central elements, or Casimirs, can be considered as fixed scalars, i.e. the differential $\d$ acts on them as a zero. To find the Hamiltonians of the confluent dynamic we have to calculate the limit of the ``time-dependent" part of the symplectic structure as $\varepsilon$ goes to zero. In other words, we have to find 
\begin{equation}\label{eq:Lie-conf}
\d \tilde{H}_i := \underset{\varepsilon\rightarrow 0}{\lim}\, \d H_i,\quad \d \tilde{H}_{n-1} := \underset{\varepsilon\rightarrow 0}{\lim}\, \d  (H_{n-1}+H_{n}),\quad \d \tilde{H}_{n} := \underset{\varepsilon\rightarrow 0}{\lim}\, \varepsilon \d H_{n}.
\end{equation}
To compute these limits, we can treat the Hamiltonians up to addition of Casimirs. This allows us to use the Casimirs to regularise parts of the Hamiltonains that are singular in $\varepsilon$. Therefore all $=$ signs in the rest of the proof are intended as equal up to Casimirs. For $i<n-2$ we have
\begin{equation}
\tilde{H}_i := \underset{\varepsilon\rightarrow 0}{\lim}H_i = \sum\limits_{j\neq i}^{n-2}\frac{\tr\left(\tilde{A}^{(i)}\tilde{A}^{(j)}\right)}{u_i-u_j} + t_1\frac{\tr\left(\tilde{A}^{(n-1)}_1 \tilde{A}^{(i)}\right)}{(u_i-w)^2}+\frac{\tr\left(\tilde{A}^{(n-1)}_0 \tilde{A}^{(i)}\right)}{u_i-w},
\end{equation}
for $i=n-1$ we have
\begin{multline}
\tilde{H}_{n-1} = \underset{\varepsilon\rightarrow 0}{\lim} \left(H_{n-1}+H_n\right) = \underset{\varepsilon\rightarrow 0}{\lim}\sum\limits_{j<n-2}\tr \tilde{A}^{(j)}\left(\frac{A^{(n-1)}}{w-u_{j}}+\frac{A^{(n)}}{w+\varepsilon t_1-u_{j}}\right) = \\ = \sum\limits_{j<n-1}\tr \tilde{A}^{(j)}\left(\frac{\tilde{A}^{(n-1)}_0}{w-u_{j}}-t_1\frac{\tilde{A}^{(n-1)}_1}{(w-u_{j})^2}\right).
\end{multline}
For $i=n$
$$
\tilde{H}_{n} = \underset{\varepsilon\rightarrow 0}{\lim} \varepsilon H_n.
$$
Substituting coalescence expansions we get
\begin{equation}\label{sect2:eHn}
\varepsilon H_n = \left[\sum\limits_{j<n-2}\frac{\tr \tilde{A}^{(j)}{A}^{(n-2)}_1}{w-u_{j}}+O(\varepsilon)\right] + \frac{\tr A^{(n)}A^{(n-1)}}{t_1}.
\end{equation}
The last term in (\ref{sect2:eHn}) contains terms of order $1/\varepsilon$ and $1/\varepsilon^2$. :
\begin{multline*}
 \frac{\tr A^{(n)}A^{(n-1)}}{\tilde{u}_n} = \frac{1}{\tilde{u}_n}\left(-\frac{1}{\varepsilon^2}\tr\left(\tilde{A}^{(n-1)}_1\right)^2+\frac{1}{\varepsilon}\tr\left(\tilde{A}^{(n-1)}_1B_0-C_0\tilde{A}^{(n-1)}_1\right)+\tr(B_0C_0) \right) +\\ +  \frac{1}{\tilde{u}_n}\tr\left(\tilde{A}^{(n-1)}_1B_1-C_1\tilde{A}^{(n-1)}_1\right)
\end{multline*}
The $1/\varepsilon^2$ term is a Casimir of the Poisson structure, so we may drop it.

Let us show that also the $1/\varepsilon$-term is a Casimir and that, after eliminating the Casimirs, $\epsilon H_n\to \tilde H_n+O(\varepsilon)$ where 
\begin{equation}\label{sect2:Hn}
\tilde{H}_{n} = \sum\limits_{j<n-2}\frac{\tr \tilde{A}^{(j)}{A}^{(n-2)}_1}{w-u_{j}} + \frac{1}{t_1}\frac{\tr\left(\tilde{A}^{(n-1)}_0\right)^2}{2}.
\end{equation}
To see this, let us remind that the Casimirs of the Poisson algebra in the Fuchsian case are $\tr \left(A^{(i)}\right)^k$, so the function
$$
\frac{1}{2}\tr \left(A^{(n)}+A^{(n-1)}\right)^2
$$
differs from  the last term of (\ref{sect2:eHn})
$$
\tr A^{(n)}A^{(n-1)}.
$$
by a Casimir. Since the Hamiltonians are defined up to the addition of a Casimir, we obtain
$$
\varepsilon H_n = \sum\limits_{j<n-2}\frac{\tr \tilde{A}^{(j)}{A}^{(n-2)}_1}{w-u_{j}} + \frac{\tr \left(A^{(n)}+A^{(n-1)}\right)^2}{2t_1} +O(\varepsilon)=  \sum\limits_{j<n-2}\frac{\tr \tilde{A}^{(j)}{A}^{(n-2)}_1}{w-u_{j}} + \frac{1}{t_1}\frac{\tr\left(\tilde{A}^{(n-1)}_0\right)^2}{2} +O(\varepsilon).
$$
Taking the limit as $\varepsilon\rightarrow 0$ we obtain the Hamiltonian (\ref{sect2:Hn}).
\endproof

\subsection{Irregular singularities arising as confluence cascades.}

In this section we consider an irregular singularity of arbitrary Poincar\'e rank $r$ as the result of the confluence cascade of $r$ simple poles $v_1,v_2\dots v_r$ with some chosen simple pole $u$. At the first step, we send $v_1$ to $u$ and create second order pole as in the previous subsection. Then we do the same for $v_2$ - we collide it with the second order pole at $u$ and create a pole of order $3$. In such a way, we continue this procedure, so at the $l$-th step we collide the simple pole $v_l$ with the pole of order $l$ at $u$ to create a new pole of order $l+1$. Finally, a the final $r$-th step, we obtain a pole of order $r+1$, i.e. of Poincar\'e rank $r$.  During this procedure, the poles $v_l$'s that disappear give rise to deformation parameters $t_l$'s for the irregular isomonodromic  problem\footnote{
The confluence procedure is not symmetric in $v_i$,  however different choices of the order of the coalescence cascade will lead to the action of the permutation group on the $t_l$'s, so we fix this ordering once for all.}. These deformation parameters appear explicitly in the coefficients of the local expansion of the connection near the singularity $u$. In the sub-section \ref{suse:kp1conf}, we prove Theorem \ref{sect4:small_theorem} that tells us that this dependence is the one described in Corollary \ref{cor:gen-form}. Before attacking that proof, let us formalise the definition of confluence:

\begin{definition} The  limiting procedure described in the hypotheses of Theorem \ref{sect4:small_theorem} is called {\it $r+1$ confluence.} 
\end{definition}

Observe that as a result of the $1+1$ confluence in  subsection \ref{suse:1p1} we obtain a connection of the form \eqref{4.3_Conn} with $r=2$. We can then apply the $1+2$ confluence to this and again obtain  a connection of the form \eqref{4.3_Conn} with $r=3$ and so on. Therefore we can give the following recursive definition:

\begin{definition} The procedure of applying Theorem \ref{sect4:small_theorem} recursively $r$ times is called 
{\it confluence cascade of $r+1$ simple poles}  on the Riemann sphere. 
\end{definition}

As mentioned at the beginning of this section, the inductive hypothesis on the local form of the connection \eqref{4.3_Conn} is not restrictive. Indeed, we expect the local form of a connection with a pole of order $r$ at $u$ to be given by an element in the Takiff algebra co-adjoint orbit $\widehat{\mathcal O}_r^\star$ with some spectral parameter $z=\lambda-u$. However, if we want to keep the number of independent variables to be maximal, we need to introduce some extra variables $t_i$ by hand in such a way that they can be treated as independent variables. In Corollary \ref{cor:gen-form}, we proved that the only way to do this is by taking precisely the form  \eqref{4.3_Conn}. Therefore, Theorem \ref{sect4:small_theorem} is equivalent to the following result:

\begin{theorem}\label{sect3:MasterTheorem}
Assume that $u$ is a singularity of Poincar\'e rank $r$ obtained by the confluence of $r+1$ simple poles. Then the coefficients of the local expansion  in a disk around $u$
$$
A(\lambda)\sim \sum\limits_{i=0}^{r}\frac{B_i(t_1,\dots t_{{r}})}{(\lambda-u)^{{i+1}}}+ O((\lambda-u)^0),
$$
take the form
$$
B_i(t_1,t_2,\dots t_r) = \sum\limits_{j=i}^{r} B^{[j]}\mathcal{M}^{(r)}_{{i},j}(t_{1},t_2,\dots t_{r}),
$$
where
$$
\mathcal{M}^{(r)}_{k,j} = \sum\limits_{w(\alpha)=j}^{|\alpha|=k}\frac{k!}{\alpha_1!\alpha_2!\dots\alpha_r!}\left(\prod\limits_{i=1}^{r}t_i^{\alpha_i}\right),\quad |\alpha| = \sum\limits_{i=1}^r \alpha_i,\quad w(\alpha) = \sum\limits_{i=1}^r i \cdot \alpha_i, \quad \mathcal{M}^{(r)}_{k>j} :=0,\quad 
$$
and $B_i^{[j]}$'s hold the following Poisson relations
\begin{equation}\label{doublePoisson}
\left\{\left(B^{[k]}\right)_\alpha, \left(B^{[p]}\right)_\beta \right\} = \left\{\begin{array}{ll}
    -\chi^{\gamma}_{\alpha\beta} \left(B^{[k+p]}\right)_{\gamma} &  \,k+p\leq r\\
    0 &  k+p> r,
\end{array}\right.
\end{equation}
 where $\chi_{\alpha\beta}^{\gamma}$ are the structure constants of the corresponding Lie algebra.
\end{theorem}


We want to underline here that the Poisson structure (\ref{doublePoisson})  gives rise to the Takiff co-algebra Poisson structure on the coefficients of the local expansion, i.e
\begin{equation}\label{sect3:TakiffT}
\left\{B_{i}(t_1,\dots t_r)_\alpha,B_{j}(t_1,\dots t_r)_\beta \right\} = - \sum\limits_{\gamma}\chi_{\alpha\beta}^{\gamma}B_{i+j}(t_1,\dots t_r)_\gamma.
\end{equation}
However, in formula \eqref{sect3:TakiffT} the dependence on the deformation parameters is implicit, while (\ref{doublePoisson}) contains information about the explicit dependence on the variables $t_i$'s.

To motivate the constructions which appear in the statement of this theorem we introduce some preliminaries on the confluence procedure and algebraic structures which appear during coalescence before the proof.
 
\subsubsection{The algebra of the weighted monomials and associated polynomials.}\label{suse:kp1conf}
The aim of this subsection is to collect some useful algebraic relations involving the coefficients $t_1,\dots,t_r$ that  arise during the confluence procedure and describe the general elements of the Takiff co-algebra with respect to the Poisson automorphisms.

In order to prove Theorem \ref{sect4:small_theorem}, in the neighborhood of a simple pole $v_r$ with a pole $w$ of Poincar\'e rank $r$ we take the following expansion
$$
v_r = w +  P_r(t,\varepsilon) = w + \sum\limits_{i=1}^rt_i\varepsilon^i.
$$
The powers of the polynomial $P_r(t,\varepsilon)$ play a significant role since they appear in the following expansion
$$
\frac{C}{\lambda - v_r} = \frac{C}{\lambda - w}\left(1 - \frac{P_r(t,\varepsilon)}{\lambda-w}\right)^{-1} = \frac{C}{\lambda - w}\left(1+\frac{P_r(t,\varepsilon)}{\lambda-w}+\dots+\frac{P_r(t,\varepsilon)^j}{(\lambda-w)^j} + \dots + \frac{P_r(t,\varepsilon)^r}{(\lambda-w)^r}\right) + O(\varepsilon^{r+1}).
$$
Each power of $P_r(t,\varepsilon)$  may be seen as a polynomial in $\varepsilon$ with coefficients in $\mathbb{C}[t_1,t_2\dots t_r]$
$$
P_r(t,\varepsilon)^r = t_1^r\varepsilon^r +O(\varepsilon^{r+1}).
$$
Because the aim of the confluence is to create a pole of Poincer\'e rank $r+1$, we need the coefficeints $(\lambda-w)^{-r-2}$ to survive, therefore,  we have to require $C$ to be a Laurent polynomial in $\varepsilon$ starting from the power $-r$.
Taking this fact into account, it is important to understand how each power of $P_r(t,\varepsilon)$ expands via $\varepsilon$
$$
P_r(t,\varepsilon)^i \mod \varepsilon^{r+1}= \sum\limits_{j=i}^r\mathcal{M}^{(r)}_{i,k}(t_1,\dots t_{r})\varepsilon^k = \mathcal{M}^{(r)}_{i,i}(t_1,\dots t_{r})\varepsilon^i+\dots + \mathcal{M}^{(r)}_{i,r}(t_1,\dots t_{r})\varepsilon^r,
$$
where
$$
\mathcal{M}^{(r)}_{i,k}(t_1,\dots t_{r}) = \frac{1}{k!}\frac{d}{d\varepsilon}P_r(t,\varepsilon)^i\Big\vert_{\varepsilon=0}.
$$
The following simple Lemma calculates an explicit formula for $\mathcal{M}^{(r)}_{i,j}(t_1,\dots t_{r})$ and gives some useful identites.

\begin{lemma}
$\mathcal{M}^{(r)}_{i,k}$ is a homogeneous polynomial in $t_1\dots t_r$ of degree $i$ for any $k$ given by 
$$
\mathcal{M}^{(r)}_{i,k}(t_1,\dots t_{r})= \sum\limits_{w(\alpha)=k}^{|\alpha|=i}\frac{i!}{\alpha_1!\alpha_2!\dots \alpha_r !} \left(\prod\limits_{j=1}^r t_j^{\alpha_j} \right),\quad w(\alpha) = \sum\limits_{j=1}^rj\alpha_j,\quad |\alpha|=\sum\limits_{j=1}^r\alpha_j.
$$
The polynomials $\mathcal{M}^{(r)}_{i,k}$ satisfy the following identities 
\begin{equation}\label{eq:mod_mon}
\mathcal{M}^{(r+1)}_{i,k} = \mathcal{M}^{(r)}_{i,k},\quad \forall k\leq r, \qquad \mathcal{M}^{(r+1)}_{i,r+1} = t_{r+1}^i
\end{equation}
and
\begin{equation}\label{eq:mod_mon1}
\mathcal{M}^{(r)}_{j,k} = \sum\limits_{p=0}^{k}\left[\mathcal{M}^{(r)}_{j-1,p}\cdot\mathcal{M}^{(r)}_{i,k-p}+\mathcal{M}^{(r)}_{j-i,k-p}\cdot\mathcal{M}^{(r)}_{i,p}\right],\quad \forall i\leq j.
\end{equation}
\end{lemma}

Note that the function $w(\alpha)$, that we call {\it weight,} can be calculated by the following formula
\begin{equation}
w\left(\prod\limits_{i=1}^n t_i^{\alpha_i}\right) = (\alpha_1,\alpha_2,\dots\alpha_n)\left(\begin{array}{c}
    1  \\
    2  \\
    .. \\
    n
\end{array}\right) = \sum\limits_{k=1}^n k\alpha_k.
\end{equation}
The weights are elements in the semi-group of homomorphism from the semi-group of monomials in the variables $t_1,\dots t_r$ to the $(\mathbb{Z}_{\geq 0}, +)$, in fact:
$$
w(\theta \cdot \eta) = w(\theta) + w(\eta).
$$

\begin{remark}
Instead of considering  the  polynomials $P_r(t_1,\dots t_r)$, we might consider the formal power series 
$$
P^{(\infty)}(\varepsilon,t) = \sum\limits_{i=1}^{\infty}t_i \varepsilon^i,
$$
and truncate all expansions at $\varepsilon^{r+1}$. The result will be the same,  but such approach probably clarifies the recursive nature of the confluence procedure. In similar way, the upper triangular matrix $\mathcal{M}^{(r)}$ with entries $\mathcal{M}^{(r)}_{i,k}$  given in (\ref{sect3:monom}) can be considered as
as a sub-matrix of size $r\times r$ in the upper left corner of some infinite dimensional upper triangular ``master" matrix $\mathcal{M}^{(\infty)}$ with entries given by
$$
\mathcal{M}^{(\infty)}_{j,r} = \frac{1}{r!}\frac{d^r}{d\varepsilon^r} P^{(\infty)}(\varepsilon, t)^j\Big\vert_{\varepsilon=0}.
$$
\end{remark}

\subsubsection{Proof of the Theorems \ref{sect3:MasterTheorem} and \ref{sect4:small_theorem}}
We use induction here to prove the theorem. Here we will start with the proof of the explicit dependence of the local expansion on $t_i$'s and then we will handle the Poisson structure.

The statement of the theorem is true for  $r=0,1$, i.e. the Fuchsian case and the case of a pole of order 2. This was proven in subsection \ref{suse:1p1}. Now let the statement be true for the irregular singularity of  Poincar\'e rank $r-1$. Adding another simple pole $v_r$, we consider the following connection
$$
A = \sum\limits_{i=0}^{r-1} \frac{B_i(t_1,t_2\dots t_{r-1})}{(\lambda- w)^{i+1}}+\frac{C}{\lambda - v_r}+\dots,\quad B_i = \sum\limits_{j=i}^{r-1} B^{[j]} \mathcal{M}_{i,j}^{(r-1)},
$$
where the dots denote regular terms in $(\lambda- w)$ and $(\lambda- v_1)$,
with the following behaviour of $v_r$
$$
v_r = w + \sum\limits_{j=1}^{r}t_j\varepsilon^j,\quad C = \sum\limits_{j=-r}^{\infty}C^{[i]}\varepsilon^i.
$$
Expanding $A$ with respect to $\varepsilon$ at $r$'th order we obtain
$$
A = \frac{C^{[-r]}t_1^{r}}{(\lambda-w)^{r+1}} + \sum\limits_{i=0}^{r-1}\frac{B_i+C P_r(t,\varepsilon)^{i}}{(\lambda-w)^{i+1}}+\dots.
$$
Using the formula \eqref{eq:mod_mon}, the coefficients $B_i$ expand via polynomials $\mathcal{M}^{(r)}_{i,j}$ giving the following
$$
B_i+C P_r(t,\varepsilon)^{i} = C^{[-r]}\mathcal{M}_{i,r}^{(r)} + \sum_{j=i}^{r-1}\left(B^{[j]}+\varepsilon^jC\right)\mathcal{M}_{i,j}^{(r)}.
$$
Since the confluence procedure requires the existence of the limit $\varepsilon\rightarrow 0$, the negative powers of $\varepsilon$ should vanish, so we obtain the expansions for the coefficients $A^{(r)}_j$ in the form
$$
B^{[k]}= -\sum\limits_{m=-r}^{-(k+1)}\frac{C^{[m]}}{\varepsilon^{m+k}}+B^{[k,0]}+\sum\limits_{l=1}^{\infty}B^{[k,l]}\varepsilon^l.
$$
Using these expansions and taking the $\varepsilon\rightarrow 0$ limit, we obtain
$$
A = \sum\limits_{i=1}^{r+1}\frac{\tilde{B}_i(t_1\dots t_r)}{(\lambda-u)^i}+\dots,
$$
where $\tilde{B}_i$'s are given by \eqref{sect3:IPA}, which finishes the proof of the first part of the theorem.

Now we prove that the Poisson structure for the coefficient of the local expansion of the connection near an irregular singularity, which arises after confluence procedure is the Poisson structure given by (\ref{doublePoisson}).

Again we use induction. The statement is obvious in case when $r=0$, and we have previously proved that it holds for $r=1$. Using the previous results  we consider the same coalescence
$$
A =  \sum\limits_{i=0}^{r-1} \frac{B_i(t_1,t_2\dots t_{r-1})}{(\lambda- w)^{i+1}}+\frac{C}{\lambda - v_r}+\dots=\sum\limits_{i=0}^{r} \frac{\tilde{B}_i}{(\lambda- w)^{i+1}}
$$
where $\tilde{B}_i$ is given by (\ref{sect3:IPA}). The expansions take the same form
$$
C= \sum\limits_{j=-r}^{\infty}C^{[i]}\varepsilon^i,\quad
B^{[k]}= -\sum\limits_{m=-r}^{-(k+1)}\frac{C^{[m]}}{\varepsilon^{m+k}}+B^{[k,0]}+\sum\limits_{l=1}^{\infty}B^{[k,l]}\varepsilon^l.
$$
In order to get rid of the indices let us use the following notation
$$
V = B^{[k]},\quad W=B^{[l]},\quad U = B^{[k+l]}.
$$
In the case when the indices on the right hand sides are out of bound we assume that the values are zero. The Poisson relations are
$$
\{V_{\alpha},W_{\beta}\} = -\chi^{\gamma}_{\alpha\beta} U_{\gamma}, \quad \{C_{\alpha},C_{\beta}\}=-\chi^{\gamma}_{\alpha\beta} C_{\gamma},\quad \{C_{\alpha},V_{\beta}\}=\{C_{\alpha},W_{\beta}\}=\{C_{\alpha},U_{\beta}\}=0
$$
and the expansions are
\begin{eqnarray}
&& C  = \frac{C^{[-r]}}{\varepsilon^r}+\frac{C^{[-r+1]}}{\varepsilon^{r-1}}+\dots +\frac{C^{[-1]}}{\varepsilon} + C^{[0]} + \sum\limits_{i=1}^{\infty}C^{[i]}\varepsilon^i  \nn\\
&& V  =-\frac{C^{[-r]}}{\varepsilon^{r-k}}-\frac{C^{[-r+1]}}{\varepsilon^{r-k-1}}-\dots -\frac{C^{[-k+1]}}{\varepsilon} + V^{[0]}+\sum\limits_{i=1}^{\infty}V^{[i]}\varepsilon^i \dots  \nn \\
&&  W  =-\frac{C^{[-r]}}{\varepsilon^{r-l}}-\frac{C^{[-r+1]}}{\varepsilon^{r-l-1}}-\dots -\frac{C^{[-l+1]}}{\varepsilon} + W^{[0]}+\sum\limits_{i=1}^{\infty}W^{[i]}\varepsilon^i \dots  \nn \\
&& U  =-\frac{C^{[-r]}}{\varepsilon^{r-k-l}}-\frac{C^{[-r+1]}}{\varepsilon^{r-k-l-1}}-\dots -\frac{C^{[-k-l+1]}}{\varepsilon} + U^{[0]}+\sum\limits_{i=1}^{\infty}U^{[i]}\varepsilon^i. \nn
\end{eqnarray}
Due to the confluence formula we want to prove that the following Poisson relation holds
\begin{multline*}
\left\{V^{[0]}_\alpha+C^{[-k]}_\alpha,W^{[0]}_\beta+C^{[-l]}_\beta \right\} = \left\{V^{[0]}_{\alpha},W^{[0]}_{\beta}\right\} + \left\{V_{\alpha}^{[0]}, C_{\beta}^{[-l]}\right\} + \left\{C_{\alpha}^{[-k]},W_{\beta}^{[0]}\right\} + \left\{C^{[-k]}_{\alpha}, C^{[-l]}_{\beta}\right\} = \\ =  -\chi^{\gamma}_{\alpha\beta}\left(U^{[0]}_{\gamma}+C^{[-k-l]}_{\gamma}\right).
\end{multline*}
Taking the corresponding $\varepsilon$-terms of the expansions of the Poisson relations we get
\begin{eqnarray}
\label{AB}
\underset{\varepsilon=0}{\operatorname{Res}}\,\varepsilon^{0-1} \left\{V_{\alpha},W_{\beta}\right\}: && 
\left\{V^{[0]}_{\alpha},W^{[0]}_{\beta}\right\}  -\mathunderline{blue}{\sum\limits_{i=1}^{r-l}\left\{V^{[i]}_{\alpha},C^{[-i-l]}_{\beta}\right\}} -\mathunderline{red}{\sum\limits_{i=1}^{r-k}\left\{C_{\alpha}^{[-i-k]},W^{[i]}_{\beta}\right\}}= 
-\chi^{\gamma}_{\alpha\beta} U^{[0]}_{\gamma}, \nn \\
\label{AW}
\underset{\varepsilon=0}{\operatorname{Res}}\,\varepsilon^{l-1}\left\{V_{\alpha},C_{\beta}\right\} : &&
\left\{V_{\alpha}^{[0]}, C_{\beta}^{[-l]}\right\} + 
\mathunderline{blue}{\sum\limits_{i=1}^{r-l}\left\{V_{\alpha}^{[i]}, C_{\beta}^{[-i-l]}\right\}} - 
\mathunderline{green}{\sum\limits_{i=1}^{r-k}\left\{C_{\alpha}^{[-i-k]},C^{[i-l]}_{\beta}\right\}} = 0, \nn \\
\label{WB}
\underset{\varepsilon=0}{\operatorname{Res}}\,\varepsilon^{k-1}\left\{C_{\alpha},W_{\beta}\right\} : &&
\left\{C_{\alpha}^{[-k]},W_{\beta}^{[0]}\right\}+
\mathunderline{red}{\sum\limits_{i=1}^{r-k}\left\{C_{\alpha}^{[-i-k]},W_{\beta}^{[i]} \right\}} - 
\mathunderline{black}{\sum\limits_{i=1}^{r-l}\left\{C_{\alpha}^{[i-k]},C^{[-i-l]}_{\beta}\right\}}=0  \\
\label{WW}
\underset{\varepsilon=0}{\operatorname{Res}}\,\varepsilon^{k+l}\frac{\left\{C_{\alpha},C_{\beta}\right\}}{\varepsilon} : && 
\left\{C^{[-k]}_{\alpha}, C^{[-l]}_{\beta}\right\} + \mathunderline{green}{\sum\limits_{i=1}^{r-k}\left\{C^{[-k-i]}_{\alpha}, C^{[i-l]}_{\beta}\right\}} + \mathunderline{black}{\sum\limits_{i=1}^{r-l}\left\{C^{[-k+i]}_{\alpha}, C^{[-i-l]}_{\beta}\right\}} = \nn \\ && =-\chi^{\gamma}_{\alpha\beta} C^{[-l-k]}_\gamma. \nn
\end{eqnarray}
Finally, taking the sum of the relations in (\ref{WB}) we get the desired Poisson structure
$$
\left\{V^{[0]}_{\alpha},W^{[0]}_{\beta}\right\} + \left\{V_{\alpha}^{[0]}, C_{\beta}^{[-l]}\right\} + \left\{C_{\alpha}^{[-k]},W_{\beta}^{[0]}\right\} + \left\{C^{[-k]}_{\alpha}, C^{[-l]}_{\beta}\right\}  = -\chi^{\gamma}_{\alpha\beta}\left(U^{[0]}_{\gamma}+ C^{[-l-k]}_\gamma\right)
$$

\subsection{Confluent Hamiltonians}

As we saw in subsection \ref{suse:1p1},  in the case of the $1+1$ confluence, the confluent Hamiltonians can be obtained as the limits of some functions on a phase space - linear combinations of the initial Hamiltonians with coefficients depending on a small parameter $\varepsilon$. Moreover the procedure of taking this limit requires the introduction of some shifts by Casimirs, since the Hamiltonians are defined up to Casimir element of the Poisson algebra. Such Casimir normalisation may be exploited in the case of the confluence for the higher order poles, however, this procedure becomes very heavy. In this section, 
we calculate these limits using residue calculus. Let us start by explaining these limits of the Hamiltonians in the $1+1$ confluence procedure; we want to calculate limits in \eqref{eq:Lie-conf}:
$$\d \tilde{H}_i := \underset{\varepsilon\rightarrow 0}{\lim}\, \d H_i,\quad \d \tilde{H}_{n-1} := \underset{\varepsilon\rightarrow 0}{\lim}\, \d  (H_{n-1}+H_{n}),\quad \d \tilde{H}_{n} := \underset{\varepsilon\rightarrow 0}{\lim}\, \varepsilon \d H_{n}.
$$
The Fuchsian Hamiltonians \eqref{sect1:SchlesingerHam} can be written in the following form
\begin{equation}\label{eq:ham-int}
H_{u_i}  = \frac{1}{2}\oint\limits_{\Gamma_{u_i}} \tr \left(A^2\right) \d \lambda,
\end{equation}
where $\Gamma_{u_i}$ is a contour that contains no singularities except $u_i$. Since the matrix $A$ admits a finite limit as $\varepsilon\to 0$,
the integrand has a finite limit. When $u_i\not= v_1, u$, we can always deform $\Gamma_{u_i}$ in such a way that the coalescence of $v_1$ and $u$ does not affect the contour of integration. This allows us to switch the limit and the integration operations, which gives the formula for the confluent Hamiltonian
\begin{equation}
\tilde{H}_{u_i}  = \frac{1}{2}\oint\limits_{\Gamma_{u_i}} \tr \left(\tilde{A}^2\right) \d \lambda, \qquad u_i\not= t_1, u.
\end{equation}

\begin{figure}
\begin{center}
\begin{tikzpicture}[scale=0.8]
\coordinate (A) at (-2,0);
\coordinate (B) at (0,2);
\draw[rounded corners,thick, color=Blue] (A) rectangle (B) {};
\draw[->,thick, color=Blue] (A)+(0,1.7)--(-2,1) ;
\draw[<-,thick, color=Blue] (B)+(0.0,-0.8);
\draw[->,thick, color=Blue] (A)+(0.3,0) -- (-1,0);
\draw[->, thick, color=Blue] (B)+(-0.3,0) -- (-1,2);
\node[color=red] at (-1,1) {\textbullet};
\node[color=black] at (-1,-0.5) {$\Gamma_u$};
\node[color=black] at (-1,-5.5) {$\Gamma_u$};
\node[color=black] at (-1,0.7) {$u$};


\coordinate (A) at (0.5,0);
\coordinate (B) at (2.5,2);

\draw[rounded corners,thick, color=darklavender] (A) rectangle (B) {};
\draw[->,thick, color=darklavender] (A)+(0,1.7)--(0.5,1) ;
\draw[<-,thick, color=darklavender] (B)+(0.0,-0.8);
\draw[->,thick, color=darklavender] (A)+(0.3,0) -- (1.5,0);
\draw[->, thick, color=darklavender] (B)+(-0.3,0) -- (1.5,2);
\node[color=red] at (1.5,1) {\textbullet};
\node[color=black] at (1.5,0.7) {$v_1$};
\node[color=black] at (1.5,-0.55) {$\Gamma_{v_1}$};
\node[color=black] at (1.5,-5.55) {$\Gamma_{v_1}$};
\node[color=black] at (10.25,-5.55) {$\Gamma_u\cup \Gamma_{v_1}$};

\draw[->, dashed,thick, color=Blue] (0.25,-0.5)--(0.25,-2.5);
\node[] at (-1.2,-1.7) {$\begin{array}{c}
    \operatorname{Contour}  \\
    \operatorname{deformation}
\end{array}$ };
\draw[->, dashed,thick, color=Blue] (4,-4)--(7,-4);
\node[] at (5.5,-3.3) {$\begin{array}{c}
    \operatorname{Contour}  \\
    \operatorname{merge}
\end{array}$ };

\draw[->, dashed,thick, color=Blue] (10.25,-2.5)--(10.25,-0.5);
\node[] at (8.7,-1.5) {$\begin{array}{c}
    \operatorname{Singularity}  \\
    \operatorname{coalescence}
\end{array}$};

\draw[Blue,thick,rounded corners] (-2, -5) rectangle (0.25, -3) {};
\draw[thick, darklavender,rounded corners] (0.25, -5) rectangle (2.5, -3) {};
\draw[<-,thick, color=darklavender] (B)+(0.0,-5.8);
\draw[<-,thick, color=darklavender] (1.3,-3)--(2.1,-3);
\draw[<-,thick, color=darklavender] (1.5,-5)--(0.4,-5);

\draw[<-,thick, color=Blue] (-2,-4.5)-- (-2,-4);
\draw[->,thick, color=Blue] (-1.7,-5)--(-1,-5);
\draw[<-,thick, color=Blue] (-1,-3)--(0.1,-3);
\node[color=red] at (0,-4) {\textbullet};
\node[color=black] at (0.7,-4.3) {$v_1$};
\node[color=red] at (0.5,-4) {\textbullet};
\node[color=black] at (-0.2,-4.3) {$u$};

\draw[rounded corners, color=amber,thick] (8, -5) rectangle (12.5, -3) {};
\draw[->,thick, color=amber] (8.3,-5)--(10.25,-5);
\draw[<-,thick, color=amber] (10.25,-3)--(12,-3);
\draw[<-,thick, color=amber] (8,-4)--(8,-3.5);
\draw[->,thick, color=amber] (12.5,-4.7)--(12.5,-4);
\node[color=red] at (10,-4) {\textbullet};
\node[color=black] at (10,-4.3) {$u$};
\node[color=black] at (10.56,-4.33) {$v_1$};
\node[color=red] at (10.5,-4) {\textbullet};

\draw[rounded corners,color=amber,thick] (8, 0) rectangle (12.5, 2) {};
\draw[->,thick, color=amber] (8.3,0)--(10.25,0);
\draw[<-,thick, color=amber] (10.25,2)--(12,2);
\draw[<-,thick, color=amber] (8,1)--(8,1.7);
\draw[->,thick, color=amber] (12.5,0.3)--(12.5,1);
\draw[color=Blue] (10.25,1) circle (0.25cm);
\node[color=red] at (10.5,1) {\textbullet};
\node[color=black] at (9.8,0.7) {$u$};
\node[color=green] at (10,1) {\textbullet};
\node[color=black] at (11,-0.5) {$\widetilde\Gamma_{u}$};

\end{tikzpicture}
\end{center}
\caption{Poles and contours confluence procedure}
\label{Fig:contours}
\end{figure}
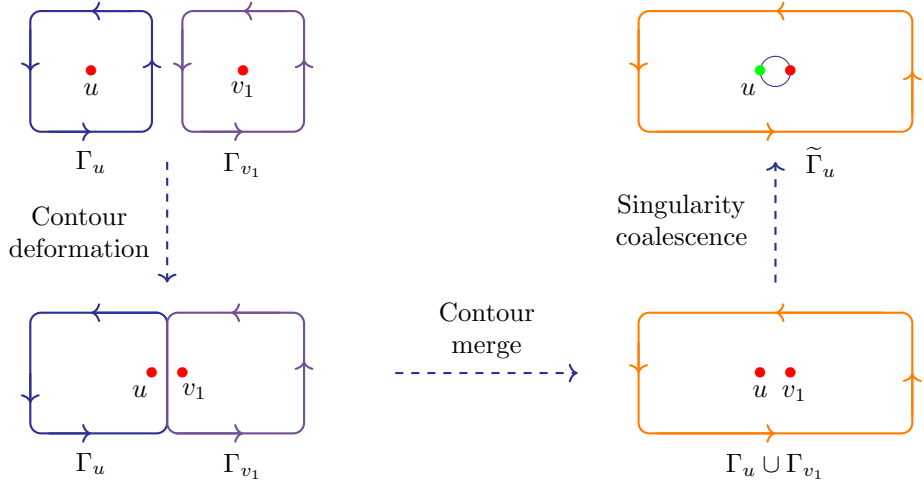

Let us now deal with the limit of $H_w+H_1$. Because both contours $\Gamma_w$ and $\Gamma_1$ will depend on $\varepsilon$, we cannot calculate the limits of  $H_w$ and $H_1$ separately. However, we can calculate the limit of the sum due to 
$$
H_u+ H_{v_1} = \frac{1}{2}\oint\limits_{\Gamma_{w}} \tr \left(A^2\right) \d \lambda+\frac{1}{2}\oint\limits_{\Gamma_{v_1}} \tr \left(A^2\right) \d \lambda = \frac{1}{2}\oint\limits_{\Gamma_{w}\cup \Gamma_{v_1}} \tr \left(A^2\right) \d \lambda,
$$
where the last equality holds since the integrands in both integrals are the same and $\Gamma_{u} \cup \Gamma_{v_1}$ denotes the contour obtained by merging $\Gamma_{u} $ and $\Gamma_{1}$ as illustrated in Fig. \ref{Fig:contours}. Such contour may be deformed to the contour $\tilde{\Gamma}_{u}$, such that the coalescent singularities are located inside this contour and the confluence doesn't affect the contour itself. Using the same arguments as before we obtain that
\begin{equation}
    \tilde{H}_u =\lim_{\varepsilon\to 0} \frac{1}{2}\oint\limits_{\Gamma_{w}\cup \Gamma_{v_1}} \tr \left(A^2\right) \d \lambda,
=\frac{1}{2}\oint\limits_{\Gamma_{w}} \tr \left(\tilde{A}^2\right) \d \lambda.
\end{equation}

In order to obtain $\tilde{H}_{1}$ we consider the following sum of Casimirs
$$
\frac{1}{2}\oint\limits_{\Gamma_{u}}(\lambda - u) \tr \left(A^2\right)\d \lambda +\frac{1}{2} \oint\limits_{\Gamma_{v_1}}(\lambda -  {v_1}) \tr \left(A^2\right)\d \lambda 
$$
which may be put to zero since the Hamiltonians are defined up to Casimirs. Expanding $v_1$  in $\varepsilon$ we obtain
$$
\frac{1}{2}\oint\limits_{\Gamma_{u}}(\lambda - u) \tr \left(A^2\right)\d \lambda +\frac{1}{2} \oint\limits_{\Gamma_{v_1}}(\lambda - u) \tr \left(A^2\right)\d \lambda - t_1 \varepsilon \frac{1}{2} \oint\limits_{\Gamma_{v_1}} \tr \left(A^2\right)\d \lambda  =  \frac{1}{2}\oint\limits_{\Gamma_{u}\cup \Gamma_{v_1}}(\lambda - u) \tr \left(A^2\right)\d \lambda - t_1 \varepsilon H_1 = 0.
$$
The relation written above finally gives us
\begin{equation}
\tilde{H}_1 = \frac{1}{2\tilde{t}_1}\oint\limits_{\Gamma_{u}}(\lambda - u) \tr \left(\tilde{A}^2\right)\d \lambda = \underset{\lambda=u}{\operatorname{Res}} \left[\frac{(\lambda - u)}{t_1} \tr\left(\frac{\tilde{A}^2}{2}\right)\right].
\end{equation}

Let us now add one more simple poles to $w$ using the confluence, by a similar computation as above, we obtain the following isomonodromic Hamiltonians
$$
H_{u_i} = \underset{\lambda=u_i}{\operatorname{Res}} \tr\left(\frac{A^2}{2}\right),\quad H_u =\underset{\lambda=u}{\operatorname{Res}} \tr\left(\frac{A^2}{2}\right)
$$
$$
H_1 = \underset{\lambda=u}{\operatorname{Res}}\left(\left[\frac{(\lambda-u)}{t_1}-\frac{t_2(\lambda-u)^2}{t_1^3} \right]\tr\left(\frac{A^2}{2}\right)\right), \quad H_2=\underset{\lambda=u}{\operatorname{Res}}\left[\frac{(\lambda-u)^2}{t_1^2}\tr\left(\frac{A^2}{2}\right)\right].
$$
The form of the Hamiltonians for $t_1$ and $t_2$ looks quite bizarre, but we may obtain them by solving the following linear system
$$
\mathcal{M}^{(2)}\left(\begin{array}{c}
    H_1  \\
    H_2
\end{array}\right)=\left(\begin{array}{cc}
    t_1 & t_2 \\
    0 & t_1^2
\end{array}\right)\left(\begin{array}{c}
    H_1  \\
    H_2
\end{array}\right) = \left(\begin{array}{c}
    S_1  \\
    S_2
\end{array}\right),\quad S_i = \frac{1}{2}\oint\limits_{\Gamma_u}(\lambda-u)^i\tr A^2 \d \lambda.
$$
Here $\mathcal{M}^{(2)}$ is a matrix which entries were already introduced in (\ref{sect3:monom}).\\

We now prove Theorem \ref{th:mainHam}:\\

\begin{theorem*}
Let $u$ be a pole of a connection $A$ with Poincar\'e rank $r$, which is the  result of confluence of $r$ simple poles with the simple pole $u$. Then the confluent Hamiltonians $H_1,\dots,H_{r}$ which correspond to the times $t_1,\dots t_r$ are defined as follows:
\begin{equation}\label{eq:ham-spec}
    \left(\begin{array}{c}
    H_1 \\
    H_2 \\
    \dots \\
    H_r
\end{array}\right) = \left(\mathcal{M}^{(r)}\right)^{-1}\left(\begin{array}{c}
    S^{(u)}_1 \\
    S^{(u)}_2 \\
    \dots \\
    S^{(u)}_r
\end{array}\right),
\end{equation}
where 
\begin{equation}
S_k^{(u)}= \frac{1}{2}\oint\limits_{\Gamma_u}(\lambda-u)^k\tr A^2 \d \lambda
\end{equation}
are spectral invariants of order $i$ in $u$ 
and the matrix $\mathcal{M}^{(r)}$  has entries $\mathcal{M}^{(r)}_{k,j}$  given by  (\ref{sect3:monom}). The Hamiltonian $H_u$ corresponding to the time $u$ is instead given by the standard formula
$$
{H}_{u_i} = \frac{1}{2}\underset{\lambda=u_i}{\res}\tr{A}(\lambda)^2.
$$
\end{theorem*}

\begin{proof}
To prove this theorem we again use induction. We showed that it holds for $r=2$ and it is trivial in case when $r=1$. Now we want to prove that
if the statement of the theorem holds for rank $r$ it is also true for rank $r+1$. The confluence expansion 
$$
v_{r+1}=u+P_{r+1}(t,\varepsilon)= u + \sum\limits_{i=1}^{r+1}\varepsilon^i t_i
$$
sends the extended symplectic form
$$
\omega = \d t_1\wedge \d H_1+ \dots +  \d t_r\wedge \d H_r + \d u\wedge \d H_u  + \d v_{r+1}\wedge \d H_{r+1} +\dots
$$
where the last set of dots corresponds to the terms that are not changed in the confluence procedure, to
$$
  \d t_1\wedge \d (H_1+\varepsilon H_{r+1}) + \dots+ \d t_r\wedge \d (H_r+\varepsilon^r H_{r+1}) + \d t_{r+1}\wedge \d (\varepsilon^{r+1} H_{r+1}) +\d u\wedge \d (H_u+H_{r+1}) +\dots=
  $$
  that must be equal to 
  $$
 \d t_1\wedge \d \tilde{H}_1 \dots  \d t_r\wedge \d \tilde{H}_r + \d t_{r+1}\wedge \d \tilde{H}_{r+1}+ \d u\wedge \d \tilde{H}_u + \dots.
$$
Therefore, we must find the limits
$$
\tilde{H}_{u} = \underset{\varepsilon\rightarrow0}{\lim} \left(H_u + H_{r+1}\right),\quad \tilde{H}_{k} = \underset{\varepsilon\rightarrow0}{\lim}  \left(H_k + \varepsilon^k H_{r+1}\right),\quad k=1\dots r,\quad \tilde{H}_{r+1} = \underset{\varepsilon\rightarrow0}{\lim}  \varepsilon^{r+1} H_.
$$
The first limit is quite simple and may be obtained via the union of contours which we already described before. To find the other limits, let us consider the relation
$$
\oint\limits_{\Gamma_{u}} (\lambda-u)^i\tr \frac{A^2}{2} \d \lambda  + \oint\limits_{\Gamma_{r+1}} (\lambda-v_{r+1})^i\tr \frac{A^2}{2} \d \lambda = S^{(u)}_i \mod \left(\operatorname{Casimirs}\right)
$$
expanding $v_{n+1}$, we obtain
$$
\oint\limits_{\Gamma_{u}\cup \Gamma_{r+1}} (\lambda-u)^i\tr \frac{A^2}{2} \d \lambda - \oint\limits_{\Gamma_{r+1}}\phi(\lambda)\tr \frac{A^2}{2} = S^{(u)}_i \mod \left(\operatorname{Casimirs}\right),
$$
where $\phi(\lambda)$ is a holomorphic function inside $\Gamma_{r+1}$ which is given by
$$
\phi(\lambda) =(\lambda-u)^i- (\lambda - u - P_{r+1}(t,\varepsilon))^{i}.
$$
Since $\phi(\lambda)$ has no zeros at $v_{r+1}$
we have
$$
\oint\limits_{\Gamma_{r+1}}\phi(\lambda)\tr \frac{A^2}{2} = \phi(u+P_{r+1}(t,\varepsilon))\oint\limits_{\Gamma_{r+1}}\tr \frac{A^2}{2} = P_{r+1}(t,\varepsilon)^{i}\oint\limits_{\Gamma_{r+1}}\tr \frac{A^2}{2} = P_{r+1}(t,\varepsilon)^{i}H_{r+1}.
$$
Finally, we obtain the following identity:
\begin{equation}\label{eq:S_rel}
\oint\limits_{\Gamma_{u}\cup \Gamma_{r+1}} (\lambda-u)^i\tr \frac{A^2}{2} \d \lambda - P_{r+1}(t,\varepsilon)^{i}H_{r+1}= S^{(u)}_i \mod \left(\operatorname{Casimirs}\right).
\end{equation}
In the case $i=r+1$, $S^{(u)}_{r+1}$ is a Casimir {due to the formula \eqref{eq:quadCas}}, therefore we have
$$
\oint\limits_{\Gamma_{u}\cup \Gamma_{r+1}} (\lambda-u)^{r+1}\tr \frac{A^2}{2} \d \lambda = P_{r+1}(t,\varepsilon)^{r+1}H_{r+1} \mod \left(\operatorname{Casimirs}\right).
$$
The left hand side of this identity has a finite limit when $\varepsilon$ goes to $0$. Indeed, since the contour contains both $u$ and $v_{r+1}$ the confluence procedure doesn't affect it and the only dependence in $\varepsilon$  is in $A$. According to Theorem \ref{sect4:small_theorem}, $A$ has a finite limit, the same has $\tr A^2$, so we have that
$$
\underset{\varepsilon\rightarrow0}{\lim}\oint\limits_{\Gamma_{u}\cup \Gamma_{r+1}} (\lambda-u)^{r+1}\tr \frac{A^2}{2} \d \lambda =\oint\limits_{\Gamma_{u}\cup \Gamma_{r+1}} (\lambda-u)^{r+1}\underset{\varepsilon\rightarrow0}{\lim}\tr \frac{A^2}{2} \d \lambda = \oint\limits_{\Gamma_{u}} (\lambda-u)^{r+1}\tr \frac{\tilde{A}^2}{2} \d \lambda = \tilde{S}^{(u)}_{r+1},
$$
where $\tilde{S}^{(u)}_{r+1}$ is a spectral invariant of the confluent system with connection $\tilde{A}$. Since after the confluence, the order of pole increases to $r+2$, such spectral invariant is not a Casimir for the confluent system. This means, that the limit of $P_{r+1}(t,\varepsilon)^{r+1}H_{r+1}$ up to Casimirs exists and equals to $ \tilde{S}^{(u)}_{r+1}$. On the other hand we have
$$
P_{r+1}(t,\varepsilon)^{r+1}H_{r+1} = \left(t_1^{r+1}\varepsilon^{r+1}+ O(\varepsilon^{r+2})\right) H_{r+1},
$$
and since the limit exists up to Casimirs we get that
$$
H_{r+1} \mod (\operatorname{Casimirs}) = \frac{\tilde{S}^{(u)}_{r+1}}{\varepsilon^{r+1}}+\sum\limits_{i=-r}^{\infty}H_{r+1}^{[i]}\varepsilon^i,
$$
so in principle $H_{r+1}$ may have terms of lower order than $1/\varepsilon^{r+1}$, but these terms have to be Casimirs.
Considering  the relations \eqref{eq:S_rel} for $i=1\dots r+1$ as a linear system, we obtain

\begin{equation}\label{sect3:linSysPr}
\left(\begin{array}{c}
    \tilde{S}_1 \\
    \tilde{S}_2 \\
    \tilde{S}_3 \\
    \dots \\
    \tilde{S}_{r+1}
\end{array}\right) - \mathcal{M}^{(r+1)}\left(\begin{array}{c}
    \varepsilon \\
    \varepsilon^2 \\
    \varepsilon^3 \\
    \dots \\
    \varepsilon^{r+1}
\end{array}\right) H_{r+1}= \left(\begin{array}{c}
    S_1 \\
    S_2 \\
    \dots \\
    S_r \\
    0
\end{array}\right) \mod \left(\operatorname{Casimirs}\right)
\end{equation}
where
$$
 \tilde{S}_i = \oint\limits_{\Gamma_{u}\cup \Gamma_{r+1}} (\lambda-u)^i\tr \frac{A^2}{2} \d \lambda,\quad S_i = \oint\limits_{\Gamma_{u}} (\lambda-u)^i\tr \frac{A^2}{2} \d \lambda.
$$
Note that the contours in the definition of $\tilde{S}_i$  are not affected by the confluence procedure. Using the same arguments as above, we compute the limits of these integrals, which are 
$$
\underset{\varepsilon\rightarrow0}{\lim}\tilde{S}_i = \tilde{S}^{(u)}_{i} = \oint\limits_{\Gamma_{u}} (\lambda-u)^i\tr \frac{\tilde{A}^2}{2} \d \lambda,
$$
where $\tilde{S}^{(u)}_{i}$ denote the spectral invariants of the confluent system with connection $\tilde{A}$.

The crucial point here is that the matrix $\mathcal{M}^{(r+1)}$ contains $\mathcal{M}^{(r)}$ as $r+1,r+1$ minor, i.e.
$$
\mathcal{M}^{(r+1)} = \left(\begin{array}{cc}
    \mathcal{M}^{(r)} & \begin{array}{c}
         t_{r+1}  \\
         \vdots  \\
         \vdots \\
        \,
    \end{array} \\
    \begin{array}{cccc}
    0 &  0 & \dots & 0
    \end{array} &  t^{r+1}_1
\end{array}\right).
$$
Now let us consider the following matrix
$$
\mathcal{C} = \left(\begin{array}{cc}
    \left(\mathcal{M}^{(r)}\right)^{-1} & \begin{array}{c}
         0  \\
         \vdots  \\
         \vdots \\
        0
    \end{array} \\
    \begin{array}{cccc}
    0 &  0 & \dots & 0
    \end{array} &  1
\end{array}\right)
$$
and let us act via $\mathcal{C}$ on the equation (\ref{sect3:linSysPr}) from the left. Then we obtain
\begin{multline*}
\mathcal{C}\tilde{S}  = \left(\begin{array}{cc}
    \ID_r & \begin{array}{c}
         \vdots \\
         \vdots \\
         \, 
    \end{array} \\
    \begin{array}{cccc}
    0 &  0 & \dots & 0
    \end{array} &  t_1^{r+1}
\end{array}\right) \left(\begin{array}{c}
    \varepsilon H_{r+1} \\
    \varepsilon^{2} H_{r+1} \\
    \dots \\
    \varepsilon^{n} H_{r+1} \\
    \varepsilon^{n+1} H_{r+1}
\end{array}\right)+ \left(\begin{array}{c}
    H_1 \\
    H_2 \\
    \dots \\
    H_r \\
    0
\end{array}\right) = 
\\
= \left(\begin{array}{cc}
    \ID_r & \begin{array}{c}
         \vdots \\
         \vdots \\
         \, 
    \end{array} \\
    \begin{array}{cccc}
    0 &  0 & \dots & 0
    \end{array} &  t_1^{r+1}
\end{array}\right) \left(\begin{array}{c}
    \varepsilon H_{r+1} +H_1\\
    \varepsilon^{2} H_{r+1} + H_2 \\
    \dots \\
    \varepsilon^{r} H_{r+1} + H_{r} \\
    \varepsilon^{r+1} H_{r+1}
\end{array}\right) = \mathcal{C} \mathcal{M}^{(r+1)}\left(\begin{array}{c}
    \varepsilon H_{r+1} +H_1\\
    \varepsilon^{2} H_{r+1} + H_2 \\
    \dots \\
    \varepsilon^{r} H_{r+1} + H_{r} \\
    \varepsilon^{n+1} H_{r+1}
\end{array}\right).
\end{multline*}
In this way we have arranged the entries of equation (\ref{sect3:linSysPr}) in such a way that  the left hand side  has a nice limit as $\varepsilon$ goes to zero (the confluence of points is inside the contour of integration for $\tilde{S}_i$). On the right hand side we have the functions whose limits we want to find. Finally, multiplying by $\mathcal{C}^{-1}$ from the left we obtain
$$
\mathcal{M}^{(r+1)}\left(\begin{array}{c}
    \tilde{H}_{1}\\
    \tilde{H}_{2} \\
    \dots \\
    \tilde{H}_{r} \\
    \tilde{H}_{r+1}
\end{array}\right) = \mathcal{M}^{(r+1)} \underset{\varepsilon\rightarrow0}{\lim}\left(\begin{array}{c}
    \varepsilon H_{r+1} +H_1\\
    \varepsilon^{2} H_{r+1} + H_2 \\
    \dots \\
    \varepsilon^{r} H_{r+1} + H_{n} \\
    \varepsilon^{r+1} H_{r+1}
\end{array}\right) = \underset{\varepsilon\rightarrow0}{\lim}\left(\begin{array}{c}
    \tilde{S}_1 \\
    \tilde{S}_2 \\
    \tilde{S}_3 \\
    \dots \\
    \tilde{S}_{r+1}
\end{array}\right) =\left(\begin{array}{c}
    \tilde{S}_1 \\
    \tilde{S}_2 \\
    \tilde{S}_3 \\
    \dots \\
    \tilde{S}_{r+1}
\end{array}\right) 
$$
which concludes the proof. \end{proof}

\begin{remark}
The matrix $\mathcal{M}^{(r)}$ is automatically upper triangular matrix, so it is quite easy to solve such a system for any reasonable $n$.
\end{remark}

We may consider Hamiltonians which are related to poles locations as  spectral invariants $S^{(u)}_0$. Then it is easy to extend  formula from \eqref{eq:ham-spec} as follows
$$
{\mathcal{N}}^{(r)}\left(\begin{array}{c}
     H_0  \\
     H_1 \\
     \dots\\
     H_r
\end{array}\right) = \left(\begin{array}{c}
     S^{(u)}_0  \\
     S^{(u)}_1 \\
     \dots\\
     S^{(u)}_r
\end{array}\right),\quad \mathcal{N}^{(r)} = \left(\begin{array}{cc}
    1 & \begin{array}{ccc}
        0 & \dots & 0
    \end{array}  \\
\begin{array}{c}
    0  \\
    \vdots  \\
    0
\end{array} & \mathcal{M}^{(r)}
\end{array}\right).
$$

\subsection{Examples of Hamiltonians}
In order to see all the features of the obtained Hamiltonians, we consider a connection with $3$ simple poles at $0,1$ and $\infty$, and one irregular singularity at some point $u$. The simplest example is
\begin{equation}\label{ex:ham-32}
    A(\lambda) = \frac{A^{(0)}}{\lambda}+\frac{A^{(1)}}{\lambda-1}+\frac{A^{(u)}_0}{\lambda-u}+\frac{t_1A^{(u)}_1}{(\lambda-u)^2}.
\end{equation}
The explicit formulas for the Hamiltonians are
\begin{eqnarray}
      && H_u= \tr \left[A^{(u)}_0\left(\frac{A^{(0)}}{u}  + \frac{A^{(1)} }{u - 1}\right) - t_1 A^{(u)}_1 \left( \frac{ A^{(1)}}{(u - 1)^2} + \frac{A^{(0)}}{u^2}\right)\right]\\
      && H_1= \tr\left[  A^{(u)}_1\left(\frac{A^{(0)} }{u} + \frac{ A^{(1)}}{u - 1}\right)+\frac{A^{(u)}_0A^{(u)}_0}{2 t_1 }\right]
\end{eqnarray}
In the lifted Darboux coordinates, the Hamiltonians take the form
\begin{multline}
     H_u= \tr\left(Q_0^{(u)}P_0^{(u)}+Q_1^{(u)}P_1^{(u)}\right) \left(\frac{Q^{(0)}P^{(0)}}{u}+\frac{Q^{(1)}P^{(1)}}{u-1}\right)-\\  -t_1Q_0^{(u)}P_1^{(u)} \left(  \frac{ Q^{(1)}P^{(1)}}{(u - 1)^2} - \frac{Q^{(0)}P^{(0)}}{u^2}\right)\qquad\qquad
\end{multline}
\begin{multline}
       H_1=\tr  \frac{Q^{(0)}P^{(0)}  Q^{(u)}_0P^{(u)}_1}{u} + \frac{Q^{(1)}P^{(1)}  Q^{(u)}_0P^{(u)}_1}{u - 1}+\frac{\left(Q_0^{(u)}P_0^{(u)}+Q_1^{(u)}P_1^{(u)}\right)^2}{2 t_1 }.\qquad\qquad\qquad\quad
\end{multline}
Using the lifted Darboux coordinates, it is a straightforward computation to check that the Hamiltonians Poisson commute, moreover we may check that the cross-derivative w.r.t. $u$ and $t_1$ is zero, i.e.
$$
\frac{\partial}{\partial t_1} H_u - \frac{\partial}{\partial u} H_1 =0,
$$
which tells us that the $\tau$-function
\begin{multline*}
\d \ln \tau = \tr\left( A^{(u)}_0\left(\frac{A^{(0)}}{u}  + \frac{A^{(1)} }{u - 1}\right) - t_1 A^{(u)}_1 \left( \frac{ A^{(1)}}{(u - 1)^2} - \frac{A^{(0)}}{u^2}\right)\right) \d u + \\ +   \tr \left(  A^{(u)}_1\left(\frac{A^{(0)} }{u} + \frac{ A^{(1)}}{u - 1}\right)+\frac{A^{(u)}_0A^{(u)}_0}{2 t_1 }\right)\d t_1
\end{multline*}
is defined correctly. If we go further and consider a pole of Poincar\'e rank 2 at $u$, the matrix $A(\lambda)$ takes form
\begin{equation}\label{sec-ex}
    A(\lambda) = \frac{A^{(0)}}{\lambda}+\frac{A^{(1)}}{\lambda-1}+\frac{A^{(u)}_0}{\lambda-u}+\frac{t_1A^{(u)}_1+t_2A^{(u)}_2}{(\lambda-u)^2}+\frac{t_1^2A^{(u)}_2}{(\lambda-u)^3}.
\end{equation}
Then, the Hamiltonians write as
\begin{multline*}
H_u = \tr\Bigg[A^{(u)}_0\left(\frac{A^{(0)}}{u} + \frac{A^{(1)}}{u - 1}\right)   -t_1A^{(u)}_1 \left(\frac{A^{(0)}}{u^2} + \frac{A^{(1)}}{(u - 1)^2}\right)   +\\+ A^{(u)}_2\left( \frac{t_1^2A^{(0)}}{u^3} + \frac{t_1 ^2 A^{(1)}}{(u - 1)^3}-\frac{t_2 A^{(0)}}{u^2} - \frac{t_2 A^{(1)}}{(u - 1)^2} +\right)\Bigg] 
\end{multline*}
\begin{multline*}
H_1=\tr\Bigg[A^{(u)}_1\left(\frac{A^{(0)}}{u} + \frac{A^{(1)}}{u - 1}\right)   -t_1 A^{(u)}_2 \left(\frac{A^{(0)}}{u^2} + \frac{A^{(1)}} {(u - 1)^2}\right)  + \frac{A^{(u)}_0 A^{(u)}_0}{2t_1} - t_2\frac{A^{(u)}_0  A^{(u)}_1 }{t_1^2} - t_2^2\frac{A^{(u)}_0  A^{(u)}_2}{t_1^3}\Bigg]
\end{multline*}
\begin{multline*}
H_2=\tr\Bigg[A^{(u)}_2\left(\frac{A^{(0)}}{u} + \frac{A^{(1)}}{u - 1}\right)  +\frac{A^{(u)}_0  A^{(u)}_1}{t_1}  + t_2  \frac{A^{(u)}_0  A^{(u)}_2}{t_1^2}\Bigg].\qquad\qquad\qquad\qquad\qquad\qquad\qquad\qquad
\end{multline*}    
As in the previous example, the cross-derivatives are zero
$$
\frac{\partial}{\partial t_1} H_u - \frac{\partial}{\partial u} H_1 = \frac{\partial}{\partial t_2} H_u - \frac{\partial}{\partial u} H_2 = \frac{\partial}{\partial t_1} H_2 - \frac{\partial}{\partial t_2} H_1 =0,
$$
so the $\tau$-function is defined correctly.

In the case of $\mathfrak{sl}_2$, the isomonodromic deformations of \eqref{ex:ham-32}  correspond to the confluence of the two-time Garnier system and belong to the list of the so-called fourth-order Painlev\'e equations introduced in \cite{Kaw1,Kaw2}.  The second example \eqref{sec-ex} is more complicated -- the Hamiltonian reduction  gives a Hamiltonian system with $3$ degrees of freedom, which corresponds to a sixth order Painlev\'e equation.

\subsection{Confluence of the higher order poles.} 
The confluence of two poles $w_1$ and $w_2$ of Poincar\'e rank $r_1$ and $r_2$ respectively can be treated a the confluence of $r_1+r_2+2$ simple poles. Indeed, we have seen at the beginning of section 5.2 that the generic connection with a Poincar\'e rank $r$ singularity can be obtained as confluence of $r+1$ simple poles. Therefore, a connection with  two poles $w_1$ and $w_2$ of Poincar\'e rank $r_1$ and $r_2$ respectively is obtained by coalescing $r_1+1$ and $r_2+1$ simple poles.

\section{The non-ramified Painlev\'e equations}
\subsection{General scheme}
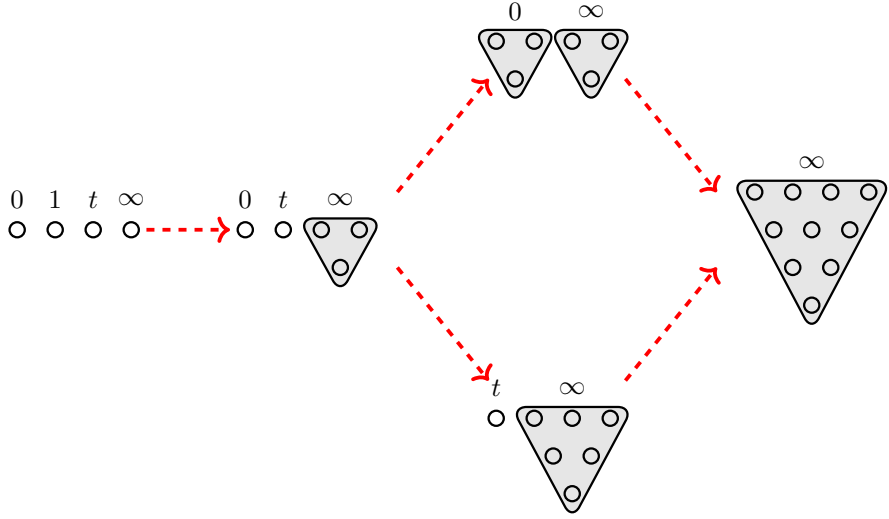
\begin{figure}
    \centering
\begin{tikzpicture}
\draw[line width=0.3mm] 
      (-13.5,0) circle (1mm)
      (-13,0) circle (1mm)
      (-12.5,0) circle (1mm)
      (-12,0) circle (1mm);
     
    \node[] at (-13.5,0.4) {$0$};
    \node[] at (-13,0.4) {$1$};
    \node[] at (-12.5,0.4) {$t$};
    \node[] at (-12,0.4) {$\infty$};
      
    \draw[line width=0.5mm, dashed, red, ->](-11.8,0)--(-10.7,0);
    \draw[line width=0.3mm,rounded corners=2mm,fill=gray!20, opacity=5] (-9.8,0.15)--(-9.25,-0.85)--(-8.7,0.15)--cycle;

    \draw[line width=0.3mm] (-10.5,0) circle (1mm)
      (-10,0) circle (1mm)
      (-9.5,0) circle (1mm)
      (-9,0) circle (1mm)
      (-9.25,-0.5) circle (1mm);

    \node[] at (-10.5 ,0.4) {$0$};
    \node[] at (-10,0.4) {$t$};
    \node[] at (-9.25,0.4) {$\infty$};
    
    \draw[line width=0.5mm, dashed, red,  ->](-8.5,0.5)--(-7.3,2.0);
    \draw[line width=0.5mm, dashed, red,  ->](-8.5,-0.5)--(-7.3,-2.0);
    \draw[line width=0.3mm,rounded corners=2mm,fill=gray!20, opacity=5] (-7.5,2.65)--(-6.95,1.65)--(-6.4,2.65)--cycle
      [rounded corners=2mm,fill=gray!20, opacity=5] (-6.5,2.65)--(-5.95,1.65)--(-5.4,2.65)--cycle;
    \draw[line width=0.3mm] (-7.2,2.5) circle (1mm)
      (-6.7,2.5) circle (1mm)
      (-6.2,2.5) circle (1mm)
      (-5.7,2.5) circle (1mm)
      (-5.95,2) circle (1mm)
      (-6.95,2) circle (1mm);
      
    \node[] at (-6.95 ,2.9) {$0$};
    \node[] at (-5.95, 2.9) {$\infty$};
      
    \draw[line width=0.5mm, dashed, red, ->](-5.5,2.0)--(-4.3,0.5);
    \draw[line width=0.5mm, dashed, red, ->](-5.5,-2.0)--(-4.3,-0.5);
    \draw[line width=0.3mm,rounded corners=2mm,fill=gray!20, opacity=5] (-4.1,0.65)--(-3.05,-1.35)--(-2.0,0.65)--cycle;
    \draw[line width=0.3mm] (-3.8,0.5) circle (1mm)
      (-3.3,0.5) circle (1mm)
      (-2.8,0.5) circle (1mm)
      (-2.3,0.5) circle (1mm)
      (-3.55,-0) circle (1mm)
      (-3.05,-0) circle (1mm)
      (-2.55,-0) circle (1mm)
      (-3.3,-0.5) circle (1mm)
      (-2.8,-0.5) circle (1mm)
      (-3.05,-1) circle (1mm);
    \node[] at (-3.05, 0.9) {$\infty$};

    \draw[line width=0.3mm,rounded corners=2mm,fill=gray!20, opacity=5] (-7,-2.35)--(-6.2,-3.85)--(-5.4,-2.35)--cycle;
    \draw[line width=0.3mm] (-7.2,-2.5) circle (1mm)
      (-6.7,-2.5) circle (1mm)
      (-6.2,-2.5) circle (1mm)
      (-5.7,-2.5) circle (1mm)
      (-5.95,-3) circle (1mm)
      (-6.45,-3) circle (1mm)
      (-6.2,-3.5) circle (1mm);
     \node[] at (-7.2 ,-2.1) {$t$};
    \node[] at (-6.2,-2.1) {$\infty$};

     draw[line width=0.5mm, dashed, green, ->](-5.3,-2.5)--(-4.3,-2.5);
    
\end{tikzpicture}
\caption{Confluence scheme for Painlev\'e equations. Each triangle at this diagram corresponds to the Takiff algebra Darboux coordinates which were introduced in Fig. \ref{fig1:PQ}.}
    \label{fig3:Congluence}
\end{figure}

The confluence diagram of the Darboux parametrisations in the case of rank $2$ non-ramified connections with $4$ points is given at Fig. \ref{fig3:Congluence}. We start by illustrating the general scheme of  reduction which works for any rank.

The Hamiltonians of the isomonodromic problem with irregular singularities of Poincar\'e rank $r_i$ at a point $u_i$  admit additional symmetries with respect to the inner action (choice of the spectral invariants at each singularity) and outer action (gauge group action). Using the Darboux parametrisation of the co-adjoint orbit for the $\mathfrak{sl}_2$--Takiff algebras, we automatically fix the spectral invariants, i.e. reduce with respect to the inner action. The only symmetry which still needs to be taken into account is the gauge freedom which leads to the fully reduced phase space. In all the examples of this section, we write down the Darboux coordinates by immediately diagonalizing the leading terms at one of the irregular singularities.  Therefore, the number of the intermediate coordinates in all examples is $4$ and not $6$ (because we have eliminated $2$ coordinates by diagonalisation). Such coordinates are in correspondence with the Darboux coordinates which were used in \cite{HR}. In order to reduce to the smallest dimension of the system (namely $2$), we have to reduce with respect to the stabilizer action, which in all examples is equivalent to the additional action of the Cartan torus, since we consider only the unramified situation, which corresponds to the case when  Katz index is an integer. The ramified situation will be considered in the next work of the first author.

\subsection{Painlev\'e V} The isomonodromic problem takes the form
\begin{equation}
    \renewcommand{\arraystretch}{1.7}
    \begin{array}{l}
        \frac{d}{d\lambda }\Psi = \left(\frac{A^{(0)}}{\lambda}+\frac{A^{(t)}}{\lambda-t}+B_1\right)\Psi   \\
        \frac{d}{d t}\Psi = -\frac{A^{(t)}}{\lambda-t}\Psi.
    \end{array}
\end{equation}
the deformation equations are
\begin{equation}
    \frac{d}{dt}A^{(0)} = \frac{1}{t}[A^{(t)}, A^{(0)}],\quad \frac{d}{dt}A^{(t)} = \frac{1}{t}[A^{(0)}, A^{(t)}] + [B_1, A^{(t)}],\quad \frac{d}{dt} B_1 = 0.
\end{equation}
The Poisson algebra is given by
\begin{equation}
      \left\{A^{(i)}\ocomma A^{(i)}\right\}=[\Pi, \ID \otimes A^{(i)}],\quad \left\{A^{(0)}\ocomma A^{(t)}\right\}=\left\{A^{(0)}\ocomma B_1\right\}=\left\{A^{(t)}\ocomma B_1\right\}=\left\{B_1\ocomma B_1\right\}=0
\end{equation}
the isomonodromic Hamiltonian writes as
\begin{equation}
    H_{V} = \underset{\lambda=t}{\operatorname{res}}\tr \left(\frac{A(\lambda)^2}{2}\right) = \tr\left(A^{(t)}B_1 + \frac{1}{t}A^{(t)}A^{(0)}\right).
\end{equation}
In the $\mathfrak{sl}_2$ case, the Darboux parametrisation of the elements of the coadjoint orbit takes the  form
$$
A^{(0)} = \left(\begin{array}{cc}{p_0 q_0-\theta_0} & {-(p_0 q_0-2 \theta_0) p_0} \\ {q_0} & {-p_0 q_0+\theta_0}\end{array}\right),\quad A^{(t)} = \left(\begin{array}{cc}{p_t q_t-\theta_t} & {-(p_t q_t-2 \theta_t) p_t} \\ {q_t} & {-p_t q_t+\theta_t}\end{array}\right),
$$
with the symplectic form
$$
\omega = \d p_t\wedge \d q_t + \d p_0\wedge \d q_0.
$$

Using the gauge freedom, we set the constant matrix $B_1$ to be diagonal:
$$
B_1 = \left(\begin{array}{cc}
    k & 0 \\
    0 & -k
\end{array}\right).
$$
In such a parametrisation, the Hamiltonian takes form
\begin{multline*}
H= \underset{\lambda=t}{\operatorname{res}}\tr\left(\frac{A(\lambda)^2}{2}\right) = 4k(p_tq_t-\theta_{t})-\frac{2}{t}(q_{t}q_{0}(p_t-p_0)^2-2(q_{0}\theta_{t}-q_{t}\theta_{0})(p_t-p_0)-2\theta_{0}\theta_t )= \\
= 4k(p_tq_t-\theta_{t}) - \frac{2q_{t}q_{0}}{t}\left(p_t-p_0-\frac{\theta_{t}}{q_{t}}+\frac{\theta_{0}}{q_{0}}\right)^2+\frac{2}{t}\left(\theta_{t}^2\frac{q_{0}}{q_{t}}+\theta_{0}^2\frac{q_{t}}{q_{0}}\right).
\end{multline*}
This Hamiltonian is invariant under the following rescaling 
$$
p_i\rightarrow p_i\alpha,\quad q_i\rightarrow \frac{q_i}{\alpha}
$$
which is the same as the gauge $SL(2)$ action via diagonal matrix. The moment map is
$$
q_{0}p_{0}+q_{t}p_{t}.
$$
The change of coordinates
$$
I = q_{0}p_{0}+q_{t}p_{t},\quad \phi = \ln(q_{0}),\quad u=-\frac{q_{t}}{q_{0}},\quad v=p_tq_0,
$$
is a canonical transformation. Resolving it with respect to the $q$'s and $p$'s we obtain
$$
q_{0}=e^{\varphi}, \quad q_t = -e^{\varphi}u,\quad p_{0} = e^{-\varphi}(I+uv)\quad p_t=e^{-\varphi}v,
$$
and the symplectic form goes to
$$
\omega = \d p_t\wedge \d q_t + \d p_0\wedge \d q_0 = \d I\wedge \d \varphi + \d v\wedge \d u.
$$
The Hamiltonian in these coordinates writes as
$$
H=-4k(uv+\theta_t)+2\frac{u}{t}\left(v-I-uv+\frac{\theta_{t}}{u}+\theta_{0}\right)^2-\frac{2}{t}\left(\theta_{t}^2\frac{1}{u}+\theta_{0}^2u\right)
$$
and it is obvious that $I$ and $\varphi$ are the part of the action-angle variables, so we may decrease the number of degrees of freedom by $1$ and consider the following Hamiltonian system
$$
H=-4k(uv+\theta_t)+2\frac{u}{t}\left(v-a-uv+\frac{\theta_{t}}{u}+\theta_{0}\right)^2-\frac{2}{t}\left(\theta_{t}^2\frac{1}{u}+\theta_{0}^2u\right),\quad \omega = \d v \wedge \d u,\quad a=\operatorname{const}.
$$
The equations of motion take form
$$
\dot{u} = \frac{\partial H}{\partial v} = -4ku+\frac{4u}{t}(1-u)\left(v-a-uv+\frac{\theta_t}{u}+\theta_0\right)
$$
$$
\dot{v} = -\frac{\partial H}{\partial u} = 4kv-\frac{2}{t}\left((v-uv-a+\frac{\theta_{t}}{u}+\theta_{0})^2-2u(v-uv-a+\frac{\theta_{t}}{u}+\theta_{0})\left(v+\frac{\theta_{t}}{u^2}\right)+\frac{\theta_{t}^2}{u^2}-\theta_{0}^2\right).
$$
Writing second order ODE for $u$ we obtain
$$
\frac{d^2 u}{dt^2}= \left(\frac{1}{u-1}+\frac{1}{2u}\right)\left(\frac{du}{dt}\right)^2-\frac{1}{t}\frac{du}{dt}+8\theta_0\frac{(u-1)^2}{t^2}\left(u-\left(\frac{\theta_t}{\theta_0}\right)^2\frac{1}{u}\right)+4k(4(a-\theta_{0}-\theta_{t})-1)\frac{u}{t}-8k^2\frac{u(u+1)}{u-1}
$$
which is the Gambier's form of the Painlev\'e V equations and the constants are given by
$$
\theta_0=\frac{\alpha}{8},\quad \theta_t^2=-\frac{\alpha\beta}{64},\quad k^2=-\frac{\delta}{8},\quad 4k(4(a-\theta_{0}-\theta_{t})-1)=\gamma.
$$
The following canonical transformation
$$
u=\frac{x}{x-1},\quad v=-((x-1)y+a-2\theta_0)(x-1),\quad \d v \wedge \d u = \d y \wedge \d x,
$$
sends $H$ to the following form
$$
tH = 2x(x-1)y^2+4(ktx(x-1)+x(\theta_t-\theta_0)-\theta_t)y+4(xkt(a-2\theta_0)-\theta_t(kt-\theta_0))
$$
which was introduced in \cite{KMNOY}. The example of the Painlev\'e V equation as a system written on the co-adjoint orbits of the Takiff algebra was recently studied by \cite{KalininBabich} in more details.

\subsection{Painlev\'e IV}
The connection is
\begin{equation}
    A(\lambda) = \frac{A^{(t)}}{\lambda-t}-B_1-B_2\lambda
\end{equation}
and the deformation one-form is
\begin{equation}
    \Omega = - \frac{A^{(t)}}{\lambda-t}\d t.
\end{equation}
The deformation equations are
$$
\dot{A}^{(t)} = [A^{(t)}, B_1+B_2 t],\quad \dot{B}_1 = [B_2, A^{(t)}],\quad \dot{B}_2 = 0.
$$
The Poisson structure is
$$
\{A^{(t)} \ocomma A^{(t)}\} = [\Pi, 1\otimes C],\quad \{B_1 \ocomma B_1\} = [\Pi, 1\otimes B_2],\quad \{B_1 \ocomma B_2\}=\{B_2 \ocomma B_2\}=0.
$$
the Hamiltonian writes as
\begin{equation}\label{sect5:PIVH}
H= \underset{\lambda=t}{\res} \tr\frac{A^2}{2} = -\tr\left(A^{(t)}B_1+tA^{(t)}B_2\right).
\end{equation}
Since $B_3$ is a constant of motion, the same holds for the transition matrix to the eigenbasis of $B_3$. This allows us to consider the gauge, which is equal to this transition matrix without changing the Poisson structure of $A^{(t)}$. In the case of $\mathfrak{sl_2}$ we have
\begin{equation}
    A^{(t)} = \left( \begin {array}{cc} p_{{t}}q_{{t}}-\theta_{{t}}&- \left( p_{{t}
}q_{{t}}-2\theta_{{t}} \right) p_{{t}}\\ q_{{t}}&-
p_{{t}}q_{{t}}+\theta_{{t}}\end {array} \right),\quad -B_2\lambda-B_1=-\lambda\theta_3\left(\begin{array}{cc}
    1 & 0 \\
    0 & -1
\end{array}\right)-\left(\begin {array}{cc} \theta_{{2}}&-2\theta_{{3}}q_{{3}}
\\ \noalign{\medskip}p_{{3}}&-\theta_{{2}}\end {array} 
\right).
\end{equation}
The Hamiltonian writes as
\begin{equation}
    H =  \left( p_{{t}}q_{{t}}-2\theta_{{t}} \right) p
_{{t}}p_{{3}}-2 \left( p_{{t}}q_{{t}}-\theta_{{t}} \right)  \left( t\theta_{{3}}+
\theta_{{2}} \right) +2\theta_{{3}}q_{{3}}q_{{t}}.
\end{equation}
Since $B_3$ is a diagonal matrix (has no Jordan blocks) the stabilizer is the Cartan torus of $SL_2$,i.e.
$$
S = \left(\begin{array}{cc}
    h & 0 \\
    0 & 1/h
\end{array}\right).
$$
The additional action of the stabilizer of $B_3$ leads to the following action on the reduced phase space
$$
q_t\rightarrow \frac{q_t}{h^2},\quad p_t\rightarrow h^2p_t,\quad q_3\rightarrow h^2q_3,\quad p_3\rightarrow \frac{p_3}{h^2},
$$
which is Hamiltonian with the following moment map 
$$
I = q_3p_3-q_tp_t.
$$
Using the symplectic change of coordinates
\begin{equation}\label{sect5:action-}
    q_3=e^{\phi},\quad q_t = e^{-\phi} u,\quad p_3 = e^{-\phi}(I+uv),\quad p_t=e^{\phi}v,\quad \d p_3\wedge \d q_3 + \d p_t \wedge \d q_t = \d I\wedge \d \phi + \d v\wedge \d u,
\end{equation}
and fixing the level set of moment map $I=I_0=\operatorname{const}$ we reduce to the system with one degree of freedom
\begin{equation}
   H= \left( uv-2\theta_{{t}} \right) v \left( uv+I_0 \right) - 2\left( uv-\theta_{{t}} \right)  \left( t\theta_{{3}}+\theta_{{2}}
 \right) +  2\theta_{{3}} u.
\end{equation}
Finally, using the change of variables
$$
u = x(xy-I_0), \quad v=\frac{1}{x},\quad \d v\wedge \d u = \d y \wedge \d x
$$
sends Hamiltonian to the Okamoto form of $P_{\operatorname{IV}}$
\begin{equation}
    H = 2y x^{2}+ \left( \theta_{{3}}{y}^{2}+ \left( -2t\theta_{{3}}-2\,
\theta_{{2}} \right) y-2 I_{0}\right) x+ \left( -I_{0}\theta_{{3}}-2\theta
_{{3}}\theta_{{t}} \right) y.
\end{equation}
Taking 
$$
\theta_3=-1,\quad \theta_2=0,\quad I_0=-\theta_0,\quad \theta_t = -\frac{1}{2}(\theta_{\infty}+\theta_0)
$$
we obtain the  $P_{\operatorname{IV}}$ Hamiltonian
$$
H= 2 y x^2 - (y^2+2ty+2\theta_0)x+\theta_{\infty}y.
$$

\subsection{Painlev\'e III}

The connection takes form
\begin{equation}
    A = \frac{B_0}{\lambda}+t\frac{B_1}{\lambda^2} + C
\end{equation}
with deformation one-form
\begin{equation}
    \Omega = -\frac{B_1}{\lambda}\d t.
\end{equation}
the Poisson algebra is 
\begin{equation}
    \{C\ocomma C\}=\{C\ocomma B_{0,1}\}=\{B_1\ocomma B_1\} = 0,\quad \{B_0\ocomma B_0\}=[{\Pi},1\otimes B_0],\quad \{B_0\ocomma B_1\} = [{\Pi},1\otimes B_1]
\end{equation}
the Hamiltonian is given by
\begin{equation}
    H = \frac{1}{2}\underset{\lambda=0}{\res} \tr \frac{\lambda}{t}A^2 = \tr\left(CB_1+\frac{B_0^2}{2t}\right).
\end{equation}
In the case of $\mathfrak{sl}_2$, choosing the gauge such that $C$ is diagonal, we have the following Darboux parametrisation
\begin{multline}
    B_0 =\left( \begin {array}{cc} p_{1}q_{1}-p_{2}q_{2}+\theta_{1}&
-p_{1}q_{1}^{2}+ \left( 2q_{1}q_{2}+1 \right) p_{2}-2
\theta_{1}q_{1}\\ p_{1}&-p_{1}q_{1}+p_{2}q_{2}-\theta_{1}\end {array} \right) \\
B_1 = \left( \begin {array}{cc} 2q_{1}q_{2}\theta_{2}+\theta_{2}&
-2\theta_{2} \left( q_{1}q_{2}+1 \right) q_{1}
\\ 2\theta_{2}q_{2}&-2q_{1}q_{2}\theta_{2}-\theta_{2}\end {array} \right),\quad 
C = \left(\begin{array}{cc}
        \theta_3 & 0 \\
        0 & -\theta_3
    \end{array}\right).
\end{multline}
the Hamiltonian writes as
\begin{equation}
    tH = p_{2}^{2}q_{2}^{2}+4t\theta_{2}\theta_{3}q_{1}q_{2}-2\theta_{1}p_{2}q_{2}+p_{1}p_{2}. 
\end{equation}
The action of the stabilizer of $C$  gives the following integral of motion
$$
I = q_1p_1-q_2p_2.
$$
In order to reduce the number of degrees of freedom, we use the change of variables 
$$
q_1=e^{\phi},\quad q_2 = -e^{-\phi} u,\quad p_1 = e^{-\phi}(I+uv),\quad p_2=-e^{\phi}v,\quad \d p_1\wedge \d q_1 + \d p_2 \wedge \d q_2 = \d I\wedge \d \phi + \d v\wedge \d u
$$
which leads to the following Hamiltonian
\begin{equation}\label{eq:P3c}
    tH = v^{2}u^{2} - \left( v^{2}+2\theta_{1}v+4t\theta_{2}\theta_{3} \right) u-I_0v
\end{equation}
where $I_0$ is given value of the first integral $I$. The obtained Hamiltonian corresponds to the Painlev\'e III equation of type $D_6$ after some choice of constants. To obtain further degenerations to $D_7$ and $D_8$ we have to consider nilpotent orbits.

\subsection{Painlev\'e II. Jimbo-Miwa}\label{suse:PIIJM}
The connection takes form
\begin{equation}
    A(\lambda) = \frac{B_3}{\lambda^4}+\frac{B_2}{\lambda^3}+\frac{B_1+tB_3}{\lambda^2}+\frac{B_0}{\lambda}.
\end{equation}
The deformation one form is
\begin{equation}
    \Omega = -\frac{B_3}{\lambda}\d t
\end{equation}
The deformation equations are
\begin{equation}
    \frac{d}{d t}B_3 = [B_2, B_3],\quad \frac{d}{dt}B_2 = [B_1, B_3],\quad \frac{d}{dt}B_1 = [B_0-tB_2,B_3],\quad \frac{d}{dt}B_0 = 0.
\end{equation}
The Poisson structure is given by
\begin{equation}\label{PS:PII}
\left\{B_i\ocomma B_j\right\} = [{\Pi},\ID \otimes B_{i+j-1}]
\end{equation}
The Hamiltonian takes the form
\begin{equation}
    H = \underset{\lambda=0}{\operatorname{res}}\tr\lambda^3 \frac{A(\lambda)^2}{2}=\tr\left(\frac{B_1^2}{2}+B_0B_2+tB_1B_3\right),
\end{equation}
where we drop the term $\tr B_3^2$ because it is a Casimir. Since we assume that for Painlev\'e II there is no singularity at $\infty$, the value of the gauge group moment map should be put to zero, i.e.
\begin{equation}
    B_0 = 0.
\end{equation}
Such reduction has to be viewed as a Hamiltonian reduction written on the co-adjoint orbit of the Takiff algebra $\hat{\mathfrak{g}}_3$, so we have to change not only the Hamiltonian, but also the Poisson structure. However, usually the second Painlev\'e equation isomonodromic problem writes in a chart where the only singularity is at $\infty$. In this case, the connection  takes form
\begin{equation}\label{eq:nahm}
    A(\lambda) = B_3\lambda^2 + B_2\lambda + B_3 t + B_1.
\end{equation}
Here we already resolved the gauge group moment map, by setting the residue at $\infty$ to be zero. The deformation one-form then my be written as
$$
\Omega = \left(B_3\lambda + B_2 \right)\d t.
$$
The deformation equations are
$$
\dot{B}_3=0,\quad \dot{B}_2=[B_3,B_1],\quad \dot{B}_1 = t[B_2,B_3]+[B_2,B_1].
$$
The deformation equations are Hamiltonian, with Hamiltonian written as
\begin{equation}
H = \underset{\lambda=0}{\res}\tr \frac{A^2}{2\lambda} = \tr\left(\frac{B_1^2}{2}+tB_1B_3\right).
\end{equation}
To obtain the Painlev\'e II equation, we consider the $\mathfrak{sl}_2$ case. the Darboux parametrisation is given by
$$
B_3=\left(\begin{array}{cc}
    \theta_4& 0 \\
    0 & -\theta_4
\end{array}\right),\quad B_2=\left(\begin {array}{cc} \theta_{{3}}&-2\theta_{{4}}q_{{3}}
\\ \noalign{\medskip}2\theta_{{4}}q_{{4}}&-\theta_{{3}}
\end {array} 
\right),
$$
$$
B_1 = \left( \begin{array}{cc}
  2\theta_{{4}}q_{{3}}q_{{4}}+\theta_{{2}
}&-\theta_{{4}}q_{3}^{3}q_{4}^{2}+ \left( \theta_{{3}}-4
\theta_{{4}} \right) q_{{4}}q_{3}^{2}-\theta_{{4}}q_{{3}}+p_{{4}
}\\ \noalign{\medskip}-\theta_{{4}}q_{3}^{2}{q_{{4}}}^{3}+
 \left( \theta_{{3}}-4\theta_{{4}} \right) q_{4}^{2}q_{{3}}+
 \left( 2\theta_{{3}}-\theta_{{4}} \right) q_{{4}}+p_{{3}}&-2
\theta_{{4}}q_{{3}}q_{{4}}-\theta_{{2}}
\end{array}\right).
$$
The Hamiltonian takes form
\begin{multline}
H = -(2\theta_4q_3q_4+\theta_2)^2-2t(2\theta_4q_3q_4+\theta_2)\theta_4-(  (\theta_3-4\theta_4)q_4 q_3^2-\theta_4q_3+p_4\\-\theta_4 q_3^3 q_4^2)((\theta_3-4\theta_4)q_4^2q_3+(2\theta_3-\theta_4)q_4+p_3-\theta_4q_3^2q_4^3).
\end{multline}
The action of stabilizer of $B_4$ gives us the moment map
$$
I=p_3q_3-p_4q_4
$$
which gives us the following change of variables $(p_3,p_4,q_3,q_4)\to (I,v,\phi,u)$
$$
p_3 = e^{-\phi}(I+uv),\quad p_4 = e^{\phi}v,\quad q_3 = e^{\phi},\quad q_4 = e^{-\phi}u,\quad \d p_3\wedge \d q_3+\d p_4\wedge \d q_4 = \d v\wedge\d u + \d I\wedge \d \phi.
$$
The Hamiltonian then writes as
$$
H=-(2\theta_4 u+\theta_2)^2-2t(2\theta_4 u+\theta_2)\theta_4-(v-\theta_4 u^2+(\theta_3-4\theta_4)u-\theta_4)(uv+(\theta_3-4\theta_4)u^2+(2\theta_3-\theta_4)u+I-\theta_4 u^3)
$$
The change of variable
$$
v = w+\frac{1}{2u}(2\theta_4 u^3-2u^2\theta_3 +8\theta_4 u^2-2u\theta_3 +2\theta_4 u-I),\quad w=-\frac{p}{q},\quad u=-\frac{q^2}{2}
$$
gives us
$$
H = \frac{p^2}{2} -\theta_4^2 q^4+\left(2\theta_4 ^2 t+2\theta_2 \theta_4 -\frac{\theta_3 ^2}{2}\right)q^2-\frac{I^2}{2q^2}
$$
which is the Hamiltonian of $P_{34}$ equation, which is equivalent to Painlev\'e II in case when $I=0$.

\begin{remark}
The isomonodromic problem with connection matrix \eqref{eq:nahm} corresponds to the non-autonomous version of the famous Nahm top which first appeared in \cite{Nahm}. Treating the variable $t$ as a constant, we obtain an integrable system with Lax matrix \eqref{eq:nahm} which is gauge equivalent to the Lax matrix for the Nahm equation. This gives the explicit Hamiltonian formulation of the Nahm equation in terms of the coadjoint orbits of the Takiff algebras. This should coincide with the Hamiltonian formalism for the Nahm equations introduced in \cite{Sask}.
\end{remark}
\section{Quantisation}\label{se:quant}

In this section we give a general formula for the confluent KZ equations with singularities of arbitrary Poincar\'e rank in any dimension. 
Moreover,
 we use the lifted Darboux coordinates in order to generalise an observation by Reshetikhin that the quasiclassical solution of the standard  KZ equations (i.e. with simple poles) is expressed via the isomonodromic $\tau$-function \cite{Res}. Here we propose an easy proof which is valid for any configuration of the points of the divisor on the Riemann sphere. Firstly we review Reshetikhin's approach for the quantum isomonodromic problems and then produce our proof which is based on the generalisation of an observation by Malgrange \cite{Mal1}.
 
Throughout this section we work with the canonical quantisation of the linear Poisson brackets that prescribes the standard correspondence principle
\begin{equation}\label{eq:quant}
\{f,g\} \xrightarrow{\quad\quad} [\hat{f},\hat{g}]= i \hbar \widehat{\{f,g\}},
\end{equation}
where the symbol $\widehat{}$ denotes the quantum operator, i.e. $\widehat f$ is the quantum operator corresponding to the classical function $f$, and $\hbar$ is a formal deformation parameter. More accurately, one can speak about the so--called {\it Rees deformation} that assigns to a filtered vector space 
$R=\cup_{i} R_i$ a canonical deformation of its associated graded algebra ${\rm gr}(R)$ over the affine line $\mathbb A_1$ considered as the spectrum ${\rm Spec}(\mathbb C[\hbar])$ of the polynomials $\mathbb C[\hbar]$. The fiber at the point $\hbar$ is isomorphic to $R$ if $\hbar\neq 0$ and to ${\rm gr}(R)$  for $\hbar= 0$. The corresponding $\mathbb C[\hbar ]$--module here is the direct sum $\oplus_i R_i$ on which $\hbar$  acts by mapping each $R_i$ to $R_{i+1}$ \cite{Eisen}.
In our case the Rees construction gives a one-parameter family of algebras $U_{\hbar}(\mathfrak{g})$, with the associated graded algebra $U_{0}(\mathfrak{g})$ being the symmetric algebra $S(\mathfrak g)$. The $\hbar$ deformation re-scales the bracket by $\hbar$, so that the $\hbar$ linear terms define the standard Poisson bracket on $S(\mathfrak g)$. 

\subsection{Finite-dimensional representation}
In this sub-section, we recall the basic ideas at the basis of Reshetikhin's approach to quantum isomonodromic problems for Fuchsian systems
and then adapt it to the irregular case. We fix $\hbar=1$ for simplicity. 

In the case of  Fuchsian systems, we are dealing with the canonical quantisation of the direct product of the co-algebras $\mathfrak{g}^{\star}$.
The quantisation functor sends the functions on the phase space of the classical system to the differential operators which act on some Hilbert space in a way that \eqref{eq:quant} holds. In principle, a choice of finite dimensional representation may be seen as a choice of the special subspace of the Hilbert space of functions on which the algebra of quantum  operators acts. However,  we may avoid such complicated construction of finite dimensional representation when the classical Poisson algebra is given by a linear Poisson bracket. Indeed, for $\mathfrak{g}^{\star}$, the \KKS bracket endows the space of functions with the structure of a Lie algebra so that the structure constants of this Poisson algebra are identified with the structure constants of the Lie algebra $\mathfrak{g}$.

In general, the quantisation procedure for the phase space of the Fuchsian system may be viewed as a map from
$$
\underbrace{\mathfrak{g}^{\star}\times\mathfrak{g}^{\star}\times\dots \times \mathfrak{g}^{\star}}_{n}
$$
to the differential operators which act on the  tensor product of Hilbert spaces $\mathcal{H}_i$:
$$
\mathcal{H}_1 \otimes\mathcal{H}_2\otimes\dots \otimes \mathcal{H}_n.
$$
However, the isomonodromic nature of the Hamiltonian systems we consider gives additional information which may be used to define a quantum problem in a uniform way. Following \cite{JNS}, we  quantise the connection that becomes the generating function for the quantum Hamiltonians. Considering the connection as a matrix whose entries are functions on  $\mathfrak{g}^{\star}\times\mathfrak{g}^{\star}\times \dots \times \mathfrak{g}^{\star}$, we obtain the following quantisation for the Fuchsian case:
$$
\hat{A}(\lambda) = \sum\limits_{i=1}^n\frac{\hat{A}^{(i)}}{\lambda - u_i},\quad \hat{A}^{(i)} = \sum\limits_{\alpha} e_{\alpha}^{(0)}\otimes e_{\alpha}^{(i)},\quad e_{\alpha}^{(i)} = 1  \otimes\dots\otimes \underset{i}{e_{\alpha}}\otimes\dots\otimes 1,
$$
where each $e^{(i)}_{\alpha}$, $i=1,\dots,n$ is a basis of the representation we choose for the quantisation and the first $e_\alpha^{(0)}$ corresponds to the auxiliary space $\mathcal H_0$ given by the connection. The Schlesinger Hamiltonians then transform to
\begin{equation}
    \hat{H}_i = \sum\limits_{j \neq i }\frac{\tr^{(0)}(\hat{A}^{(i)}\hat{A}^{(j)})}{u_i-u_j},
\end{equation}
where $\tr^{(0)}$ is the trace in the auxiliary space  $\mathcal{H}_0$.
The quantum Schlesinger Hamiltonians $ \hat{H}_i $ are the solutions for the classical Yang-Baxter equations and may be written as
\begin{equation}
\hat{H}_i = \sum\limits_{j \neq i }\frac{r_{ij}}{u_i-u_j},
\end{equation}
where $r_{ij}$ is a solution of the classical Yang-Baxter equation
$$
[r_{ij},r_{ik}]+[r_{ij},r_{jk}]+[r_{ik},r_{jk}]=0.
$$
The corresponding set of  Schr\"odinger equations are called \KZ equations and take form
$$
\nabla_i\Psi =\left(\frac{\partial}{\partial u_i} - \sum\limits_{j \neq i }\frac{r_{ij}}{u_i-u_j}\right)\Psi = 0.
$$
Moreover, the \KZ operators commute, i.e. 
$$
[\nabla_i,\nabla_j] = 0 \quad \Longleftrightarrow \quad 
    \frac{\partial}{\partial u_i}\hat{H}_j=\frac{\partial}{\partial u_j}\hat{H}_i , \quad
    [\hat{H}_i,\hat{H}_j]=0.
$$
Reproducing the same scheme for the Takiff co-algebras, we obtain the  quantisation map that acts by replacing the co-algebra with the Lie algebra
\begin{equation}
    \hat{\mathfrak{g}}_{r_1}^{\star}\times \hat{\mathfrak{g}}_{r_2}^{\star}\times\dots \hat{\mathfrak{g}}_{r_n}^{\star}\times \hat{\mathfrak{g}}_{r_\infty}^{\star}  \quad \xrightarrow{\quad\quad}\quad \hat{\mathfrak{g}}_{r_1}\otimes\dots \hat{\mathfrak{g}}_{r_n}\otimes \hat{\mathfrak{g}}_{r_\infty}.
\end{equation}
The quantum connection then takes the form
$$
\hat{A}(\lambda) = \sum\limits_{i}^n\left(\sum\limits_{j=0}^{r_i}\frac{\hat{B}^{(i)}_j\left(t^{(i)}_1,t^{(i)}_2\dots t^{(i)}_{r_i}\right)}{(\lambda-u_i)^{j+1}}\right),\
$$
where $\hat{B}^{(i)}$'s are given by
$$
\hat{B}^{(i)}_j(t^{(i)}_1,\dots t^{(i)}_{r_i}) = \sum\limits_{k=j}^r \hat{A}^{(i)}_k\mathcal{M}^{(r_i)}_{j,k}(t^{(i)}_1,t^{(i)}_2\dots t^{(i)}_{r_i}),\quad \hat{A}_k^{(i)} = \sum\limits_{\alpha}e^{(0)}_{\alpha}\otimes e_{\alpha}^{(i)}\otimes z_i^k,\quad  e_{\alpha}^{(i)} = 1  \otimes\dots\otimes \underset{i}{e_{\alpha}}\otimes\dots\otimes 1.
$$
The Hamiltonians which correspond to the position of poles are given as in the Fuchsian case, i.e.
$$
\hat{H}_{u_i} = \frac{1}{2}\underset{\lambda=u_i}{\res}\tr_0 \hat{A}(\lambda)^2,
$$
where $\tr_0$ is the trace in the $0$-th space, so we now have to choose a quantum ordering, for example lexigraphical ordering.
The irregular Hamiltonians have to be calculated according to the Theorem \ref{th:mainHam} at each irregular singularity changing $\tr$ by $\tr_0$. Again we will choose a quantum ordering. Thus, we obtain that the irregular Hamiltonians are given by
$$
\mathcal{M}^{(r_i)}\left(\begin{array}{c}
    \hat{H}^{(i)}_1 \\
    \hat{H}^{(i)}_2 \\
    \dots \\
    \hat{H}^{(i)}_{r_i}
\end{array}\right) = \left(\begin{array}{c}
    \hat{S}^{(u_i)}_1 \\
    \hat{S}^{(u_i)}_2 \\
    \dots \\
    \hat{S}^{(u_i)}_{r_i}
\end{array}\right),\quad \hat{S}_k^{(u_i)} = \frac{1}{2}\oint\limits_{\Gamma_{u_i}}(\lambda-u_i)^k\tr_0 \hat{A}^2 \d \lambda
$$
at the point $u_i$ with the Poincare rank $r_i$. To prove 
 Theorem \ref{th:mainKZ} we need to show that the confluent KZ gives a quantum integrable system, namely that the differential operators defined in \eqref{eq:KZopu}, \eqref{eq:KZopt}
$$
\nabla_{u_j}:=\frac{\partial}{\partial u_j}-\widehat H_{u_j},\quad j=1,\dots,n
$$
$$
\nabla_{k}^{(i)}:=\frac{\partial}{\partial t^{(i)}_k}-\widehat H_{k}^{(i)},\quad i=1,\dots,n,\infty, \quad k=1,\dots,r_i
$$
commute. This is a simple consequence of the fact that in the quantisation process the derivatives remain commutative, i.e. for example $[\frac{\partial}{\partial u_j},\frac{\partial}{\partial t^{(i)}_k}]=0$, and that the quantum Hamiltonians are linear combinations of the quantum Gaudin spectral invariants $\hat{S}^{(u_i)}_k$, $k=0,\dots,r_i$, which commute as proved in \cite{MTV}.
 
We have to mention that for the Fuchsian times the isomonodromic Hamiltonian depends on each phase space $\mathfrak{g}_{r_i}$ linearly -- which means that it may be written as
$$
\hat{H}_{u_i} \in \hat{\mathfrak{g}}_{r_1}\otimes\dots \hat{\mathfrak{g}}_{r_n}\otimes \hat{\mathfrak{g}}_{r_\infty} \subset U\left(\hat{\mathfrak{g}}_{r_1}\oplus\dots  \oplus \hat{\mathfrak{g}}_{r_\infty}\right).
$$
In the case of irregular poles, the Hamiltonians become more complicated -- there are quadratic terms which contain elements from the same space and in general we have that
$$
\hat{H}^{(i)}_k \in U(\hat{\mathfrak{g}}_{r_1})\otimes\dots U(\hat{\mathfrak{g}}_{r_n})\otimes U( \hat{\mathfrak{g}}_{r_\infty}).
$$
The problem of calculating the explicit form of the Hamiltonians introduced in this paper involves the  representation theory of $U(\hat{g}_{r_i})$, that is rather complicated. In order to avoid this representational theoretic problems, we write down the quantum Hamiltonians for the irregular isomonodromic deformations using the intermediate Darboux coordinates. We deal with  the classical examples of the Painlev\'e equations in the next section, where we provide invariant subspaces for these Hamiltonians. These subspaces give finite dimensional representations for the Hamiltonians which are the quantum reduction of the irregular Hamiltonians introduced in this section.

\subsection{Intermediate Quantum Hamiltonians for Painlev\'e equations.}
In this subsection we write quantum Hamiltonians for the Painlev\'e equations in Darboux coordinates before the reduction with respect to the gauge group action. In the case of Painlev\'e VI, the gauge group action is not taken into account. For the other cases, we partly resolve the gauge group action by diagonalising the leading term, but we do not finish reduction by ignoring the additional Cartan torus action (otherwise the quantisation becomes very cumbersome). Because of this, in the Painlev\'e VI example the number of coordinates for $\mathfrak{sl}_2$ for $4$ punctures is $6$ while in the other examples the number of intermidiate coordinates is $4$ ($2$ moments $+$ $2$ positions). Since we are dealing with Darboux coordinates, the quantisation process becomes fairly straightforward. In this subsection, we show that for each of the non-ramified Painlev\'e differential equations, there is a choice of quantisation  such that the quantum operator acts nicely on the space of homogeneous polynomials. More precisely, we show that the invariant subspaces for the  quantum Hamiltonians are given by the homogeneous polynomials in several variables ($3$ for Painlev\'e VI and $2$ for others) with fixed degree. In this section we keep $\hbar$ explicit as that makes it clearer how to extract semi-classical limits.

\subsubsection{Painlev\'e VI}
For the $\mathfrak{sl}_2$ Fuchsian system we have that the Hamiltonians in the intermediate coordinates take form
\begin{equation}
\begin{aligned}
H_i = \sum\limits_{j\neq i}\frac{h_{ij}}{u_i-u_j},\quad  h_{ij}=2 p_{i} p_{j} q_{i} q_{j}- & p_{i}^{2} q_{i} q_{j}-p_{j}^{2} q_{i} q_{j}- \\ & 2 \theta_{j} p_{i} q_{i}-2 \theta_{i} p_{j} q_{j}+2\theta_{i} p_{i} q_{j} +2\theta_{j} p_{j} q_{i}+2 \theta_{i} \theta_{j}
\end{aligned}
\end{equation}
The quantisation problem is not trivial because we have to choose the quantum ordering. There are three standard ways of the ordering, which are given by
$$
:\widehat{p_i} \widehat{q_j}: = :\widehat{q_j} \widehat{p_i}: = \widehat{q_j}\widehat{p_i} + \delta_{ij}\varepsilon^{(i)},\quad 
\varepsilon^{(i)} = \left\{\begin{array}{ll}
    0,\, &  \text{left,}\\
    i\hbar,\,& \text{right,} \\
    \frac{i\hbar}{2},\, & \text{Weyl.}
\end{array}\right.
$$
This leads to the following forms of Hamiltonians
$$
\hat{H_i}=\sum\limits_{j\neq i}\frac{\hat{h}_{ij}}{u_i-u_j}
$$
where
$$
\begin{aligned}
  \hat{h}_{ij} = 2  \hat{q}_{i} \hat{q}_{j}\hat{p}_{i} \hat{p}_{j}- \hat{q}_{i} \hat{q}_{j}\hat{p}_{i}^{2}- \hat{q}_{i} \hat{q}_{j}\hat{p}_{j}^{2}-2 (\theta_{j}-\varepsilon^{(j)})  \hat{q}_{i}\hat{p}_{i} - 2 (\theta_{i}-\varepsilon^{(i)}) \hat{q}_{j}\hat{p}_{j}& + 2(\theta_{i}-\varepsilon^{(i)})  \hat{q}_{j}\hat{p}_{i}+ \\ +2(\theta_{j}-\varepsilon^{(j)})  \hat{q}_{i}\hat{p}_{j} &
 + 2(\theta_{i}-\varepsilon^{(i)}) (\theta_{j}-\varepsilon^{(j)}).
 \end{aligned}
$$
Here we see that different choices of the ordering lead to different shifts of the local monodromies $\theta_i\rightarrow\theta_i-\varepsilon^{(i)}$.  Thanks to this fact, and the fact that the shifts are of order $\hbar$, we may fix the left ordering without loss of generality.

The most remarkable property is that the Hamiltonians $\hat{H}_i$ leave invariant the space of homogeneous polynomials of $q_i$ with fixed degree in the following choice of the quantisation $\hat{p}_i=-i\hbar\frac{\partial}{\partial x_i}\cdot$, $\hat{q}_i=x_i\cdot$. So we may look for a solutions for the set of quantum Schrodinger equations
\begin{equation}\label{eq:schr}
i\hbar \partial_{u_i}\Psi = \hat{H}_i\Psi
\end{equation}
in the following form
$$
\Psi^{(n)} = \sum_{|\alpha|= n} w_{\alpha}(u_1,..,u_i..,u_m) \prod_{i=1}^m x_i^{\alpha_k},\quad |\alpha| = \sum\limits_{i=1}^m\alpha_i 
$$
which will lead to the non-autonomous linear system of ODE for the $w_{\alpha}(\bf{u})$-s. The resulting equations in fact are KZ equations, since the equations for $w_\alpha$ inherit the singularities of $\hat{H}_i$. Let's consider the vector
$$
W^{(n)} = \left(\begin{array}{c}
    w_{\alpha_1} \\
    w_{\alpha_2} \\
    .. \\
    .. \\
    w_{\alpha_N}
\end{array}\right)
$$
where $\alpha_i$ are the distinct partitions of $n$ with height $m$ (with zero entries). Then $W^{(n)}$ satisfies the equations
$$
i\hbar \frac{\partial}{\partial u_i}W^{(n)}  - \sum\limits_{j\neq i}\frac{M_n^{(i,j)}}{u_i-u_j} W^{(n)} = 0
$$
where $M_n^{(i,j)}$ is the action of $\hat{h}_{ij}$ on homogeneous polynomials of degree $n$. These equations are \KZt equations. 

In the case of the Painlev\'e VI equation, we deal with 4-punctured sphere $0,1,t,\infty$. The quantum Hamiltonian then writes as
\begin{equation}
\begin{aligned}
   \hat{H} = & \frac{1}{t}\left( 2  \hat{q}_{1} \hat{q}_{2}\hat{p}_{1} \hat{p}_{2}- \hat{q}_{1} \hat{q}_{2}\hat{p}_{1}^{2}- \hat{q}_{1} \hat{q}_{2}\hat{p}_{2}^{2}-2 \theta_{2}  \hat{q}_{1}\hat{p}_{1} - 2 \theta_{1} \hat{q}_{2}\hat{p}_{2} + 2\theta_{1}  \hat{q}_{2}\hat{p}_{1}+  +2\theta_{2}  \hat{q}_{1}\hat{p}_{2} + 2\theta_{1} \theta_{2}\right) + \\
   & \frac{1}{t-1}\left( 2  \hat{q}_{1} \hat{q}_{3}\hat{p}_{1} \hat{p}_{3}- \hat{q}_{1} \hat{q}_{3}\hat{p}_{1}^{2}- \hat{q}_{1} \hat{q}_{3}\hat{p}_{3}^{2}-2 \theta_{3}  \hat{q}_{1}\hat{p}_{1} - 2 \theta_{1} \hat{q}_{3}\hat{p}_{3} + 2\theta_{1}  \hat{q}_{3}\hat{p}_{1}+  +2\theta_{3}  \hat{q}_{1}\hat{p}_{3} + 2\theta_{1} \theta_{3}\right).
\end{aligned}
\end{equation}
Let's consider simple case where $|\alpha|=1$. Substituting the following function
$$
\Psi^{(1)} = w_1x_1+w_2x_2+w_3x_3
$$
into the Schrodinger equation \eqref{eq:schr} gives the following system
\begin{equation}
\begin{gathered}
i\hbar\frac{\mathrm{d}}{\mathrm{d} t} w_{1} = \frac{2 i \hbar \theta_{2}(w_2-w_1)+2 \theta_{1} \theta_{2} w_{1}}{t}+\frac{2 i \hbar \theta_{3}(w_3-w_1)+2 \theta_{1} \theta_{3} w_{1}}{t-1} \\
i \hbar \frac{\mathrm{d}}{\mathrm{d} t} w_{2} = \frac{2 i \hbar  \theta_{1}(w_1- w_{2})+2 \theta_{1} \theta_{2} w_{2} }{t}+\frac{2 \theta_{1} \theta_{3} w_{2} }{t-1} \\ 
i  \hbar\frac{\mathrm{d}}{\mathrm{d} t} w_{3}=\frac{2 \theta_{1} \theta_{2} w_{3} }{t}+\frac{2 i  \hbar  \theta_{1}(w_1-w_{3})+2 w_{3}  \theta_{1} \theta_{3}}{t-1}
\end{gathered}
\end{equation}
whose solution is given by the hypergeometric function in the following way
$$
\begin{aligned}
w_1=& C_1 t^{-\frac{2 i}{\hbar}\theta_1 \theta_2}  (t-1)^{-\frac{2 i}{\hbar} \theta_1 \theta_3}+ C_2 t^{-\frac{2 i}{\hbar} \theta_1 \theta_2} (t-1)^{-\frac{2 i}{\hbar} \theta_1 \theta_3-2(\theta_1+\theta_3)}{}_2F_1(2\theta_2, -2\theta_3+1;2(\theta_1+\theta_2)+1; t) \\ 
 & C_3  t^{-\frac{2 i}{\hbar} \theta_1 \theta_2-2(\theta_1+\theta_2)} (t-1)^{-\frac{2 i}{\hbar} \theta_1 \theta_3-2(\theta_1+\theta_3)}{}_2F_1(-2\theta_1, -2(\theta_1+\theta_2+\theta_3)+1;-2(\theta_1+\theta_2)+1; t).
\end{aligned}
$$

\subsubsection{Painlev\'e V} the Hamiltonian in the intermediate coordinates is given by
\begin{equation}
    tH = 2 t \theta_{\infty} q_{1} p_1 - q_{0} q_{1}p_{0}^{2}+2q_0 q_1 p_{0} p_{1}  - q_{0} q_{1}p_{1}^{2}+2 \theta_0 q_{1}p_{0} +2  \theta_{1} q_0p_{1} +2 \theta_{0} \theta_{1}.
\end{equation}
Using the same argument as in the previous case, we consider left ordering. Moreover, we see that if quantise in the following way
\begin{equation}\label{eq:stQ}
\hat{q}_i = x_i\cdot ,\quad \hat{p}_i= -i\hbar \frac{\partial }{\partial x_i},
\end{equation}
the space of homogeneous polynomials in $x_0$ and $x_1$ is invariant under the action of the Hamiltonian. Considering the example of the degree 2
$$
\Psi^{(2)} = w_1 x_1^2+w_2 x_0^2+w_3 x_0x_1,
$$
we get the following system of ordinary differential equations for the coefficients
\begin{equation}
\begin{gathered}
i t \hbar \frac{d}{d t}w_{1} =-4i t \theta_\infty \hbar   w_{1} -2 \,i\theta_0 \hbar  w_{3} \\ 
i t \hbar  \frac{d}{d t}w_{2} = -2 i  \theta_{1} \hbar w_{3}
\\
i t \hbar  \frac{d}{d t}w_{3}  =2 \hbar^{2}( w_{1}+w_{2}) - 2 i \theta_{\infty} t \hbar w_{3}  - 4 \theta_1 i\hbar  w_{1} - 4 \theta_0 i  \hbar w_{2} \theta_{0}-2 \hbar^{2} w_{3}.
\end{gathered}
\end{equation}
These linear equations can be solved explicitly in terms of modified Bessel functions.

\subsubsection{Painlev\'e IV}
the Hamiltonian in the intermediate coordinates takes form
\begin{equation}
H =   q_{{t}}p_{{t}}^2p_{{3}}-2\theta_{{t}} p_{{t}}p_{{3}}-2 \left( p_{{t}}q_{{t}}-\theta_{{t}} \right)  \left(t\theta_{{3}}+\theta_{{2}} \right) +2\theta_{{3}}q_{{3}}q_{{t}}.
\end{equation}
In general, the choice of the Lagrangian submanifold for the quantisation procedure defines the properties of the quantum Hamiltonian. Here the quantum Hamiltonian will not preserve the homogeneous polynomials if we choose the standard quantisation \eqref{eq:stQ}. However, the choice of the Lagrangian sub-manifold is irrelevant when we deal with the Darboux coordinates and corresponds to the integral transformation on the quantum level. If we choose the following quantisation
$$
 \hat{q}_{{3}}=x\cdot,\quad \hat{p}_3=\hbar\frac{\partial}{\partial x},\quad\hat{q}_t=\hbar\frac{\partial}{\partial y},\quad \hat{p}_t = y\cdot
$$
the quantum Hamiltonian will preserve degree of the homogeneous polynomials. Moreover the choice of the ordering shifts the monodromy parameter $\theta_t$ by $\hbar$-small values. the Hamiltonian writes as
\begin{equation}
    \hat{H} = y^2\frac{\partial^2}{\partial x\partial y}-2\theta_t y\frac{\partial}{\partial x} - 2\left(t\theta_{{3}}+\theta_{{2}} \right)\left(y\frac{\partial}{\partial y}-\theta_t\right)+2\theta_3x\frac{\partial}{\partial y}.
\end{equation}
Writing down the system for the second order polynomial wave function
$$
\Psi^{(2)}=w_1x^2+w_2y^2+w_3xy,
$$
we obtain the system
\begin{equation}
    \frac{i\hbar}{2}\frac{d}{dt}\left(\begin{array}{c}
        w_1  \\
        w_2  \\
        w_3
    \end{array}\right) = \left(\begin{array}{ccc}
        -(t  \theta_3 +\theta_2)   \theta_t & 0 & -\theta_3 \\
         & -(t  \theta_3 +\theta_2)   \theta_t  & 2\theta_t-1\\
         2\theta_t & -2\theta_3 & -(t  \theta_3 +\theta_2)   \theta_t 
    \end{array}\right)\left(\begin{array}{c}
        w_1  \\
        w_2  \\
        w_3
    \end{array}\right),
\end{equation}
which may be solved via exponential functions.

\subsubsection{Painlev\'e III}
the Hamiltonian is
$$
tH = p_{2}^{2}q_{2}^{2}+4t\theta_{2}\theta_{3}q_{1}q_{2}-2\theta_{1}p_{2}q_{2}+p_{1}p_{2}. 
$$
We quantise as follows:
$$
\hat{q}_1 = x\cdot,\quad \hat{p}_1= i\hbar \frac{\partial}{\partial x},\quad \hat{q}_2= i \hbar\frac{\partial}{\partial y},\quad \hat{p}_2 =y\cdot ,
$$
leading to the quantum Hamiltonian (up to $\hbar$ shifts of $\theta_1$) takes form
\begin{equation}\label{eq:P3q}
    \hat{H}=y^2 \frac{\partial^2}{\partial y^2} -2\theta_1y \frac{\partial}{\partial y}+4t\theta_2x \frac{\partial}{\partial y}+y \frac{\partial}{\partial x}.
\end{equation}
Writing down system for the second order polynomial wave function
$$
\Psi^{(2)}=w_1x^2+w_2y^2+w_3xy,
$$
we obtain the system
\begin{equation}
    i\hbar t \frac{d}{dt}\left(\begin{array}{c}
        w_1  \\
        w_2  \\
        w_3
    \end{array}\right) = \left(\begin{array}{ccc}
        0  & 0 & 4t \\
        0 & 2-   4\theta_1  & 1\\
         2& 8t & -2  \theta_1
    \end{array}\right)\left(\begin{array}{c}
        w_1  \\
        w_2  \\
        w_3
    \end{array}\right).
\end{equation}
These equations can be solved in terms of confluent hypergeometric functions of type $_{1}F_{2}$.

\subsubsection{Painlev\'e II}
The intermediate Darboux coordinates Hamiltonian is
\begin{multline}
   H= (q_3 ^5 q_4 ^5 + 8 q_3 ^4 q_4 ^4 + 18 q_3 ^3 q_4 ^3 + 12 q_3 ^2 q_4 ^2 + (4 t+1) q_3  q_4  ) \theta_4 ^2 +\\ + (-2 q_3 ^4 q_4 ^4 - 10 q_3 ^3 q_4 ^3 - 10 q_3 ^2 q_4 ^2 - 2 q_3  q_4 ) \theta_3 \theta_4  -\\- (p_3  q_3 ^3 q_4 ^2 - p_4  q_3 ^2 q_4 ^3 - 4 p_3  q_3 ^2 q_4  - 4 p_4  q_3  q_4 ^2 + 4 q_3  q_4  \theta_2  + 2 t \theta_2  - p_3  q_3  - p_4  q_4 ) \theta_4 +\\  + (q_3 ^3 q_4 ^3 + 2 q_3 ^2 q_4 ^2) \theta_3  ^2 + (p_3  q_3 ^2 q_4  + p_4  q_3  q_4 ^2 + 2 p_4  q_4 ) \theta_3  + p_3  p_4.
\end{multline}
Choosing of the following quantisation
$$
\hat{q}_3=-i\hbar \frac{\partial}{\partial x},\quad \hat{p}_3 = x\cdot,\quad \hat{q}_4 = y\cdot,\quad \hat{p}_4 = -i\hbar \frac{\partial}{\partial y}
$$
leads to the invariance of the degree of the homogeneous polynomials with respect to the Hamiltonian action. Indeed, after quantisation, the Hamiltonian is mapped  to the operator with the same number of derivatives and multiplications in each member. We do not provide the explicit form of the quantum Hamiltonian and the action on the eigenspaces since the calculation is straightforward but the answer is too long.

\begin{remark}
In these examples, we consider the deformation quantisation of the intermediate Darboux coordinates. This means that the quantised Hamiltonians are  elements of the Weyl algebra in two variables $\mathbb{W}[x,y]=\mathbb{C}[x,\partial_x, y,\partial_y]/\langle [\partial_x,x]=1,[\partial_y,y]=1\rangle$. However, we know that the Hamiltonian we quantise allows additional symmetry, which lifts to an additional vector field $\hat{I}$ that commutes with the quantum Hamiltonian vector field. For example, in the case of the Painlev\'e III equation, the quantum Hamiltonian \eqref{eq:P3q} commutes with
$$
\hat{I} = x\frac{\partial}{\partial x}+y\frac{\partial}{\partial y}.
$$
By restricting to the eigenfunctions of $\hat{I}$ with some chosen eigenvalue $I_0$, we produce quantum Hamiltonian reduction, which is simply given by the quotient of the algebra $\mathbb{W}[x,y]/\langle \hat{I}-I_0\rangle$. As a result we obtain the following quantum Hamiltonian
$$
\hat{H}_{\operatorname{III}}=q^2\frac{\partial}{\partial q} + \left(-q^2 - 2q\theta_1 + 4 t \theta_2\right)\frac{\partial}{\partial q} + I_0 q,\quad q=\frac{y}{x},
$$
which is just the Dirac quantisation of the Hamiltonian for the Painlev\'e III equation \eqref{eq:P3c}. Such reduction may be performed for all examples, the resulting quantum Hamiltonians coincide with the quantum Hamiltonians introduced in \cite{JNS, Nag} up to change of variables and ordering.
\end{remark}

\subsection{Semi-classical solution of the confluent \KZ equation}\label{suse:semi-cl}
In this section we discuss the semi-classical solutions of the confluent \KZ equations in terms of the isomonodromic tau function. By the term ``semi-classical'' we mean the solutions $\Psi_{sc}$ of the system
\begin{equation}\label{eq:KZopu1}
\hbar \frac{\partial\Psi}{\partial u_j}=\widehat H_{u_j}\Psi,\quad j=1,\dots,n
\end{equation}
and
\begin{equation}\label{eq:KZopt1}
\hbar  \frac{\partial\Psi}{\partial t^{(i)}_k}=\widehat H_{k}^{(i)}\Psi,\quad i=1,\dots,n,\infty, \quad k=1,\dots,r_i
\end{equation}
that can be formally expressed as power series in $\hbar$ in an open set in the phase space. To characterise these solutions, we use the lifted Darboux coordinates and quantise them
according to \eqref{eq:quant}
\begin{equation}
 [\widehat P_{i_{ab}},\widehat Q_{j_{cd}}]=i \hbar\, \delta_{ij}\delta_{c b}\delta_{a d}.
\end{equation}
Such quantisation leads to the infinite dimensional representation of the isomonodromic Hamiltonians as  differential operators on a Hilbert space of functions depending on some  coordinates $x_{j_{ab}},\,j=1\dots d$, $a,b = 1\dots m$ and the isomonodromic times. In particular we put
$$
\widehat{Q}_{j_{ab}} = x_{j_{ab}}\cdot ,\quad \widehat{P}_{i_{cd}} =\hbar \frac{\partial\quad\,\,}{\partial x_{i_{dc}}}.
$$
The standard theoretical physics approach (beautifully described for example by Voros in his seminal paper \cite{Vor}) consists in performing WKB analysis of $\Psi_{cs}$ as $\hbar\to 0$.  This of course means paying careful attention to avoid the so-called turning points, or in other words, points in which the action functional expanded in $\hbar$ has zero constant term. Assuming that one stays clear of turning points, the formula for the semiclassical solutions is written as:
$$
\Psi_{sc}\sim \exp\left(\frac{i}{\hbar}{\mathcal S}\right),\quad \hbar\to 0
$$
where $\mathcal{S}$ is the classical action functional which explicitly depends on the entries of the classical variables $Q$ and the isomonodromic times. The dependence of  $\mathcal{S}$ on $P$ is implicit, since
$$
P_{i_{kl}} = \frac{\partial \mathcal S \,\,\,\,}{\partial Q_{i_{lk}}}.
$$
In this section we prove Theorem \ref{th-semicl}, namely that  $\Psi_{{cs}}$ evaluated along solutions of the classical system may be written as the isomonodromic 
$\tau$-function. This statement already appeared in \cite{Res} for the \KZ equations with Fuchsian singularities. However, our approach works also for irregular systems.\\ 

\noindent{\it Proof of Theorem \ref{th-semicl}.} We write the Hamiltonian system with Hamiltonians $H^{(i)}_{u_i}$ and $H^{(i)}_{1},\dots,H^{(i)}_{r_i}$ for $i=1,\dots,n$ in the Darboux coordinates $P_1,P_2\dots P_d,Q_1,Q_2\dots Q_d$. The action functional satisfies the following relation
\begin{equation}\label{eq:actiontau}
\d\mathcal{S} = \sum\limits_{j=1}^dP_j\d Q_j -  \sum\limits_{i}\left(H_{u_i}\d u_i+\sum_{k=1}^{r_i}H^{(i)}_{k}\d t_k^{(i)} \right)=\sum\limits_{j=1}^dP_j\d Q_j -\d \ln (\tau),
\end{equation}
along the solutions of the system. We now use a result by Malgrange to show that the logarithmic differential of the $\tau$ function is already contained in the definition of the action functional: 

\begin{lemma}(Malgrange \cite{Mal1})
If the Hamiltonians are homogeneous polynomials of degree two in $P_1,\dots,P_d$, then along solutions one has 
\begin{equation}\label{eq:Mal_S}
{\rm d}{\mathcal S} =  \sum\limits_{i}\left(H^{(i)}_{u_i}\d u_i+\sum_{k=1}^{r_i}H^{(i)}_{k}\d t_k^{(i)} \right).
\end{equation}
\label{lm:Mal}
\end{lemma}

\proof 
Evaluating the first term in \eqref{eq:actiontau} along the solutions of the isomonodromic deformation equations, we obtain
\[\begin{split}
 \sum_{j}{\rm Tr} (P_j {\rm d}Q_j) =  \sum_{j}{\rm Tr} \left(P_j 
 \sum_l \left(\frac{{\rm d}Q_j}{{\rm d} u_l} {\rm d}u_l +\sum_{k=1}^{r_l} \frac{{\rm d}Q_j}{{\rm d} t^{(l)}_k}{\rm d} t^{(l)}_k\right)\right)=\\
 =
 \sum_{j}{\rm Tr} \left(P_j 
 \sum_l \left(\frac{\partial H_{u_l}}{\partial  P_j} {\rm d}u_l +\sum_{k=1}^{r_l} \frac{\partial H^{(l)}_k}{\partial P_j}{\rm d} t^{(l)}_k\right)\right).
 \end{split}  \]
Using the fact that the Hamiltonians are homogeneous of degree two in $P_1,\dots,P_d$, we obtain that
$$
\sum_{j}{\rm Tr} \left(P_j 
 \sum_l \left(\frac{\partial H_{u_l}}{\partial  P_j} {\rm d}u_l +\sum_{k=1}^{r_l} \frac{\partial H^{(l)}_k}{\partial P_j}{\rm d} t^{(l)}_k\right)\right)=   
    2   \sum\limits_{i}\left(H^{(i)}_{u_i}\d u_i+\sum_{k=1}^{r_i}H^{(i)}_{k}\d t_k^{(i)} \right).
$$
which leads to the statement of the lemma.
\endproof

In the case of the Fuchsian isomonodromic deformation are given by \eqref{eq:LiftedSchles}
$$
H_i=\sum\limits_{j\neq i}\frac{\operatorname{Tr}(Q_iP_iQ_jP_j)}{u_i-u_j},
$$
which are homogeneous of degree $2$ in the entries of matrices $P_1,P_2\dots P_n$. The same holds for the irregular singularities - indeed, the irregular Hamiltonians are given by the quadratic spectral invariants, i.e.
\begin{equation}
    H = \sum\limits_{\alpha,\beta}\sum\limits_{i,j}C_{i,j}^{\alpha,\beta} \tr\left(A^{(\alpha)}_iA^{(\beta)}_k\right),
\end{equation}
where $\alpha,\beta$ are the indices of the singular points, while $i$ and $j$ are the indices which correspond to the coefficients of local expansion near singularity and $C_{i,j}^{\alpha,\beta}$ are coefficients which can be explicitly computed by using the formulas from section 4. Thanks to Lemma \ref{lemma:liftedTakif}, all the terms $\tr\left(A^{(\alpha)}_iA^{(\beta)}_k\right)$ are homogeneous polynomials of degree $2$ in the variables $P_i$ (as well as homogeneous polynomials of degree $2$ in the $Q_i$). This fact allows us to apply Lemma \ref{lm:Mal} to conclude that, up to constant terms, the action functional evaluated on the solutions of the isomonodromic deformation equations coincides with the $\tau$-function. Thanks to this fact, given a solution $(P_1,\dots,P_d,Q_1,\dots,Q_d)$  of the classical isomonodromic deformation equations, the corresponding semi-classical solution admits the following WKB expansion
$$
\Psi_{sc}\sim \exp\left(\frac{i}{\hbar}\log\tau\right),\quad \hbar\to 0
$$
for $u_1,\dots,u_n$, $t_k^{(i)}$, $i=1,\dots,n,\infty$, $k=1,\dots,r_i$ in a poly-disk that does not contain the zeroes of the action functional evaluated along the given solution $(P_1,\dots,P_d,Q_1,\dots,Q_d)$.\hfill{$\square$}\\

Observe that this proof depends on the coordinates we use to quantise. In general, the property of semi-classical solution to be a power of an isomonodromic $\tau$-function breaks for the reduced systems. The  classical analogue of this phenomenon is equivalent to the  statement that the reduced Hamiltonians are not hoimogeneous in moments or coordinates. This can be seen on the Painlev\'e II example - in the fully reduced coordinates the Hamiltonian writes as 
$$
H = \frac{p^2}{2}-\frac{1}{2}\left(q^2+\frac{t}{2}\right)-\theta q,
$$
while the action along the solution writes as
$$
\d \mathcal S = p\d q - H\d t = \left[p \frac{\partial q}{\partial t} - H \right ]\d t = \left[p^2 - H\right] \d t = \left[\frac{p^2}{2}+ \frac{1}{2}\left(q^2+\frac{t}{2}\right)+\theta q\right]\d t \neq H\d t.
$$
The classical action now differs from the $\tau$-function by some function that depends on time. This deviation from the classical action functional was investigated in the paper by Its and Prokhorov \cite{ItsPro} for the classical Painlev\'e equations in fully reduced coordinates. From the quantum point of view, the reduction is a restriction to the eigenspace of the Casimir operator which  partially provides separation of variables in the quantum problem. By passing to a smaller number of coordinates, the parts which were depending on the lifted coordinates vanish, so the structure of the solution changes rapidly. However, despite the fact that the theorem doesn't work in the reduced case, we still see the avatars of this statement since the $\tau$-function still enters the quasiclassical solution in some way, see paper \cite{ItsPro} and formula (2.27) in \cite{ZabrodinZotov}.

\end{document}